\documentclass[final,5p,twoside]{elsarticle}

\usepackage{epsfig}

\usepackage{amssymb}

\usepackage{natbib}
\usepackage{rotating}
\usepackage{setspace}
\usepackage{multirow}
\usepackage{graphicx}
\usepackage{subfig}
\usepackage{float}
\usepackage{amssymb}
\usepackage{mathtools}
\usepackage{pdflscape}
\usepackage{booktabs}
\usepackage{adjustbox}
\usepackage{afterpage}
\usepackage{xcolor}
\usepackage{hyperref}
\usepackage{bigstrut}
\usepackage{tikz}

\usepackage{xcolor}
\hypersetup{
    colorlinks,
    linkcolor={red!50!black},
    citecolor={blue!50!black},
    urlcolor={blue!80!black}
}

  \definecolor{darkgreen}{RGB}{67,115,41}  
  \definecolor{mycolor}{RGB}{0,160,160}  
  \definecolor{mycolor2}{RGB}{0,114,189}
  \definecolor{mycolor3}{RGB}{249,0,255}

\newcommand*{\intx}{\int_{x_{i-\f{1}{2}}}^{x_{i+\f{1}{2}}}}

\newcommand*{\f}[2]{\frac{#1}{#2}}
\newcommand*{\df}[2]{\dfrac{#1}{#2}}

\newcommand*\colvec[3][]{
    \begin{bmatrix}\ifx\relax#1\relax\else#1\\\fi#2\\#3\end{bmatrix}
}

\newcommand*\colvecfour[4]{
    \begin{bmatrix}
    #1 \\
    #2 \\
    #3 \\
    #4
    \end{bmatrix}
}
\newcommand*{\one}{excited}
\newcommand*{\two}{activation}
\newcommand*{\dt}{\partial_t}
\newcommand*{\dttotal}{\frac{d}{dt}}

\newcommand*{\dx}{\partial_x}

\newcommand*{\Q}{\mathbf{Q}}
\newcommand*{\A}{\mathbf{A}}
\newcommand*{\F}{\mathbf{F}}

\renewcommand*{\S}{\mathbf{S}}

\newcommand*{\R}{\mathbf{R}}

\begin{document}

\begin{frontmatter}

\title{A one-dimensional mathematical model of collecting lymphatics coupled with an electro-fluid-mechanical contraction model and valve dynamics}

\author[Christian]{Christian Contarino\corref{cor1}}
\ead{christian.contarino@unitn.it}
\cortext[cor1]{Corresponding author}
\author[Toro]{Eleuterio F. Toro}
\ead{eleuterio.toro@unitn.it}
\address[Christian]{Department of Mathematics, University of Trento}
\address[Toro]{Laboratory of Applied Mathematics, DICAM, University of Trento}

\begin{abstract}
We propose a one-dimensional model for collecting lymphatics coupled with a novel Electro-Fluid-Mechanical Contraction (EFMC) model for dynamical contractions, based on a modified FitzHugh-Nagumo model for action potentials. The one-dimensional model for a compliant lymphatic vessel is a set of hyperbolic Partial Differential Equations (PDEs). The EFMC model combines the electrical activity of lymphangions (action potentials) with fluid-mechanical feedback (stretch of the lymphatic wall and wall shear stress) and the mechanical variation of the lymphatic wall properties (contractions). The EFMC model is governed by four Ordinary Differential Equations (ODEs) and phenomenologically relies on: (1) environmental calcium influx, (2) stretch-activated calcium influx, and (3) contraction inhibitions induced by wall shear stresses. We carried out a complete mathematical analysis of the stability of the stationary state of the EFMC model. Contractions turn out to be triggered by the instability of the stationary state.
Overall, the EFMC model allows imitating the influence of pressure and wall shear stress on the frequency of contractions observed experimentally. Lymphatic valves are modelled using a well-established, versatile lumped-parameter model which allows simulating stenotic and regurgitant valves and quantifying their lymphodynamical effects. Modern numerical methods are employed for the one-dimensional model (PDEs), for the EFMC model and valve dynamics (ODEs).
Adopting the geometrical structure of collecting lymphatics from rat mesentery, we show numerical tests inspired by experiments and also an example of a Riemann problem for lymphatics. 
We analysed several lymphodynamical indexes of a single lymphangion for a wide range of upstream and downstream pressure combinations where positive and negative pressure gradients are taken into account.
The most influential model parameters were found by performing two sensitivity analyses for positive and negative pressure gradients.
Stenotic and regurgitant valves were modelled, and their effects are here quantified.
Results for stenotic (or obstructed) valves showed in the downstream lymphangion that for low frequencies of contractions the Calculated Pump Flow (CPF) index remained almost unaltered, while for high frequencies the CPF dramatically decreased depending on the severity of the stenosis (up to 93$\%$ for a severe stenosis).
Results for incompetent (or regurgitant) valves showed that the net flow during a lymphatic cycle tends to zero as the degree of incompetence increases, and this is representative for a defective valve unable to prevent backflows.
\end{abstract}


\end{frontmatter}

\tableofcontents


\section{Introduction}

The lymphatic system is an intricate network of vessels and nodes which connect tissues to the bloodstream. The main functions of the lymphatic system comprise the hydrodynamic homoeostasis of the fluid tissues through the drainage of the excess interstitial fluid, the transport of proteins and waste products, as well as the transport of immune cells. From the pumping action of the heart, blood flows through arteries, arrives at capillaries, and is filtered into interstitial fluids through Starling forces. A portion of it, known as lymph, enters into the lymphatic system, is pushed uphill through several mechanisms, is filtered at lymph nodes, and is drained into the venous system at the junction between the internal jugular veins and the subclavian veins.
In addition, according to the classical theory, cerebrospinal fluid drains into nasal lymphatics through the cribriform plate, flows into deep cervical lymph nodes, and enters into the subclavian veins. Moreover, interstitial fluid from the central nervous system parenchyma drains into lymph nodes within the basement membrane of cerebral capillaries and arteries and then enters into the subclavian veins \cite{Bakker:2016a, Engelhardt:2016a,Nedergaard:2013a, Brinker:2014a, Cserr:1992a}.
From there on, lymph mixes with blood again, enters into the heart and the cycle restarts. 
During these regular cycles, leukocytes travel throughout the lymphatic system and provide an immune response for the body. 

Several pathologies are associated with the impairment of the lymphatic system.  Very recently, Louveau et al. \cite{Louveau:2015aa} and Aspelund et al. \cite{Aspelund:2015a} have just discovered the lymphatic system in the brain, known as meningeal or dural lymphatics. This discovery has significant effects on the immune privilege concept of brain, on brain-waste clearance of parenchyma, on the genesis of neurological disorders \cite{Louveau:2015aa,Dissing_Olesen:2015a,Conway:2015a,Zamboni:2015b,Raper:2016a,Engelhardt:2016a}, and on possibly new therapeutic treatments for removal of macromolecules in the central nervous system \cite{Louveau:2016a}.

The building block of the lymphatic system is the lymphangion: a mini-heart like, compliant vessel, which contracts and pushes lymph centripetally, and has several mechanobiological auto-regulatory systems to provide optimal flows in various scenarios \cite{Munn:2015a,Kunert:2015a}. It is enclosed between valves which prevent reverse flow and allow for unidirectional flows. The contraction of lymphangions has been described for several animals in \cite{McHale:1992a,Benoit:1989a,Gashev:2008a}; it is regulated by Ca$^{2+}$ fluxes from extracellular and intracellular stores \cite{Weid:2001a}, is divided in systolic and diastolic phases \cite{Benoit:1989a}, and can be described using cardiac indexes \cite{Davis:2011a,Davis:2012a,Scallan:2012a,Scallan:2013a}. Several works have been devoted to studying the membrane potential in lymphangions \cite{Telinius:2015a,Weid:1997a,Ohhashi:1978a}.

The mechanobiological systems which underlie the auto-regulatory homeostatic functions of lymphatic contractions, are still being investigated \cite{Breslin:2014a,Margaris:2012a,Munn:2015a}. 
The frequency of the lymphatic contractions depends on the rate of stretch of the vessel wall. Indeed, stretch-activated calcium channels in the vessel wall change conformation in response to membrane tension \cite{Munn:2015a,Telinius:2015a}. McHale and Roddie \cite{McHale:1976a} described the dynamics of frequencies in bovine mesenteric vessels varying intraluminal pressures and showed that the more the lymphangions are stretched, the faster they contract. 
Lymphatic muscle also exhibits rate-sensitive contractile responses to stretch, as described for rat mesenteric lymphatics by Davis et al. \cite{Davis:2009a}.
The authors analysed the responses in amplitude and frequency to time-varying preload and pressure. Bursts of contraction occurred when positive ramps were imposed. 
Lymphatic contractions are also inhibited by flows. Gashev et al. \cite{Gashev:2002a} studied rat mesenteric lymphatics in response to imposed flows and showed that the frequency dropped from $9.0\pm 1.6$ min$^{-1}$ to $3.1\pm 1.4$ min$^{-1}$ when flow changed from zero to a transaxial-pressure-gradient induced flow of $7$ cmH$_2$O. In addition, the activity of the lymphatic contractions differs from region to region. 
Gashev et al. \cite{Gashev:2004a} showed that the flow-induced inhibition of contraction in the thoracic duct of rats is more evident when compared to that of femoral lymphatics. As a matter of fact, the authors showed that the frequencies of thoracic duct and femoral lymphatics are $4.6\pm 0.6$ min$^{-1}$ and $15.2\pm2.6$ min$^{-1}$, respectively, at zero flow conditions, and $0.1\pm 0.1$ min$^{-1}$ and $5.0 \pm 2.4$ min$^{-1}$, respectively, at transaxial-pressure-gradient induced flow of $5$ cmH$_2$O. An underlying mechanism of the above-mentioned flow-induced contraction inhibition is the local production of NO induced by wall shear forces \cite{Kunert:2015a,Dixon:2006a,Kornuta:2015a}, which results in blunting the Ca$^{2+}$-dependent contraction. A lymphatic vessel composed of two or more lymphangions is called a collecting lymphatic. Several studies have analysed the coordination of adjacent lymphangions \cite{Crowe:1997a,McHale:1992a,Zawieja:1993a}, wave-propagations throughout the collector \cite{Ohhashi:1980a,Zawieja:1993a,Akl:2011a} and the microstructure of the vessel wall \cite{Ohhashi:1980a,Rahbar:2012a,Caulk:2015a}. Overall, the above-mentioned processes which underlie and permit optimal lymph flows, are more complex and difficult to quantify, when compared to those of arteries and veins. For complete reviews of the mechanics of lymphangions and collectors, see \cite{Munn:2015a,Breslin:2014a,Margaris:2012a,Nipper:2011a,Weid:2004a,Bridenbaugh:2003a,Zawieja:2011a}.

The lymphatic system has two different types of valves called primary and secondary valves. The former is located at the initial lymphatics at the level of the endothelium, while the latter is located between lymphangions in collectors \cite{Schmid:2003a, Bazigou:2012a}.
Primary and secondary lymphoedema, a lymphatic disease that leads to tissue swelling, is linked to lymphatic valve deficits \cite{Kinmonth:1954a,Mellor:2011a,Rockson:2008a,Noel:2001a,Mihara:2012a}. For instance, the lack of valves in the lymphoedema distichiasis impairs lymphatic flow due to the inability to properly pump lymph forward \cite{Mellor:2011a, Petrova:2004a, Sabine:2015a,Bazigou:2012a}. Also, the chronic venous insufficiency leads to fibrotic lymph vessels due to hypertension, then it compromises the functionality of lymphatic valves, and finally result in accumulation of fluid in tissues \cite{Mortimer:2004a, Rasmussen:2016a}. In addition, in secondary lymphoedema after lymph node dissection, lymph retention and lymphatic hypertension occur, and valvular dysfunction induces retrograde lymph flow \cite{Mihara:2012a}. Despite the connections of lymphatic valve deficits and lymphoedema, to the authors' knowledge, the effect of stenotic and regurgitant valves in the lymphatic system has not properly been investigated and quantified. This is probably due to the difficulties in performing experiments on animal lymphatic valves, though the effects of genes mutations in engineered mice can be studied. 

A huge gap is currently present in the literature between mathematical models for the circulatory \cite{Sherwin:2003b,Matthys:2007aa,Liang:2009b,Alastruey:2011a,Liang:2014a,Mueller:2014a,Mueller:2014b,Blanco:2014a,Mynard:2015a,Causin:2015a,Quarteroni:2016a,Vergara:2016a}, and lymphatic systems. The first mathematical attempt to model for the lymphatic system is attributed to Reddy et al. \cite{Reddy:1974a}, and was based on one-dimensional flow equations. The model was then refined by MacDonald et al. \cite{Macdonald:2008a} including a tension and a damping term. However, the convection term, tube law and valve model were relatively simple. Extensive work has been done in lumped-parameter models \cite{Drake:1986a,Venugopal:2007a,Quick:2006a,Quick:2008a,Venugopal:2010a,Bertram:2011a, Bertram:2013a,Jamalian:2013a,Bertram:2014a,Gajani:2015a,Bertram:2016a,Jamalian:2016a,Caulk:2016a}. Jamalian et al. \cite{Jamalian:2016a} constructed a lumped-parameter model able to simulate lymph transport in a network of rat lymphangions based on experimental measurements for vessel and valve parameters \cite{Bertram:2013a}. The authors analysed the effect of pumping coordination in branched network structures. Caulk et al. \cite{Caulk:2015a} described a detailed microstructurally four-fibre family constitutive law for rat thoracic duct based on experimental measurements of collagen, elastin and solid structures of the lymphatic wall. Then, the authors combined the lumped-parameter model described by Bertram et al. \cite{Bertram:2014a} with their four-fibre family constitutive law, and studied the variation of muscle contractility in response to a sustained elevation in afterload \cite{Caulk:2016a}. Rahbar et al. \cite{Rahbar:2011a} investigated the validity of assuming Poiseuille flow to estimate the wall shear stress. A mechanobiological oscillatory model for the lymphatic contraction has been proposed by Kunert et al. \cite{Kunert:2015a}. Compared to the other existing mathematical models for lymphangions, the authors provided a biologically-based, dynamical model for the contractibility of the vessel wall. They proposed 1) evolutionary equations for Ca$^{2+}$ and NO, and 2) a coupling between wall mechanical forces, Ca$^{2+}$ and NO. The resulting model was able to control lymphatic transport via mechanobiological feedback loops, given by stretch-activated contractions and flow-induced relaxations. Very recently, Baish et al. \cite{Baish:2016a} proposed a model of a vascular oscillator and studied the interaction of Ca$^{2+}$ and NO in the context of lymphatic vessel contractions. 

With the aim of constructing a mathematical model of the entire lymphatic system, the above-mentioned mathematical models for lymphangions, except for \cite{Kunert:2015a, Baish:2016a}, are based on a relatively simple contraction dynamic. As a matter of fact, these models 1) prescribe muscle contractility dynamics by using trigonometric functions, and 2) artificially prescribe time delays between adjacent lymphangions. This result in non-dynamical models insofar as contractions are triggered by artificial parameters rather than by local dynamical factors, and contraction frequencies in those models do not depend on the rate of stretch of the lymphatic wall, nor on local shear forces.
 

In the study of cardiac contraction, there exist an extensive literature. From the pioneer, detailed model of Hodgkin and Huxley in 1952 \cite{Hodgkin:1952a} of how action potential in neurons are initiated and propagated, several extended and even simplified models have been proposed in the literature for neurons and heart contractions \cite{Aliev:1996a,Nagumo:1962a,Panfilov:1993a,Kogan:1991a}. The FitzHugh-Nagumo model \cite{Nagumo:1962a} is an example of a significantly simplified, but elegant two-parameter formulation of the original Hodgkin-Huxley model \cite{Goktepe:2009a}. It is a set of two Ordinary Differential Equations (ODEs) with a fast and a slow variable. The former represents the action potential, while the latter phenomenologically summarises all the effects of all ionic currents. The FitzHugh-Nagumo model has been adapted to simulate rhythmically discharging cells in the heart membrane, such as the atrioventricular and sinoatrial node \cite{Colli_Franzone:2014a}.
Thanks to its simplicity, it also permits a rigorous analysis of the stability. Action potential arises to be triggered by an unstable equilibrium state of the ODEs. Many studies have been done to couple modified versions of the FitzHugh-Nagumo model to the heart contractions \cite{Colli_Franzone:2014a}. However, to date no studies have been performed to model contractions of lymphangions with the over-mentioned dynamical and phenomenological set of ODEs for action potentials. 

In the present paper, we propose a one-dimensional model for collecting lymphatics coupled with an Electro-Fluid-Mechanical Contraction (EFMC) model for dynamical lymphatic contractions based on a modified FitzHugh-Nagumo model. 
The one-dimensional model for a compliant lymphatic vessel adopt a tube law which fits the experimental measurements shown in Bertram et al. \cite{Bertram:2013a} and performed by Davis et al. \cite{Davis:2011a}. The resulting system is a set of hyperbolic Partial Differential Equations (PDEs), is written in conservative formulation, allows for space-variable geometrical parameters and permits simulating contraction-wave propagations based on a prescribed space-time variation of the Young modulus. 
The EFMC model is a set of four ODEs. We carried out a complete mathematical analysis of the EFMC model on the stability of its equilibrium state. Contractions turn out to be triggered by the instability of the stationary state and phenomenologically depend on: (1) environmental calcium influx, (2) stretch-activated calcium influx, and (3) contraction inhibitions induced by wall shear stresses. 
Lymphatic valves are modelled using a well-established, versatile lumped-parameter model proposed by Mynard et al. \cite{Mynard:2012a} and analysed by Toro et al. \cite{Toro:2015aa}. This valve model also allows simulating stenotic and regurgitant valves and quantify the lymphodynamical effects. 
Here, we simulate lymphatic vessel from the mesentery of rats. We use the SLIC method \cite{Toro:2000a} to numerically solve the one-dimensional equations and an implicit Runge-Kutta solver for the systems of ODEs.
We show some numerical tests of our mathematical model inspired by experiments \cite{Davis:2011a, Davis:2012a, Scallan:2012a, Scallan:2013a}. 
Then we analyse several lymphodynamical indexes of a single lymphangion for different combinations of upstream and downstream pressures and considering both positive and negative pressure gradients. We performed two sensitivity analyses based on \cite{Griensven:2006a,Liang:2014a} for both positive and negative pressure gradients and investigate the most influential parameters affecting the lymphodynamical indexes. Finally, we quantified the effects of stenotic and regurgitant valves.

\begin{figure*}[t]
\centering
\includegraphics[width=1\textwidth]{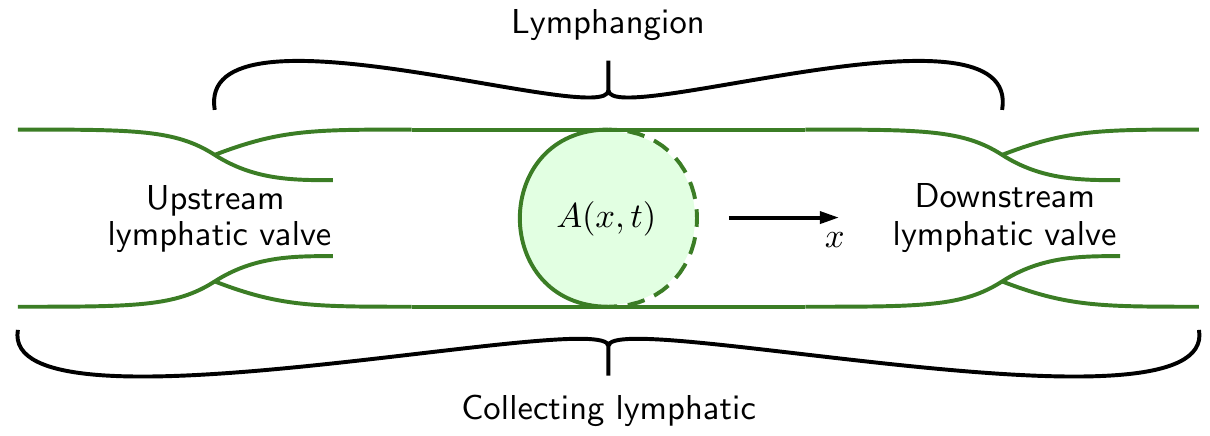} 
\caption{\scriptsize {\bf Illustration of a collecting lymphatic.} A lymphangion is a lymphatic vessel delimitated by upstream and downstream valves and a collecting lymphatic is a chain of lymphangions. The natural lymph flow direction is from the upstream side to the downstream one and the lymphatic valves prevent backflows throughout the collecting lymphatic.  The figure also shows a general cross-sectional area. } \label{fig:IllustrationLymphangion}
\end{figure*}

The rest of this paper is structured as follows: in Section \ref{sec:mathemodel} we propose a one-dimensional model for lymph flow that allows for variable parameters, we describe the EFMC model, analyse it and then we briefly review the valve model. In Section \ref{sec:NumericalMethods} we briefly review the finite volume schemes used for the one-dimensional lymph flow equations, explain how to couple valves with lymphangions, how to impose boundary pressure at the interfaces of the lymphangions, and describe the numerical methods used for the valve and EFMC models. In Section \ref{sec:results} we show some numerical results our mathematical model, analyse the lymphodynamical indexes for a wide range of boundary pressures, perform two sensitivity analyses for positive and negative pressure gradients and finally investigate the effects of defective lymphatic valves, namely stenotic and regurgitant valves. Section \ref{Sec:conclusions} gives a summary and concluding remarks.

{\footnotesize
\section{Mathematical models} \label{sec:mathemodel}

A {\em lymphangion} is a lymphatic vessel delimited by two valves, also called secondary lymphatic valves, located on the upstream and downstream boundaries of the lymphangion. A chain of lymphangions is called {\em collecting lymphatic}. Here we aim to model the dynamics of flowing lymph inside a collecting lymphatic propelled by lymphatic contractions and pressure gradients, and the dynamic of lymphatic valves. First, we present a one-dimensional mathematical model for lymph flow and a set of four ODEs to model the lymphatic contraction, and then we describe a mathematical model for lymphatic valves.
Fig. \ref{fig:IllustrationLymphangion} illustrates a collecting lymphatic, a single lymphangion and two lymphatic valves. It also shows a generic cross-sectional area at a given position $x$ and time $t$. 

\subsection{A one-dimensional model for lymph flow} 
Here we assume the lymph to be an incompressible Newtonian fluid. To derive the one-dimensional flow equation for a compliant lymphatic vessel, one can follow the procedure done for arteries and veins, where Reynold's transport theorem is used to obtain the conservation of mass and momentum in a flexible tube, see \cite{Toro:2013a,Formaggia:2009a,Toro:2016b}.
The one-dimensional flow equations for a compliant lymphatic vessel are 
\begin{equation}\label{eq:equazioniSangue}
\begin{cases}
\displaystyle
\dt A + \dx q=0\;, \\
\dt q + \dx \left(\alpha\df{q^2}{A}\right)+\df{A}{\rho}\dx p=-\df{f}{\rho}\;,
\end{cases}
\end{equation}
where $x$ is the space variable, $t$ is time, $\alpha$ is the Coriolis coefficient assumed to be $\alpha=1$, $A(x,t)$ is cross-sectional area of the vessel, $q(x,t)=A(x,t)u(x,t)$ is flow, $u(x,t)$ is velocity, $p(x,t)$ is pressure, $\rho$ is lymph density, $f(x,t)=2(\gamma+2) \pi \mu u(x,t)$ is friction force per unit length of the tube with parameter $\gamma$ chosen from the velocity profile \cite{Alastruey:2006a}, and $\mu$ is the dynamic viscosity. There are two governing partial differential equations and three unknowns, namely $A(x,t)$, $q(x,t)$ and $p(x,t)$. For this reason, an extra relation is required to close the system, the {\em tube law}, which relates pressure $p(x,t)$ and cross-sectional area $A(x,t)$. Equations \eqref{eq:equazioniSangue} have been widely described in the literature for blood flow \cite{Toro:2016b,Toro:2012a,Toro:2013a,Formaggia:2009a}, where different tube laws were proposed depending on the characteristics of the wall material. The lymphatic wall is characterized by elastin, collagen, smooth muscle cells and other extracellular matrix constituents \cite{Caulk:2015a,Rahbar:2012a}. Elastin fibres give to lymphatic vessels a compliant, elastic behaviour, while collagen prevents vessels from stretching beyond their physiological limits. The overall dynamics of elastin and collagen is reflected in highly non-linear tube laws. Several works have been proposed to describe tube laws for rats \cite{Caulk:2015a,Rahbar:2012a,Bertram:2014a}, bovine \cite{Macdonald:2008a} and human thoracic ducts \cite{Telinius:2010a}. Here we propose the following general and purely elastic tube law:
\begin{equation}\label{eq:pressure}
p(x,t)=K(x,t)\psi \left(A(x,t);A_0(x)\right)+p_e(x,t)\;,
\end{equation}
with 
\begin{equation} \label{eq:TubeLaw}
\begin{aligned}
\psi (A(x,t);A_0(x))&=\left(\df{A(x,t)}{A_0(x)}\right)^m-\left(\f{A(x,t)}{A_0(x)}\right)^n \\ & +C\left[ \left(\df{A(x,t)}{A_0(x)}\right)^z-1\right]\;,
\end{aligned}
\end{equation}
where $p_e(x,t)$ is the external pressure, $A_0(x)$ is vessel cross-sectional area at equilibrium, $K(x,t)$ is the stiffness of the vessel wall,  $m\geq 0$, $n\leq 0$, $z\geq 0$, and $C\geq 0$ are real numbers to be specified. 
The {\em transmural pressure} is defined as
\begin{equation}\label{eq:transmuralpressure}
p_{transm}(x,t):=p(x,t)-p_e(x,t)\;.
\end{equation}
To the authors' knowledge, a mechanical study of the stiffness of collecting lymphatic is not present in the literature. It is clear that lymphatics are highly deformable and their physiological pressure is comparable or smaller than that of veins. As a matter of fact, the lymphatic pressure approximately ranges from the interstitial fluid pressure to the subclavian vein pressure.
Here we assume parameter $K$ of lymphatics to be equal to that of veins \cite{Mueller:2014a}, namely 
\begin{align}
K(x,t)=\df{E(x,t)}{12(1-\nu^2)}\Bigg(\df{h_0(x)}{r_0(x)}\Bigg)^3\;,
\end{align}
where $E(x,t)$, $\nu$, $h_0(x)$ and $r_0(x)$ are respectively the Young modulus of elasticity, the Poisson ratio, the wall-thickness and cross-sectional radius at zero transmural pressure (equilibrium), respectively. The assumption to use parameter $K$ equal to that of veins, might need further investigations. 
Following \cite{Macdonald:2008a,Macdonald:2008b}, lymphatic contractions are modelled by varying the Young modulus $E(x,t)$ from a minimal value $E_{min}(x)$ to a maximum value $E_{max}(x)$ as follows
\begin{equation}\label{eq:ContractionYoungModulus}
\begin{aligned}
E(x,t)=&E_{min}(x)+s(x,t)\big(E_{max}(x)-E_{min}(x)\big)\;, \\& s(x,t)\in [0,1]\;,
\end{aligned}
\end{equation}
where $s(x,t)$ is the {\em state of contraction}. The lymphangion can be relaxed and contracted and this is described by variable $s(x,t)$. The lymphangion is contracted when $s(x,t)=1$ and the Young modulus reaches its maximum $E(x,t)=E_{max}$ while the lymphangion is relaxed when $s(x,t)=0$ and the Young modulus reaches its minimum $E(x,t)=E_{min}$. In Section \ref{sec:FitzHugh} we provide the governing equation for $s$ with a set of four ODEs. 

In the present paper, we simulated lymphatic vessel from the mesentery of rats, whose parameters are found in Table \ref{table:parameters}. The tube law used in Eq. \eqref{eq:TubeLaw} is a relationship between normalised cross-sectional area and pressure. The parameters of the tube law and the minimum and maximum Young's modulus were tuned to fit the experimental measurements shown in Bertram et al. \cite{Bertram:2013a} and performed by Davis et al. \cite{Davis:2011a}. The experimental measurements and the tube law with and without contraction are shown in Fig. \ref{fig:tubeLawfitting}. 

\begin{figure}[t]
\begin{center}
\includegraphics[width=0.45\textwidth]{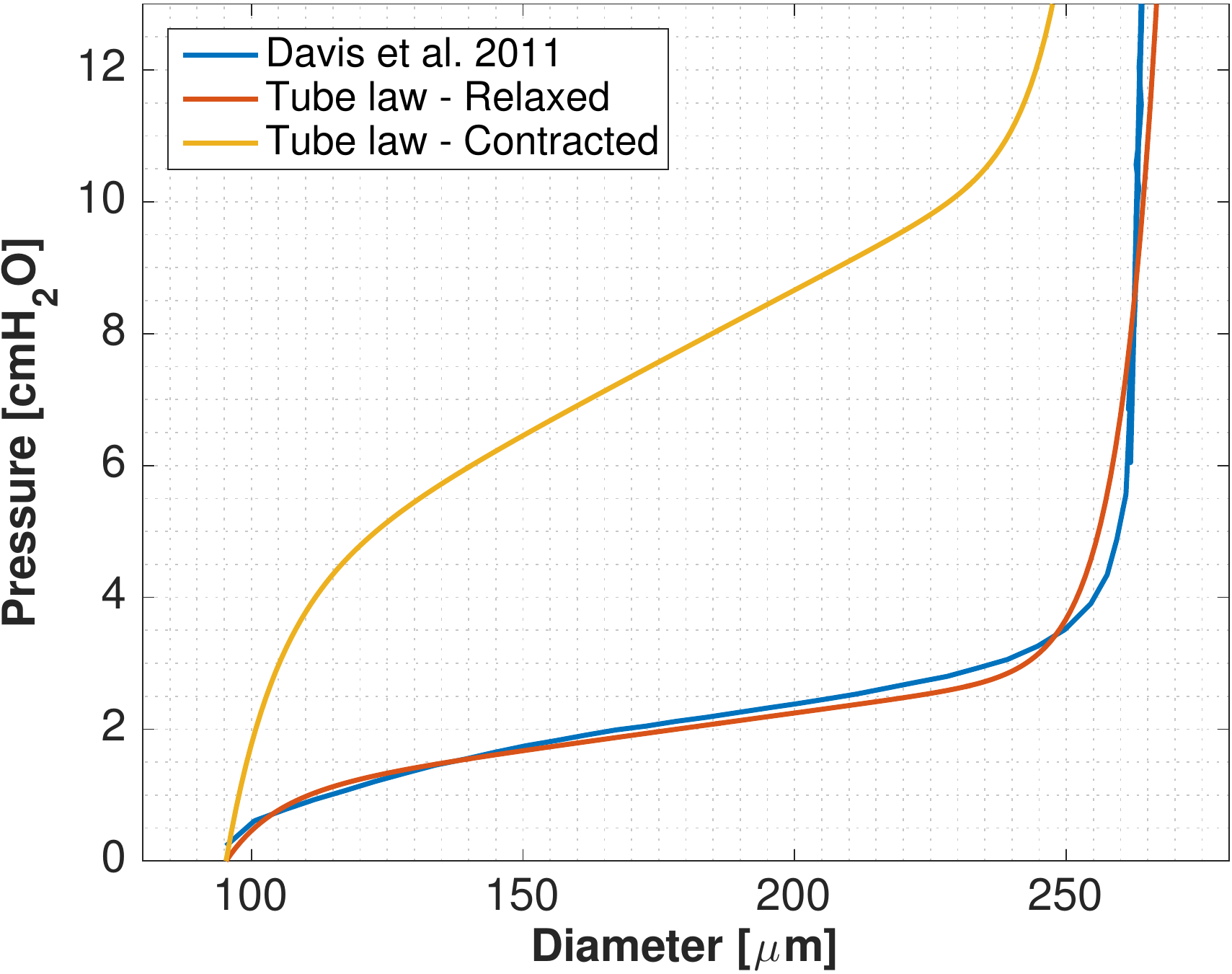}
\end{center}
\caption{\scriptsize {\bf Pressure-diameter relation}. Here we show the tube law used for the lymphatic wall. The parameters were tuned to fit the experimental measurements of Davis et al. \cite{Davis:2011a} and are found in Table \ref{table:parameters}. The figure also shows the tube law without contraction and with contraction.} \label{fig:tubeLawfitting}
\end{figure}

The Wall Shear Stress (WSS) is fundamental in connecting the fluid mechanics and the Nitric Oxide (NO) production. The fluid movement increases the shear at the lymphatic wall and increases the production of NO through the activation of local endothelium NOS. The production of NO results in vasodilating the lymphatic vessel and in decreasing the contraction frequency by blunting the $Ca^{2+}$-dependent contraction. For references, see \cite{Kunert:2015a, Dixon:2006a, Kornuta:2015a, Gashev:2002a, Breslin:2014a}. According to \cite{Alastruey:2006a}, the WSS in our formulation is 
\begin{equation}\label{eq:WSS}
\tau (x,t) = -u(x,t)\mu \f{\gamma+2}{r(x,t)}\;,
\end{equation}
where $r(x,t)$ is cross-sectional radius. 

\begin{table*}[!t] 
\scriptsize
\centering 
\begin{tabular}{c l c c} 
\toprule[1.5pt] 
{\bf Parameter} & {\bf Description} & {\bf Value} & {\bf Units} \\ 
\hline & \multicolumn{3}{ l }{{\bf Unknowns} } \\ \hline 
$A$ & Lymphatic cross-sectional area &$A(x,t)$ & $\mu$m$^2$ \\ 
$q$ & Lymph flow &$q(x,t)$ & $\mu$L min$^{-1}$ \\ 
$v$ & Excitable variable &$v(t)$ & - \\ 
$w$ & Recovery variable &$w(t)$ & - \\ 
$s$ & State of contraction ($0 \le s \le 1$) &$s(t)$ & - \\ 
$I$ & Stimulus &$I(t)$ & - \\ 
$\xi$ & State of the lymphatic valve ($0 \le \xi \le 1$) &$\xi(t)$ & - \\ 
$q_v$ & Flow across the lymphatic valve &$q_v(t)$ & $\mu$L min$^{-1}$ \\ 
\hline & \multicolumn{3}{ l }{{\bf Material parameters} } \\ \hline 
$\gamma$ & Parameter for velocity profile &2 & $-$ \\ 
$\mu$ & Lymph dynamic viscosity &1 & cP \\ 
$\rho$ & Lymph density &998 & kg m$^{-3}$ \\ 
$\nu$ & Poisson ratio &0.5 & - \\ 
$E_{max}$ & Maximum Young modulus &135000 & Pa \\ 
$E_{min}$ & Minimum Young modulus &35000 & Pa \\ 
$r_0$ & Radius at zero transmural pressure &47.7 & $\mu$m \\ 
$h_0$ & Wall-thickness at zero transmural pressure &$0.3r_0$ & $\mu$m \\ 
$A_0$ & Cross-sectional area at zero transmural pressure &$\pi r_0^2$ & $\mu$m$^2$ \\ 
$L$ & Lymphangion length &1.5 & mm \\ 
$p_e$ & External pressure &0 & cmH$_2$O \\ 
\hline & \multicolumn{3}{ l }{{\bf Tube law} } \\ \hline 
$m$ & Parameter &0.5 & - \\ 
$n$ & Parameter &-5.0 & - \\ 
$z$ & Parameter &19.0 & - \\ 
$C$ & Parameter &1.0e-16 & - \\ 
\hline & \multicolumn{3}{ l }{{\bf Electro-Fluid-Mechanical Contraction (EFMC) model} } \\ \hline 
$t_{\one}$ & Required time to perform an action potential &$\approx$ 2 $s$ & s \\ 
$t_{\two}$ & Time required to activate an action potential &Eq. \eqref{eq:ActivationTime} & s \\ 
$f_{min}$ & Minimum frequency without stretch &3.0 & min$^{-1}$ \\ 
$f_{Ca}$ & Maximum frequency at $\tilde{A}=A_{Ca}$ &20.0 & min$^{-1}$ \\ 
$R_I$ & Radius of the activation region &0.1 & - \\ 
$n_{Ca}$ & Stretch-activation parameter &10 & - \\ 
$A_{Ca}$ & Reference stretch-activation cross-sectional area &7.75$A_0$ & $\mu$m$^2$ \\ 
$k_{Ca}^{(1)}$ & Baseline increasing rate of stimulus $I$ &Eq. \eqref{eq:parameterkCa} & s$^{-1}$ \\ 
$k_{Ca}^{(2)}$ & Stretch-activated increasing rate of stimulus $I$ &Eq. \eqref{eq:parameterkCa} & s$^{-1}$ \\ 
$k_{rel}$ & Decreasing rate of the stimulus $I$ &10 & s$^{-1}$ \\ 
$a_1$ & Parameter &100 & s$^{-1}$ \\ 
$a_2$ & Parameter &0.5 & - \\ 
$a_3$ & Parameter &25.0 & - \\ 
$b_1$ & Parameter &3.0 & s$^{-1}$ \\ 
$b_2$ & Parameter &0.0 & s$^{-1}$ \\ 
$c_1$ & Increasing rate of contraction state $s$ &110 & s$^{-1}$ \\ 
$c_2$ & Decreasing rate of contraction state $s$ &3 & s$^{-1}$ \\ 
$\tilde{I}$ & Approximated stimulus required to trigger an action potential &Eq. \eqref{eq:Imean} & - \\ 
$k_{NO}$ & Contraction inhibition parameter ($0 \le k_{NO} \le 1$) &0.8 & - \\ 
$\tau_{NO}$ & Reference wall shear stress &6.0 & dyne cm$^{-2}$ \\ 
$n_{NO}$ & Wall shear stress inhibition parameter &1.2 & - \\ 
\hline & \multicolumn{3}{ l }{{\bf Valve model} } \\ \hline 
$\Delta p_{open}$ & Valve opening threshold pressure difference &0 & cmH$_2$O \\ 
$\Delta p_{close}$ & Valve closure threshold pressure difference &0 & cmH$_2$O \\ 
$K_{vo}$ & Rate coefficient valve opening &10.0 & Pa$^{-1}$ s$^{-1}$ \\ 
$K_{vo}$ & Rate coefficient valve closure &10.0 & Pa$^{-1}$ s$^{-1}$ \\ 
$B$ & Bernoulli resistance &Eq. \eqref{eq:ValveParameters} & cmH$_2$O s$^2$ $\mu$L$^{-2}$ \\ 
$L$ & Lymphatic inertia &Eq. \eqref{eq:ValveParameters} & cmH$_2$O s$^2$ $\mu$L$^{-1}$ \\ 
$R$ & Viscous resistance to flow &Eq. \eqref{eq:ValveParameters} & cmH$_2$O s $\mu$L$^{-1}$ \\ 
$M_{st}$ & Maximum valve opening ($0 \le M_{st} \le 1$) &1.0 & - \\ 
$M_{rg}$ & Minimum valve closure ($0 \le M_{rg} \le 1$) &0.0 & - \\ 
$L_{eff}$ & Effective length &0.1 & mm \\
\bottomrule[1.5pt] 
\end{tabular} \caption{\scriptsize {\bf Parameters used for the one-dimensional EFMC model for lymph flow.} We adopted the geometrical structure of collecting lymphatics from rat mesentery.} \label{table:parameters} 
\end{table*} 

\subsubsection{Conservative formulation of the one-dimensional lymph flow equations} 
It is possible to write the lymph flow equations in conservative form as follows:
\begin{equation}\label{eq:Lymph}
\dt\Q(x,t) +\dx \F(\Q(x,t),x,t)=\S(\Q(x,t),x,t)\;,
\end{equation}
where 
\begin{align}
\Q\left(x,t\right)&=\colvec{A\left(x,t\right)}{A\left(x,t\right)u\left(x,t\right)}\;, \label{eq:ConservedQuantity} \\
\F\left(\Q,x,t\right)&=\colvec{Au}{Au^2-\df{K}{\rho}A_0\partial_{A_0}\Psi}\;, \label{eq:PhysicalFlux} \\
\S\left(\Q,x,t\right)&=\colvec{0}{-\df{1}{\rho}\bigg(f+A\dx p_e +\Psi \dx K+K\dx A_0 \partial_{A_0}\Psi\bigg)}\;, 
\label{eq:SOURCE}
\end{align}
with
\begin{equation}
\begin{aligned}
\Psi&=\int_A\psi(A;A_0) \mathrm{d}A\\ &
\begin{aligned}=A_0\Bigg(&\df{1}{m+1}\bigg(\df{A}{A_0}\bigg)^{m+1}-\df{1}{n+1}\bigg(\df{A}{A_0}\bigg)^{n+1}\\ +C&\f{1}{z+1}\bigg(\df{A}{A_0}\bigg)^{z+1}\Bigg)\;,\label{eq:intTransm}
\end{aligned}
\end{aligned}
\end{equation}
and
\begin{equation}
\begin{aligned}
\partial_{A_0}\Psi&=\partial_{A_0}\Psi(A;A_0)\\ &
\begin{aligned}=-\Bigg(&\df{m}{m+1}\bigg(\df{A}{A_0}\bigg)^{m+1}-\df{n}{n+1}\bigg(\df{A}{A_0}\bigg)^{n+1}\\ +C&\f{z}{z+1}\bigg(\df{A}{A_0}\bigg)^{z+1}\Bigg)\;.\label{eq:intTransm2}
\end{aligned}
\end{aligned}
\end{equation}
The constants arising from the integrals \eqref{eq:intTransm} and \eqref{eq:intTransm2} are set to zero for consistency with \eqref{eq:equazioniSangue} and \eqref{eq:pressure}, see \cite{Elad:1991a, Elad:1999a, Toro:2016b}.

The present formulation allows for space-time variable Young modulus. This let us simulate travelling contraction-waves through the lymphatic wall by prescribing a space-time varying trigonometric contraction state $s(x,t)$. However, in the present work we consider the simpler case in which the contraction state is constant throughout the lymphangion, namely $s=s(t)$, and we also neglect the interaction between adjacent lymphangions. Then, instead of prescribing a trigonometric function for $s$, here we propose a set of governing ODEs given in Section \ref{sec:FitzHugh}.

Here, parameters $h_0(x)$, $r_0(x)$, $E_{min}(x)$, $E_{max}(x)$ and $p_e(x)$ are constant in space. As a result, the source term simplifies in
\begin{equation}
\S(\Q,x,t)=\S(\Q)=\colvec{0}{-2\left(\gamma+2\right)\pi \frac{\mu}{\rho}u}\;.
\end{equation}
The general case of variable material properties poses numerical and mathematical \cite{Han:2015a, Han:2015b} challenges and require the use of well-balanced schemes, see \cite{Pares:2006a, Mueller:2013c}.

\subsubsection{Mathematical analysis of the one-dimensional lymph flow equations}
Here we study the mathematical properties of \eqref{eq:Lymph} assuming constant parameters along the lymphatic vessel. The equations in \eqref{eq:Lymph} are a generalization of the one-dimensional blood flow equations presented in \cite{Toro:2013a}. As a matter of fact, the main difference is an additional term in the tube law \eqref{eq:TubeLaw}. For this reason, here we summarize the main mathematical structure of the lymph flow equations without proofs. System \eqref{eq:Lymph} can be written in quasi-linear form as 
\begin{equation}\label{eq:NonConservativeLymph}
\dt \Q +\A(\Q,t)\dx\Q=\S(\Q)\;,
\end{equation}
where
\begin{equation}
\A(\Q,t)=\begin{bmatrix}
       0 & 1 \\
       \df{A}{\rho}\partial_{A}\psi-u^2 & 2u \\
     \end{bmatrix}\;, \quad \S(\Q)=\colvec{0}{-\df{f}{\rho}}\;.
\end{equation}
The eigenvalues of matrix $\A$ are
\begin{equation}
\lambda_1=u-c\;, \quad \lambda_2=u+c\;, 
\end{equation}
where $c$ is the {\em wave speed}
\begin{equation}
c=\sqrt{\df{A}{\rho}\partial_{A}\psi}=\sqrt{\df{K}{\rho}\left[ m\left(\df{A}{A_0}\right)^m-n\left(\df{A}{A_0}\right)^n+Cz\left(\df{A}{A_0}\right)^z\right]}\;.
\end{equation}
We assume parameters $m\geq 0$, $n\leq 0$, $z\geq 0$, and $C\geq 0$ for the tube law. Thus, the wave speed $c$ is always real. 
The eigenvectors of $\A$ are
\begin{equation}
\R_1=\gamma_1\colvec{1}{u-c}\;, \quad \R_2=\gamma_2\colvec{1}{u+c}\;,
\end{equation}
where $\gamma_1$ and $\gamma_2$ are arbitrary scaling factors. It can be shown that system \eqref{eq:Lymph} is hyperbolic, as the eigenvalues are real and distinct and the eigenvectors $\R_1$ and $\R_2$ are linearly independent. 
Following proofs in \cite{Toro:2013a} and \cite{Toro:2016b}, the $\lambda_1$ and $\lambda_2$ characteristic fields are genuinely non-linear outside the locus of the following function
\begin{equation}
\begin{aligned}
G\left(\df{A}{A_0}\right)&=m\left(m+2\right)\bigg(\df{A}{A_0}\bigg)^m-n\left(n+2\right)\bigg(\df{A}{A_0}\bigg)^n \\&+Cz\left(z+2\right)\bigg(\df{A}{A_0}\bigg)^z\;.
\end{aligned}
\end{equation}
With the choice of parameters $m$, $n$ $z$ and $C$ in Table \ref{table:parameters}, there exist at least one solution of $G\left(\frac{A}{A_0}\right)=0$. This means that the two characteristic fields are neither genuinely non-linear nor linearly degenerate. The consequences of this are unclear to the authors, and might require further investigations. See \cite{LeFloch:2002a} and \cite{Dafermos:2016a} for details. 
The Generalized Riemann Invariants (GRIs) for $\lambda_1$ and $\lambda_2$ characteristic fields are respectively
\begin{equation}
\begin{drcases}
\lambda_1 - \text{GRI} : \quad u+\int\f{c\left(A\right)}{A}\mathrm{d}A=constant\;, \quad \\
\lambda_2 - \text{GRI} : \quad u-\int\f{c\left(A\right)}{A}\mathrm{d}A=constant\;.
\end{drcases}
\end{equation} 
In the present work, the generalized Riemann invariants will be used to couple valves with lymphangions and to impose the pressure at the terminal interfaces of the collector.

\subsection{A dynamical Electro-Fluid-Mechanical Contraction \\ (EFMC) model} \label{sec:FitzHugh}
Lymphatic contractions are preceded by action potentials \cite{Weid:2004a,Weid:1997a,Breslin:2014a}. Action potentials occur at smooth muscle cells of the collecting lymphatic wall. These are triggered by the rapid influx of $Ca^{2+}$ and other ions through different type of channels (low voltage T-channels, high-voltage L-channels, stretch-activated channels, gap junctions) \cite{Munn:2015a}. From the pioneering work of Hodgkin and Huxley in 1952 \cite{Hodgkin:1952a} on action potentials in neurons, mathematical models for action potentials have been widely used in the literature and applied in several contexts \cite{Aliev:1996a,Nagumo:1962a,Panfilov:1993a,Kogan:1991a}. 
Here we propose an Electro-Fluid-Mechanical Contraction (EFMC) model for lymphatics, based on the FitzHugh-Nagumo model for action potentials. This model combines the electrical activity of lymphangions (action potential) with fluid-mechanical feedback (stretch of the lymphatic wall and wall shear stress), and the mechanical variation of the lymphatic wall properties (contractions). 

The modelling system of ODEs is
\begin{equation}\label{eq:FHNSystem}
\dttotal \mathbf{Y}=\mathbf{L}\left(\mathbf{Y}\right)\;,
\end{equation}
where 
\begin{equation}\label{eq:FHNModel}
\mathbf{Y}\left(t\right)=\colvecfour{v\left(t\right)}{w\left(t\right)}{I\left(t\right)}{s\left(t\right)}\;,\quad \mathbf{L}\left(\mathbf{Y}\right)=\colvecfour{a_1\big[v\big(v-a_2\big)\big(1-a_3v\big)-w+vI\big]}{b_1v-b_2w}{f_I(\bar{A},\bar{\tau},v,w,I)}{f_s(v,w,s)}\;,
\end{equation}
and
\begin{equation}\label{eq:Trigger}
\resizebox{\linewidth}{!}{
$
f_I(\bar{A},\bar{\tau},v,w,I)=\begin{dcases}
\left(k_{Ca}^{(1)}+k_{Ca}^{(2)}\bigg(\df{\bar{A}}{A_{Ca}}\bigg)^{n_{Ca}}\right)f_{NO}\big(\bar{\tau}\big)\;,& \sqrt{v^2+w^2}\leq R_I\;, \\
-Ik_{rel}\;, & \sqrt{v^2+w^2}>R_I\;, 
\end{dcases}\;
$
}
\end{equation}
\begin{equation}
f_{NO}\big(\bar{\tau}\big)=1-k_{NO}\left(\df{2}{1+exp\left({-\left|\df{\bar{\tau}}{\tau_{NO}}\right|^{n_{NO}}}\right)}-1\right)\;,
\end{equation}
and
\begin{equation}\label{eq:ContractionFunction}
f_s(v,w,s)=\begin{dcases}
+c_1vw\left(1-s\right)\;, & v> 0\;, \\
-c_2s\left|1-w\right|_+\;, & v\leq 0\;,
\end{dcases}
\end{equation}
where $\left| x \right|_+$ is the positive value of $x$
\begin{equation}
\left|x\right|_+=\begin{dcases}
x\;, & x> 0\;, \\
0\;, & x\leq 0\;.
\end{dcases}
\end{equation}
The unknowns of the above system are: the {\em excitable variable} $v(t)$ (membrane potential), the {\em recovery variable} $w(t)$, the {\em stimulus} $I(t)$, and the {\em contraction state} $s(t)$ introduced in Eq. \eqref{eq:ContractionYoungModulus}. The first two equations of \eqref{eq:FHNModel} are the classical FitzHugh-Nagumo (FHN) model with a simple modification. In the classical formulation of the FHN model, the stimulus $I$ has a constant value, while here it varies in time and multiplies the excitable variable $v$. As described in the next section, this modification allows us to control the equilibrium of the stationary solution and to trigger action potentials under certain conditions on the stimulus $I$. 

The evolution in time of $I$ is controlled by the function $f_I$ and phenomenologically depends on three mechanisms: 
(1) environmental calcium influx, (2) stretch-activated calcium influx, and (3) contraction inhibitions induced by WSS. The environmental baseline influx is regulated by the parameter $k_{Ca}^{(1)}$. The stretch-activated calcium influx is regulated by the parameters $k_{Ca}^{(2)}$, $A_{Ca}$ and $n_{Ca}$. The contraction inhibitions induced by WSS are regulated by the function $f_{NO}$, which depends on parameters $k_{NO}$, $\tau_{NO}$ and $n_{NO}$. The function $f_{NO}$ is bounded by $1-k_{NO}$ and $1$, namely
\begin{equation}
\displaystyle \lim_{\left|{\bar{\tau}}\right| \to +\infty}f_{NO}(\tau)=1-k_{NO}\leq f_{NO}\leq 1=f_{NO}(0)\;.
\end{equation}
The contraction state is controlled by the function $f_s$, which depends on parameters $c_1$ and $c_2$. Functions $f_{I}$ and $f_{NO}$ are evaluated at the space-averaged cross-sectional area and at the space-averaged WSS, respectively, at the current time
\begin{equation}
\bar{A}(t)=\f{1}{L}\int_{0}^LA(x,t)\mathrm{d}x\;, \quad \bar{\tau}(t)=\f{1}{L}\int_{0}^L\tau(x,t)\mathrm{d}x\;,
\end{equation}
where $L$ is the length of the lymphangion.
Parameters of the EFMC model are found in Table \ref{table:parameters}.

\subsubsection{Mathematical analysis of the modified FitzHugh-Nagumo model}\label{sec:equilibriumFHN}
Here we analyse the modified FitzHugh-Nagumo model on which the EFMC model is based. First we find the stationary state solution, and then we study its nature depending on the stimulus $I$. 
\begin{figure*}
\begin{center}
\includegraphics[width=1\textwidth]{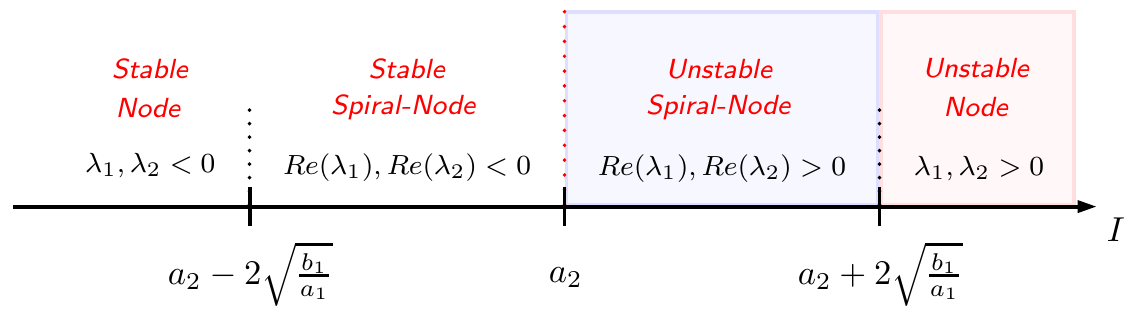}
\end{center}
\caption{\scriptsize {\bf Stability analysis of the stationary point $(0,0)$ of the modified FitzHugh-Nagumo model.} The nature of the stationary point depends on the stimulus $I$. For $I<a_2-2\sqrt{\f{b_1}{a_1}}$ and $I>a_2+2\sqrt{\f{b_1}{a_1}}$, the eigenvalues of the modified FHN model are real while for $a_2-2\sqrt{\f{b_1}{a_1}}<I<a_2+2\sqrt{\f{b_1}{a_1}}$, the eigenvalues are imaginary. Action potentials can be periodically triggered when $I>a_2$.} \label{fig:stabilityEFMC}
\end{figure*}
The stationary points are found by solving the following system
\begin{equation}\label{eq:problemasteady}
F_{FHN}(v,w)=\colvec{a_1\big[v\big(v-a_2\big)\big(1-a_3v\big)-w+vI\big]}{b_1v-b_2w}=\colvec{0}{0}\;.
\end{equation}
From now on, we assume $b_2=0$. Therefore, the only stationary point of \eqref{eq:problemasteady} is $(v,w)=(0,0)$ and is independent of $I$. To study the nature of the stationary point, one has to evaluate the Jacobian of $F_{FHN}$ at the stationary state and study the sign of its eigenvalues. The resulting eigenvalues are
\begin{equation}
\lambda_{1,2}=\df{-a_1(a_2-I)\pm \sqrt{a_1^2(a_2-I)^2-4a_1b_1}}{2}\;.
\end{equation}
Then we conclude that the stationary solution $(0,0)$ is a 

\begin{align*}
&\text{{\bf stable node}} &\text{when  } \quad & I< a_2-2\sqrt{\df{b_1}{a_1}}\;, \\
&\text{{\bf stable spiral-node} } &\text{when   }  \quad &a_2-2\sqrt{\df{b_1}{a_1}}<I< a_2\;, \\
&\text{{\bf unstable spiral-node} } &\text{when   }  \quad &a_2< I <a_2+2\sqrt{\df{b_1}{a_1}}\;, \\
&\text{{\bf unstable node}} &\text{when   }  \quad &a_2+2\sqrt{\df{b_1}{a_1}}<I\;.
\end{align*}
When $I< a_2$, the nature of the stationary solution is stable, and therefore action potentials are not automatically triggered. When $I> a_2$ the nature of the stationary solution is unstable, and therefore action potentials can be triggered. Fig. \ref{fig:stabilityEFMC} summarizes the study of the nature of the stationary solution. 

We stress that the mathematical analysis reported here is based on well-known mathematical criteria and is valid only near the stationary solution. If $I>a_2$ and the excitable and recovery variables are not near enough to the stationary solution $(0,0)$, then we might not have an action potential. In other words, when $I$ exceeds the threshold level $a_2$, the stationary solution becomes unstable, and an action potential develops as soon as the variables $v$ and $w$ are near enough to the stationary solution. 
For a time-varying stimulus $I$, we assume that the needed stimulus $\tilde{I}$ to trigger an action potential is between a minimum and a maximum value $\tilde{I}_{min}$ and $\tilde{I}_{max}$, defined as follows:
\begin{equation}\label{eq:Iminmax}
\tilde{I}_{min}:=a_2\;, \quad \tilde{I}_{max}:=a_2+2\sqrt{\df{b_1}{a_1}}\;.
\end{equation} 
These two values will be useful to estimate the frequency of contractions of the EFMC model. 

Note that the initial condition for the excitable variable $v$ or the recovery variable $w$ needs to be different from zero. Otherwise, action potentials do not occur in time.

\subsubsection{Qualitatively analysis of the EFMC model}
Here we qualitatively explain the EFMC model. Fig. \ref{fig:FHN} depicts a representative numerical result of the EFMC model. Excitable variable $v$, recovery variable $w$, contraction state $s$ and stimulus $I$ are shown for two lymphatic cycles. When the excitable variable $v$ and the recovery variable $w$ are near the stationary state ($\sqrt{v^2+w^2}<R_{I}$), the stimulus $I$ increases as given by Eq. \eqref{eq:Trigger}. When a certain value $\tilde{I}$ is reached (in this case $\tilde{I}_{min}<\tilde{I}<\tilde{I}_{max}$), an action potential is triggered: variables $v$ and $w$ perform a cycle, increase in absolute value and move far from the stationary state ($\sqrt{v^2+w^2}>R_{I}$) and consequently the stimulus $I$ exponentially decreases to zero. The state of contraction $s$ increases until $v>0$, and then decreases to zero, see Eq. \eqref{eq:ContractionFunction}. When the action potential ends, variables $v$ and $w$ return to the equilibrium point. From there on, the stimulus $I$ restarts to increase and possibly triggers a new contraction. 

Fig. \ref{fig:FHN_phase} shows the numerical result in the phase space. Three nullclines are represented. They are: the nullcline for the recovery variable $\dt w=0 : v=0$, and two representative nullclines with $I=0$ and $I=\tilde{I}$ for the excitable variable $ \dt v =0 : w= v\big(v-a_2\big)\big(1-a_3v\big)$ and $ \dt v =0 : w = v\big(v-a_2\big)\big(1-a_3v\big)+v\tilde{I}$. As soon as a contraction occurs, the stimulus $I$ decreases and the third nullcline tends on the second one. The sphere of radius $R_I$ centred at $v=w=0$, intersection of all the nullclines, divides the phase space into two regions. We call the region outside the sphere the {\em excited region}, while the latter the {\em activation region}. In the excited region, the solution quickly performs a cycle, while the stimulus $I$ exponentially decreases to zero. In the activation region, the numerical solution tends to the equilibrium state, while the stimulus $I$ increases. 

The time required for the numerical solution to perform a cycle (from the excited region into the activation region) can be numerically evaluated and is denoted with $t_{\one}$. 
The time required to activate an action potential is denoted with $t_{\two}$.

\begin{figure*}[t]
\includegraphics[width=1\textwidth]{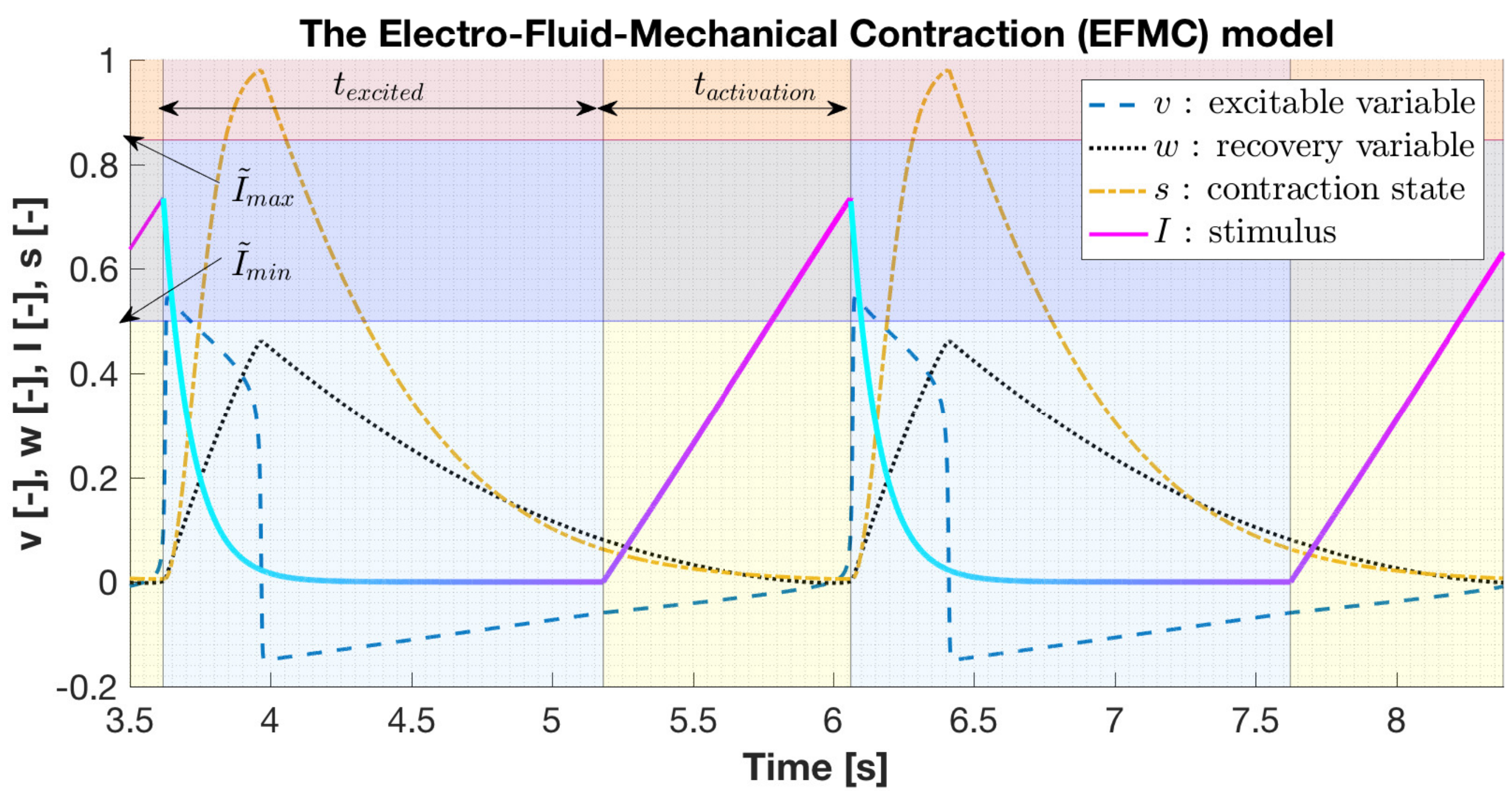} 
\caption{\scriptsize {\bf Illustration of the EFMC model in the time domain}. Results show the excitable variable, recovery variable, contraction state and the stimulus for two lymphatic cycles. Minimum and maximum triggering values $\tilde{I}_{min}$ and $\tilde{I}_{max}$ are also shown. Dark blue and red shaded areas illustrate the {\em unstable spiral-node} and the {\em unstable node} regions.
The {\em excited} and {\em activation} regions can be determined by the different shaded colors (blue and yellow). Here we solved the system of ODEs \eqref{eq:FHNSystem} with initial condition $v(0)=0.001$, $w(0)=I(0)=s(0)=0$. The parameters of the EFMC model were taken from Table \ref{table:parameters}, but here we set $f_{min}=5$ min$^{-1}$, and we assumed $\bar{A}=7.9324A_0$ and $\bar{\tau}=0$. The colour gradient of the stimulus $I$ is the same as in Fig. \ref{fig:FHN_phase}. }\label{fig:FHN}
\end{figure*}

\begin{figure}
\centering
\includegraphics[width=0.5\textwidth, height=1\textheight,keepaspectratio]{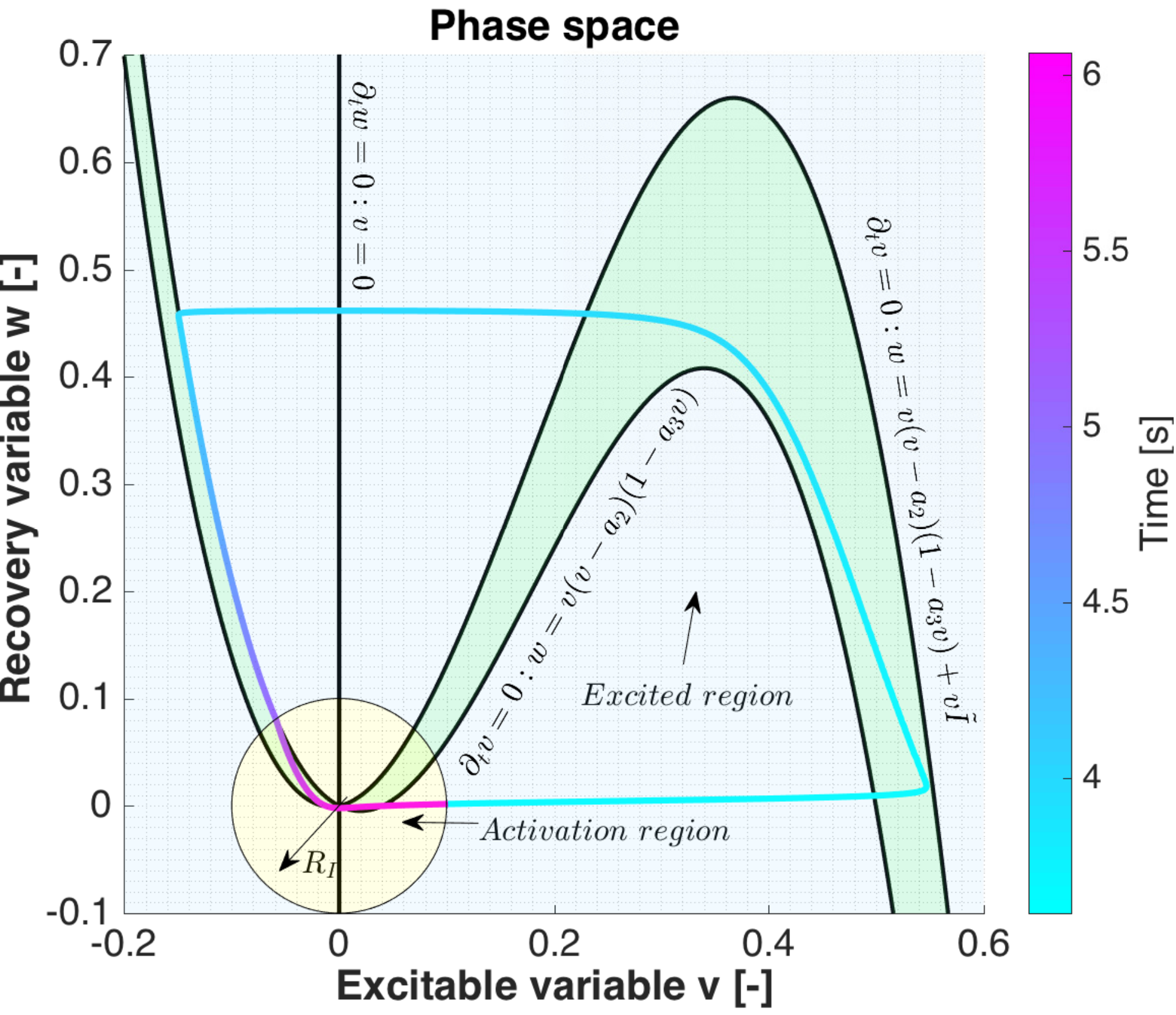}
\caption{\scriptsize {\bf Illustration of the EFMC model in phase space of a representative lymphatic cycle.} The circle of radius $R_I$ centred in the stationary point $(0,0)$ divides the space in two regions: the {\em activation} region (yellow) and the {\em excited} region (blue). Three nullclines are shown: the nullcline for the recovery variable $\dt w=0$, and two representative nullclines with $I=0$ and $I=\tilde{I}$ for the excitable variable $ \dt v =0$. As soon as a contraction occurs, the stimulus $I$ decreases and the third nullcline tends to the second one. Results in the time domain can be seen in Fig. \ref{fig:FHN}.} \label{fig:FHN_phase}
\end{figure}

\subsubsection{Analysis of the EFMC frequency}

\begin{figure*}[t]
\centering
\subfloat[][]{\label{fig:a}\includegraphics[width=0.33\textwidth]{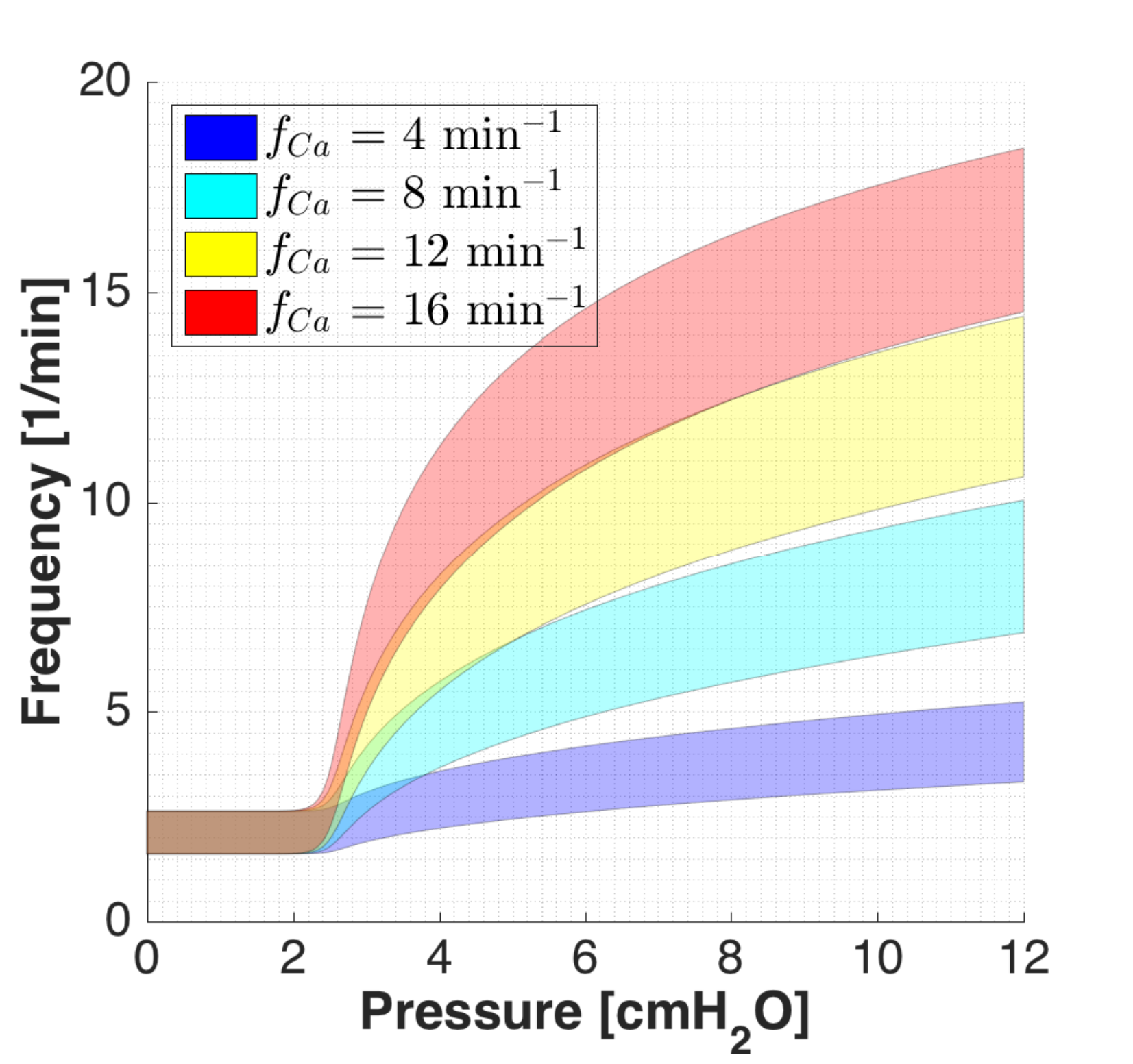}}
\subfloat[][]{\label{fig:b}\includegraphics[width=0.33\textwidth]{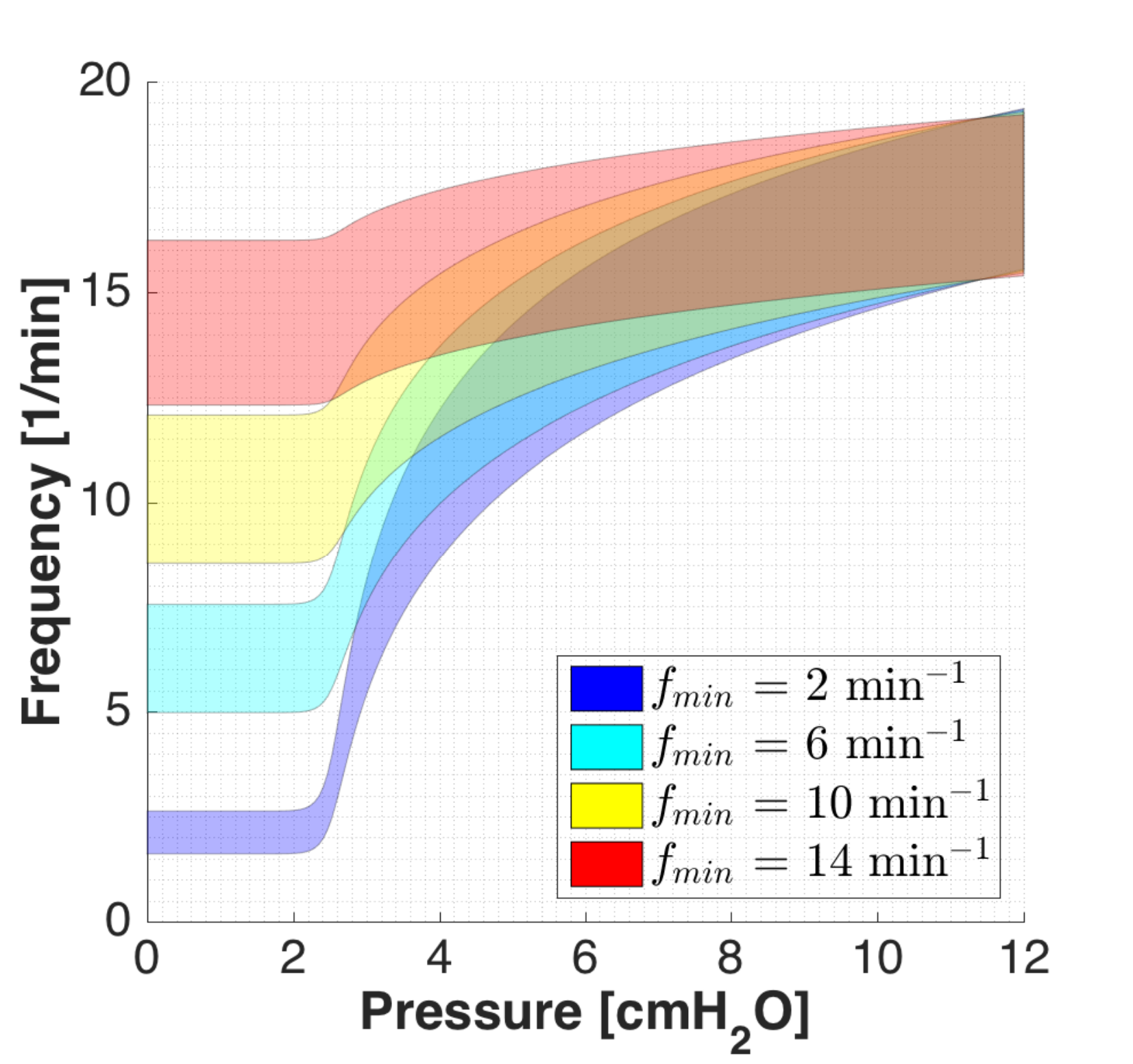}}
\subfloat[][]{\label{fig:c}\includegraphics[width=0.33\textwidth]{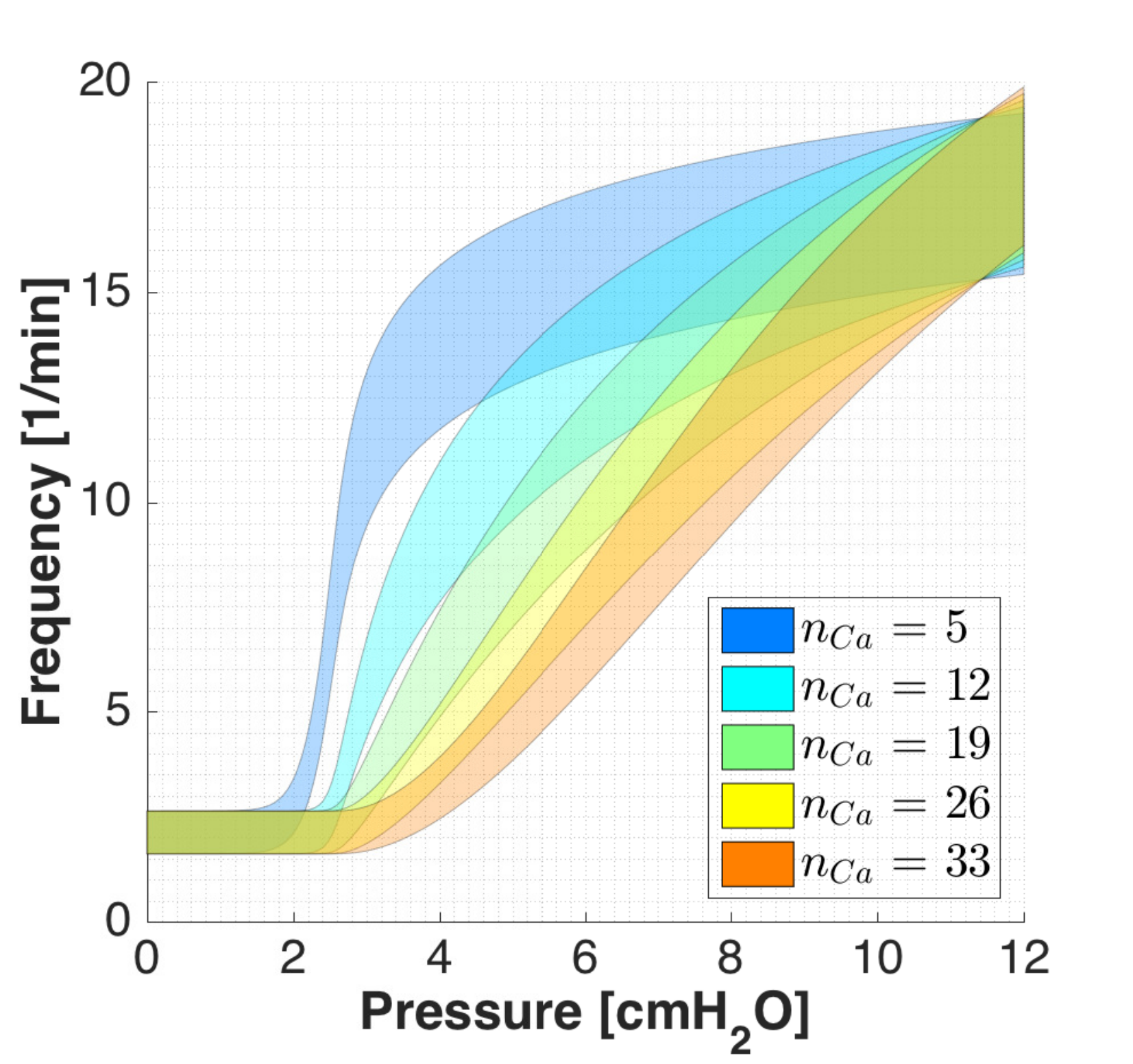}}  \\
\subfloat[][]{\label{fig:d}\includegraphics[width=0.33\textwidth]{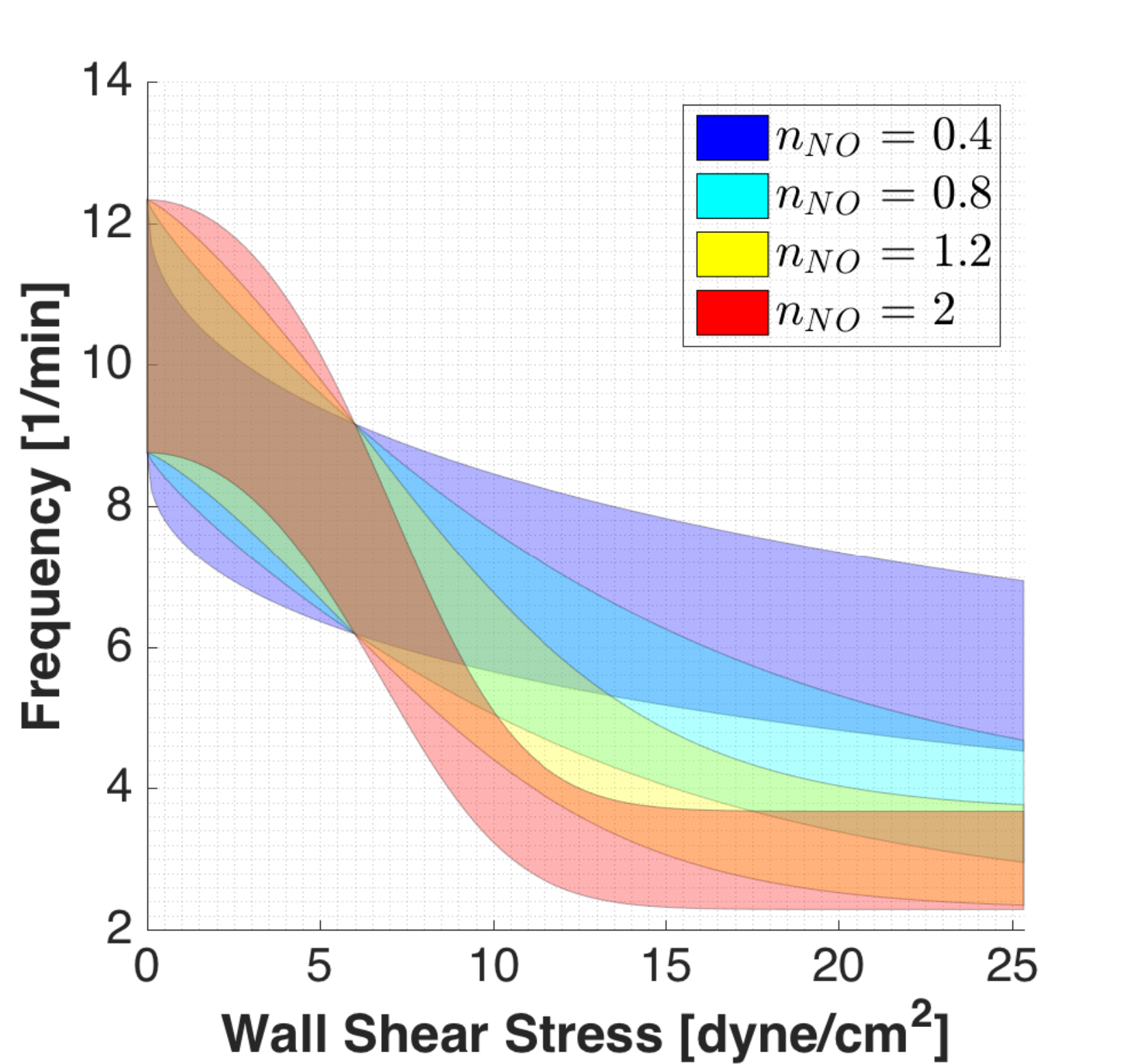}} 
\subfloat[][]{\label{fig:e}\includegraphics[width=0.33\textwidth]{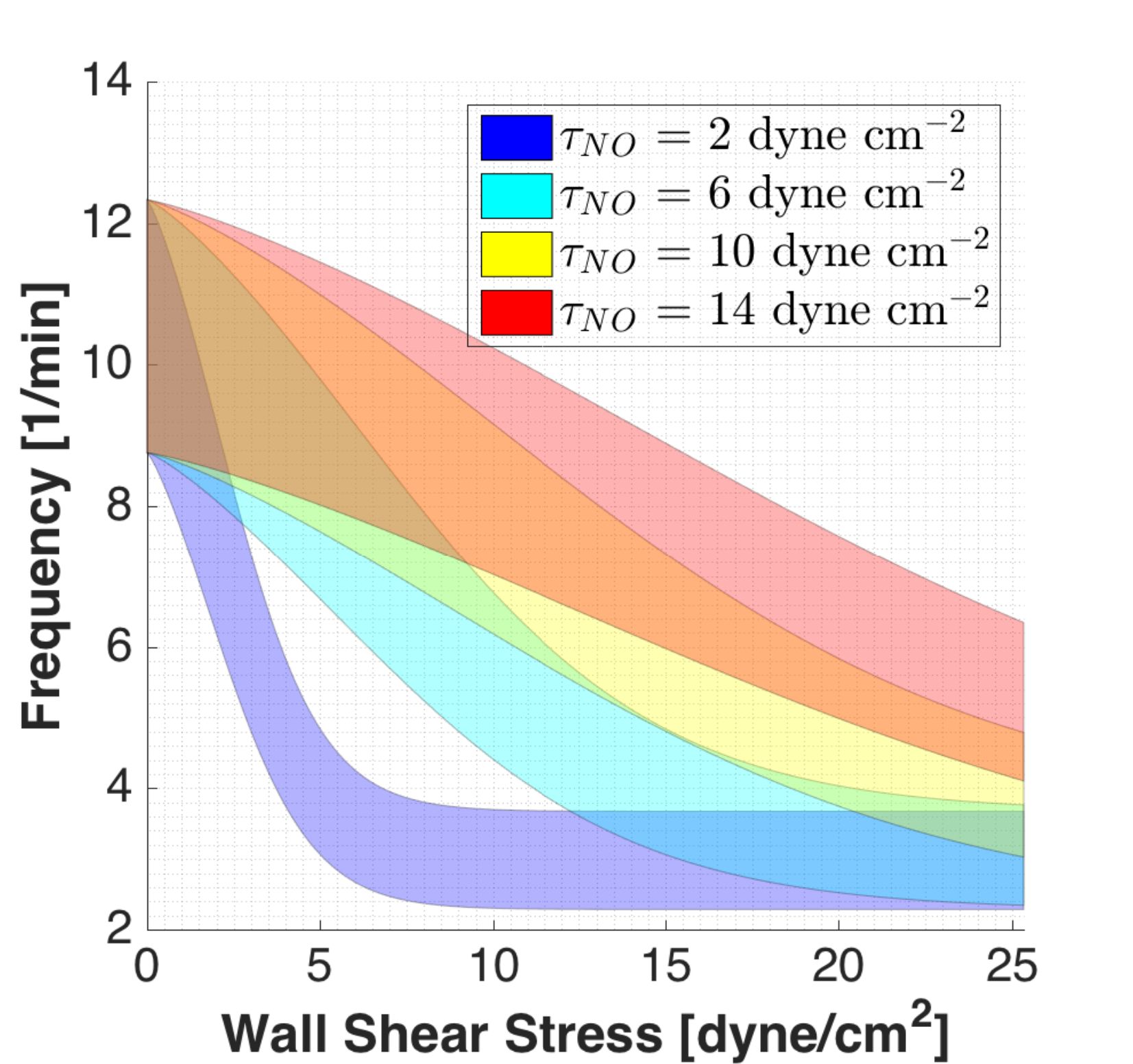}} 
\subfloat[][]{\label{fig:f}\includegraphics[width=0.33\textwidth]{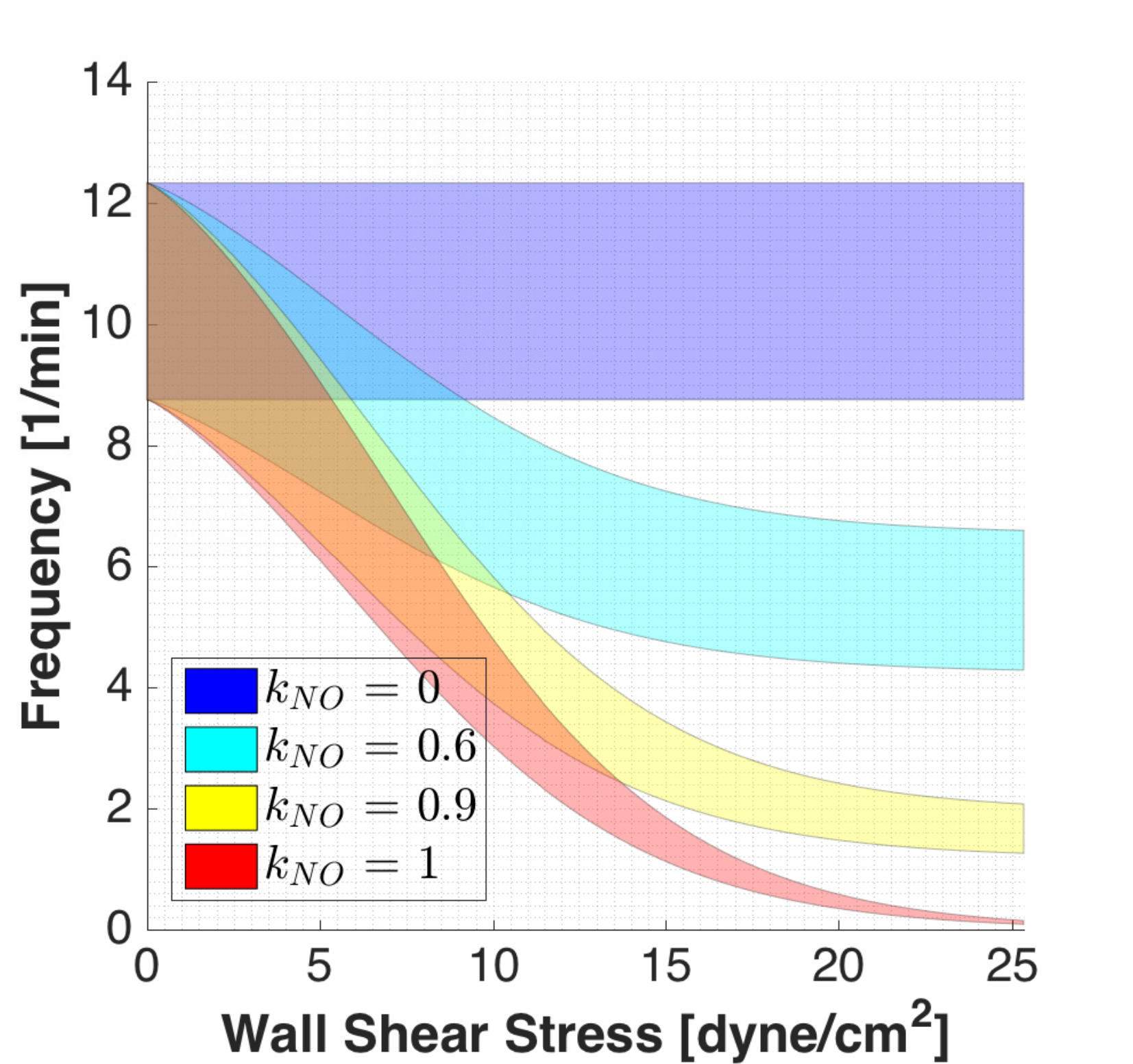}}
\caption{\scriptsize {\bf Theoretical effects of different EFMC model parameters on the contraction frequency}. In the first row we show theoretical results for pressure against frequency varying $f_{Ca}$, $f_{min}$ and $n_{Ca}$. In the second row we show theoretical results for WSS againts the frequency varying $n_{NO}$, $\tau_{NO}$ and $k_{NO}$. Results are based on Eq. \eqref{eq:frequency} and the coloured-shaded area are bounded by using triggering stimulus values $\tilde{I}_{min}$ and $\tilde{I}_{max}$.}\label{fig:EFMC}
\end{figure*}

Lymphangions contract differently according to the location of the vessel, the stretch of the lymphatic wall and wall shear stress feedback \cite{Gashev:2004a}. Here we aim to analyse the frequency of contraction of the EFMC model and to estimate parameters $k_{Ca}^{(1)}$ and $k_{Ca}^{(2)}$ in order to prescribe a baseline frequency $f_{min}$ and a frequency $f_{Ca}$ at $\bar{A}=A_{Ca}$. The time $t_{total}$ between two cycles can be written as follows:
\begin{equation}
t_{total}=t_{\one}+t_{\two}\;,
\end{equation}
and the relative frequency is 
\begin{equation}\label{eq:Frequency}
f=\f{1}{t_{total}}=\f{1}{t_{\one}+t_{\two}}\;.
\end{equation}
The excited time $t_{\one}$ can be assumed to be constant and can be evaluated numerically from the FHN model. The activation time, on the other hand, depends strongly on the rate of increase of $I$ given by \eqref{eq:Trigger}. We now estimate the activation time, namely the time required for the stimulus $I$ to attain a certain triggering value $\tilde{I}$. 
Near the stationary solution $v=w=0$, it is reasonable to assume $\sqrt{v^2+w^2}<R_I$. We solve the following initial value problem
\begin{align}
\resizebox{\linewidth}{!}{
$
\begin{drcases}
\text{ODE:}\quad &\dttotal I (t)=\left(k_{Ca}^{(1)}+k_{Ca}^{(2)}\left(\df{\bar{A}\left(t\right)}{A_{Ca}}\right)^{n_{Ca}}\right)f_{NO}\big(\bar{\tau}\left(t\right)\big)\;, \quad t\geq t_{\one}\;, \\
\text{IC:}\quad & I(t_{\one})=0\;,
\end{drcases}
$}
\end{align}
whose solution is
\begin{equation} \label{eq:IVPactivation}
\resizebox{\linewidth}{!}{
$
\begin{drcases}
I(t)=I(t_{\one})+\int_{t_{\one}}^{t} \left(k_{Ca}^{(1)}+k_{Ca}^{(2)}\left(\df{\bar{A}(t)}{A_{Ca}}\right)^{n_{Ca}}\right)f_{NO}\left(\bar{\tau}(t)\right)\mathrm{d} t\;, & t\geq t_{\one}\;, \\
I(t_{\one})=0\;.
\end{drcases}
$}
\end{equation}
We assume that during the activation time, $\bar{A}$ and $\bar{\tau}$ are constant in time because the lymphangion is already at the end of the diastolic phase. Thus, the above integral can be exactly computed and is
\begin{equation}
I(t)=t\bigg(k_{Ca}^{(1)}+k_{Ca}^{(2)}\bigg(\df{\bar{A}}{A_{Ca}}\bigg)^{n_{Ca}}\bigg)f_{NO}\left(\bar{\tau}\right)\;.
\end{equation}
Consequently, the activation time $t_{\two}$ required to attain a triggering value $\tilde{I}$ is 
\begin{equation}\label{eq:ActivationTime}
t_{\two}=\df{\tilde{I}}{\left(k_{Ca}^{(1)}+k_{Ca}^{(2)}\bigg(\df{\bar{A}}{A_{Ca}}\bigg)^{n_{Ca}}\right)f_{NO}\left(\bar{\tau}\right)}\;.
\end{equation}
The maximum activation time (when $\textstyle \left({\bar{A}}/{A_{Ca}}\right)^{n_{Ca}} \approx 0$) and the activation time at $\bar{A}=A_{Ca}$ both at zero WSS ($\bar{\tau}=0$) are
\begin{equation}\label{eq:imposedfrequency1}
t_{\two}^{max}=\f{\tilde{I}}{k_{Ca}^{(1)}}\;, \quad t_{\two}^{Ca}=\f{\tilde{I}}{k_{Ca}^{(1)}+k_{Ca}^{(2)}}\;.
\end{equation}
The maximum activation time, corresponding to the minimum frequency, depends on parameter $k_{Ca}^{(1)}$. In our model, parameter $k_{Ca}^{(1)}$ phenomenologically represents the environmental calcium influxes. Indeed, in the limiting case in which there is no environmental calcium influxes ($k_{Ca}^{(1)}=0$), the activation time becomes infinite, which means that the lymphangion does not autonomously contract. Parameter $k_{Ca}^{(2)}$, on the other hand, phenomenologically regulates the stretch-induced calcium influxes. The greater the parameter $k_{Ca}^{(2)}$, the less the activation time, and thus the greater the influence of the stretch of the lymphatic wall on the frequency contractions.

If we assume the frequencies $f_{min}$ and $f_{Ca}$, corresponding to $t^{max}_{\two}$ and $t^{Ca}_{\two}$ respectively, to be known, then we have
\begin{equation}\label{eq:imposedfrequency}
\df{1}{f_{min}}=t_{\one}+t^{max}_{\two}\;, \quad \df{1}{f_{Ca}}=t_{\one}+t^{Ca}_{\two}\;.
\end{equation}
Using \eqref{eq:imposedfrequency1} and \eqref{eq:imposedfrequency}, we can explicitly find parameters $k_{Ca}^{(1)}$ and $k_{Ca}^{(2)}$: 
\begin{equation}\label{eq:parameterkCa}
\resizebox{\linewidth}{!}{
$
k_{Ca}^{(1)}=\df{\tilde{I}}{\df{1}{f_{min}}-t_{\one}}\;,\quad k_{Ca}^{(2)}=\df{\tilde{I}}{\df{1}{f_{Ca}}-t_{\one}}-\df{\tilde{I}}{\df{1}{f_{min}}-t_{\one}}\;.
$}
\end{equation}
To assure positive activation times $t^{max}_{\two}>0$, $t^{Ca}_{\two}>0$, the following conditions need to be satisfied:
\begin{equation}
\df{1}{t_{\one}}>f_{min}\;, \quad f_{Ca}>f_{min}\;.
\end{equation}
Here we assume the triggering value $\tilde{I}$ to be the mean value of $I_{max}$ and $I_{min}$ defined in Eq. \eqref{eq:Iminmax}, namely
\begin{equation} \label{eq:Imean}
\tilde{I}_{mean}=\df{\tilde{I}_{max}+\tilde{I}_{min}}{2}=a_2+\sqrt{\df{b_1}{a_1}}\;.
\end{equation}
Numerical results confirmed that this is a good choice, even though $\tilde{I}_{min}$ and $\tilde{I}_{max}$ can be used as triggering values too. Substituting $k_{Ca}^{(1)}$ and $k_{Ca}^{(2)}$ and the activation time $t_{\two}$ defined in \eqref{eq:ActivationTime} into \eqref{eq:Frequency}, one obtain a frequency function as
\begin{equation}\label{eq:frequency}
f\left(\bar{A},\bar{\tau},\tilde{I}\right)=\f{1}{t_{\one}+t_{\two}\left(\bar{A},\bar{\tau},\tilde{I}\right)}\;.
\end{equation}
Then, one can easily prove the following inequalities
\begin{equation}\label{eq:property1}
f\left(\bar{A},0,\tilde{I}\right)>f\left(\bar{A},\bar{\tau}_1,\tilde{I}\right)>f\left(\bar{A},\bar{\tau}_2,\tilde{I}\right)\;, \quad \left|\bar{\tau}_1\right|<\left|\bar{\tau}_2\right|\;,
\end{equation}
and
\begin{equation}\label{eq:property2}
f\left(\bar{A}_1,\bar{\tau},\tilde{I}\right)<f\left(\bar{A}_2,\bar{\tau},\tilde{I}\right)\;, \quad \bar{A}_1<\bar{A}_2\;.
\end{equation}
The first property \eqref{eq:property1} says that the frequency decreases as the absolute value of the WSS increases, and maximum contractions are attained at zero WSS. The second property \eqref{eq:property2} says that the frequency increases as the cross-sectional area increases, that is the higher the intraluminal pressure, the higher the frequency of the contraction.

Fig. \ref{fig:EFMC} shows theoretical results of the EFMC model. The parameters were taken from Table \ref{table:parameters}. The top row shows pressure against frequency, for different values of $f_{Ca}$, $f_{min}$ and $n_{Ca}$, while the bottom row shows WSS against frequency for a given pressure, for different values of $n_{NO}$, $\tau_{NO}$ and $k_{NO}$. In the top row, we considered a linear variation for pressure between 0 and 12 cmH$_2$O and obtained the corresponding cross-sectional area $A$ using the inverse of the tube law \eqref{eq:pressure}. Then we depicted a coloured shaded area for the frequency in Eq. \eqref{eq:frequency} using $\tilde{I}_{min}$ and $\tilde{I}_{max}$ as triggering values. Indeed, note that
\begin{equation}
f\left(\bar{A},\bar{\tau},\tilde{I}_{max}\right)<f\left(\bar{A},\bar{\tau},\tilde{I}_{mean}\right)<f\left(\bar{A},\bar{\tau},\tilde{I}_{min}\right)\;.
\end{equation}
Numerical results of the complete model composed of the one-dimensional lymph flow equations in Eq. \eqref{eq:Lymph} and of the EFMC model in Eq. \eqref{eq:FHNSystem} shows that the resulting frequency follows the theoretical trend shown in Fig. \ref{fig:EFMC}.
These results show that it is possible to imitate the experimental measurements of a specific lymphangion by fitting parameters $f_{min}$, $f_{Ca}$ and $n_{Ca}$. 
In the bottom row, we considered a linear variation of flow $q$ between 0 and 240 $\mu$L min$^{-1}$ and then obtained the WSS from Eq. \eqref{eq:WSS} assuming $\bar{A}=7A_{0}$. These results show that given the experimental trend of the frequency of a specific lymphangion for different flows, it is possible to adjust parameters $n_{NO}$, $\tau_{NO}$ and $k_{NO}$ to imitate that trend. For instance, from experimental measurements we know that the WSS has a greater negative chronotropic effect on the thoracic duct than on femoral lymphatics \cite{Gashev:2004a}, and this can be modelled by using $k_{NO}\approx 1$ for the thoracic duct and, for example, $k_{NO}\approx 0.5$ for femoral lymphatics.

\subsection{A lumped-parameter model for lymphatic valves} \label{sec:valves}

Here we follow the mathematical model proposed by Mynard et al. \cite{Mynard:2012a}, which has been already adapted to simulate human venous valves by \cite{Mueller:2014b, Mueller:2015a,Toro:2015aa}. The time variation of the flow across the valve $q_v\left(t\right)$ is approximated as
\begin{equation}\label{eq:ValveModel}
\dttotal q_v =\f{1}{L(\xi)}\left(\Delta p\left(t\right) -R\left(\xi\right)q_v - B\left(\xi\right)q_v|q_v| \right)\;,
\end{equation}
where
\begin{equation}
\Delta p\left(t\right) =p_{u}\left(t\right)-p_{d}\left(t\right)\;.
\end{equation}
Here $p_{u}$ and $p_{d}$ are the upstream and downstream pressures, respectively. Coefficients $B$, $L$ and $R$ are the Bernoulli resistance, the lymphatic inertia and the viscous resistance to flow, given respectively as
\begin{equation}\label{eq:ValveParameters}
B(\xi)=\df{\rho}{2A_{eff}^2(\xi)}\;,\quad L(\xi)=\rho \df{L_{eff}}{A_{eff}(\xi)}\;, \quad R(\xi)=\f{2\left(\gamma +2 \right)\pi \mu}{A_{eff}^2(\xi)}L_{eff}
\end{equation}
where $L_{eff}$ is the effective length and $A_{eff}$ is the effective area, which varies from a minimum value to a maximum value as
\begin{equation}
\begin{aligned}
A_{eff}(\xi)=&A_{eff,min}+\xi (t)\big(A_{eff,max}-A_{eff,min}\big)\;, \\& \xi \in [0,1]\;.
\end{aligned}
\end{equation}
Compared to the work of Mynard et al. \cite{Mynard:2012a}, we have added the Poiseuille-type viscous losses insofar as numerical experiments showed that this term plays an important role in the lymphatic context. See \cite{Alastruey:2008b} for the derivation of $L$ and $R$.
The minimum and the maximum effective areas are evaluated as follow
\begin{equation}
A_{eff,min}=M_{rg}A_0\;, \quad A_{eff,max}=M_{st}A_0\;,
\end{equation}
where $M_{rg}$ is a parameter that controls the minimum closure. A normal closure is modelled by setting $M_{rg}=0$, while to model a regurgitant valve we use $0<M_{rg}\leq 1$. The parameter $M_{st}$ controls the maximum opening. A normal opening is modelled by setting $M_{st}=1$, while to model a stenotic valve we use $0\leq M_{st}< 0$. Then $A_0$ is here taken as the mean value between the cross-sectional areas at equilibrium of the adjacent lymphangions.
The {\em valve state} $\xi\left(t\right)$ is governed by the following ODE
\begin{equation}\label{eq:ValveState}
\dttotal \xi = f_{\xi}\left(\xi,t\right)=\begin{dcases}
K_{vo}(1-\xi)\big(\Delta p\left(t\right)-\Delta p_{open}\big)\;, & \Delta p\left(t\right)> \Delta p_{open}\;,\\
K_{vc}\xi\big(\Delta p\left(t\right)-\Delta p_{close}\big)\;, & \Delta p\left(t\right)< \Delta p_{close}\;,
\end{dcases}
\end{equation}
where $K_{vo}$ and $K_{vc}$ are the valve opening/closure rates, and $\Delta p_{vo}$ and $\Delta p_{vc}$ are the opening/closure threshold pressures. Davis et al. \cite{Davis:2011a} noticed that lymphatic valves are biased to stay open and display hysteresis. The authors observed that the threshold pressure to open and close the valve are strictly related to the transmural pressure of the lymphatic vessel. Bertram et al. \cite{Bertram:2014a} proposed formulas for $\Delta p_{open}$ and $\Delta p_{close}$ to imitate those experimental measurements. In the present work, we assume both the opening and closure threshold pressures to be zero. Hence, the system of ODEs to be solved is
\begin{equation}\label{eq:systemValve}
\dttotal \mathbf{Y}=\mathbf{L}\left(\mathbf{Y},t\right)\;,
\end{equation}
where 
\begin{equation}
\mathbf{Y}\left(t\right)=\colvec{q_v\left(t\right)}{\xi\left(t\right)}\;,\quad \mathbf{L}\left(\mathbf{Y},t\right)=\colvec{\f{1}{L\left(\xi\right)}\left(\Delta p\left(t\right) -R\left(\xi\right)q_v - B\left(\xi\right)q_v|q_v| \right)}{f_{\xi}\left(\xi,t\right)}\;.
\end{equation}
\section{Numerical methods} \label{sec:NumericalMethods}
\begin{figure*}[t]
\centering
\includegraphics[width=0.9\textwidth]{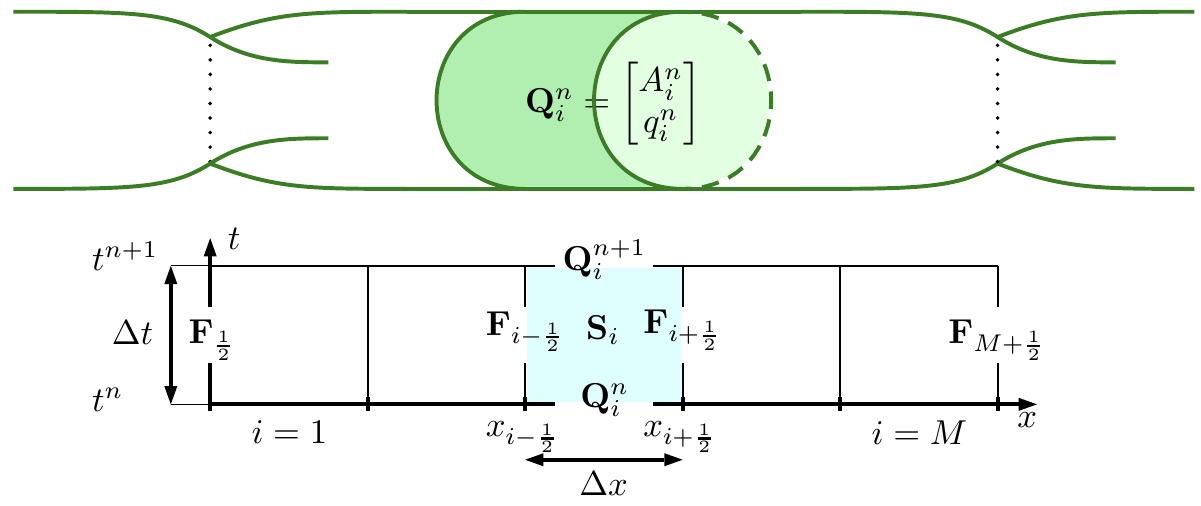} 
\caption{\scriptsize {\bf Framework for a finite volume scheme.} Top: illustratation of a computational volume for a lymphangion. Bottom: illustratation of the space-time control volume.} \label{fig:FiniteVolume}
\end{figure*}
 Here we briefly describe the finite volume schemes used for the one-dimensional lymph flow equations, explain how lymphatic valves and lymphangions are coupled, and illustrate the treatment of the boundary conditions at the terminal interfaces of the lymphangion. Then, we summarize the numerical methods used for the valves and the EFMC models.
\subsection{A finite volume method for the one-dimensional model} \label{sec:finitevolume}
Consider the system of $m$ hyperbolic balance laws
\begin{align} \label{eq:balancelaw}
\dt\Q +\dx \F(\Q)=\S(\Q)\;.
\end{align} 
By integrating \eqref{eq:balancelaw} over the control volume $V=[x_{i-\f{1}{2}},x_{i+\f{1}{2}}]\times[t^n,t^{n+1}]$ we obtain the exact formula
\begin{equation}\label{eq:FiniteVolumeFormula}
\Q^{n+1}_i=\Q^n_i-\f{\Delta t}{\Delta x}\big(\F_{i+\f{1}{2}}-\F_{i-\f{1}{2}}\big)+\Delta t \S_i\;,
\end{equation}
with definitions
\begin{equation}\label{eq:ConservativeVariable}
\Q^n_i =\f{1}{\Delta x}\intx \Q(x,t^n) \mathrm{d}x\;,
\end{equation}
\begin{equation}\label{eq:NumFluxNumSource}
\begin{drcases}
\F_{i+\f{1}{2}} = \f{1}{\Delta t}\int_{t^n}^{t^{n+1}} \F(\Q(x_{i+\f{1}{2}},t)) \mathrm{d}t\;, \\ \S_i=\f{1}{\Delta t \Delta x}\int_{t^n}^{t^{n+1}}  \intx \S(\Q(x,t)) \mathrm{d}x \mathrm{d}t\;.
\end{drcases}
\end{equation}
Eq. \eqref{eq:ConservativeVariable} gives the spatial-integral average at time $t=t^n$ of the conserved variable $\Q$ while Eqs. \eqref{eq:NumFluxNumSource} give the time-integral average at interface $x=x_{i+\f{1}{2}}$ of the physical flux $\F$ and the volume-integral average in $V$ of the source term $\S$.
Spatial mesh size and time step are $\Delta x=x_{i+\f{1}{2}}-x_{i-\f{1}{2}}$ and $\Delta t=t^{n+1}-t^n$ respectively. 
Finite volume methods for \eqref{eq:balancelaw} depart from \eqref{eq:FiniteVolumeFormula} to \eqref{eq:NumFluxNumSource}, where integrals are approximated, and then formula \eqref{eq:FiniteVolumeFormula} becomes a {\em finite volume method}, where the approximated integrals in \eqref{eq:NumFluxNumSource} are called {\em numerical flux} and {\em numerical source}, respectively. 
Here index $i$ runs from $1$ to $M$, where the cell $i=1$ is the leftmost cell with $x_{\frac{1}{2}}$ being the first interface, and the cell $i=M$ is the rightmost cell with $x_{M+\frac{1}{2}}$ being the last interface. See Fig. \ref{fig:FiniteVolume} for an illustration of the finite volume framework. 
To compute the time step $\Delta t$, the Courant-Friedrichs-Lewy condition is applied for each lymphangion 
\begin{equation}\label{eq:timeStep}
\displaystyle \Delta t^j = CFL \f{\Delta x^j}{\max\limits_{i=1,\dots, M^j} \left( |u_i^j|+c_i^j \right) }\;,
\end{equation}
with $CFL=0.9$. Superindex $j$ indicates the $j$-th lymphangion. Then, the time step $\Delta t$ to be used is the minimum of all the time steps, namely $\Delta t=\min\limits_{j}\left(\Delta t^j\right)$. 

In the present paper we used the SLIC method to evaluate the numerical fluxes within the domain ($\F_{\frac{3}{2}},\dots, \F_{M-\frac{1}{2}}$) \cite{Toro:2000a}. This method is second-order accurate in space and time and is based on the MUSCL-Hancock scheme where the Godunov upwind flux is replaced by the FORCE flux, see Section 14.5.3 of \cite{Toro:2009a} and references therein. The numerical source was approximated using a second order in space and time method, see Chapter 19 of \cite{Toro:2009a}. For the numerical fluxes at the boundaries ($\F_{\frac{1}{2}}$ and $\F_{M+\frac{1}{2}}$) we used a first-order Godunov-type method based on the solution of a classical Riemann problem at the interface.

\subsection{Coupling between valves and lymphangions} \label{sec:couplingvalves}
\begin{figure}[t]
\centering
\includegraphics[width=0.5\textwidth]{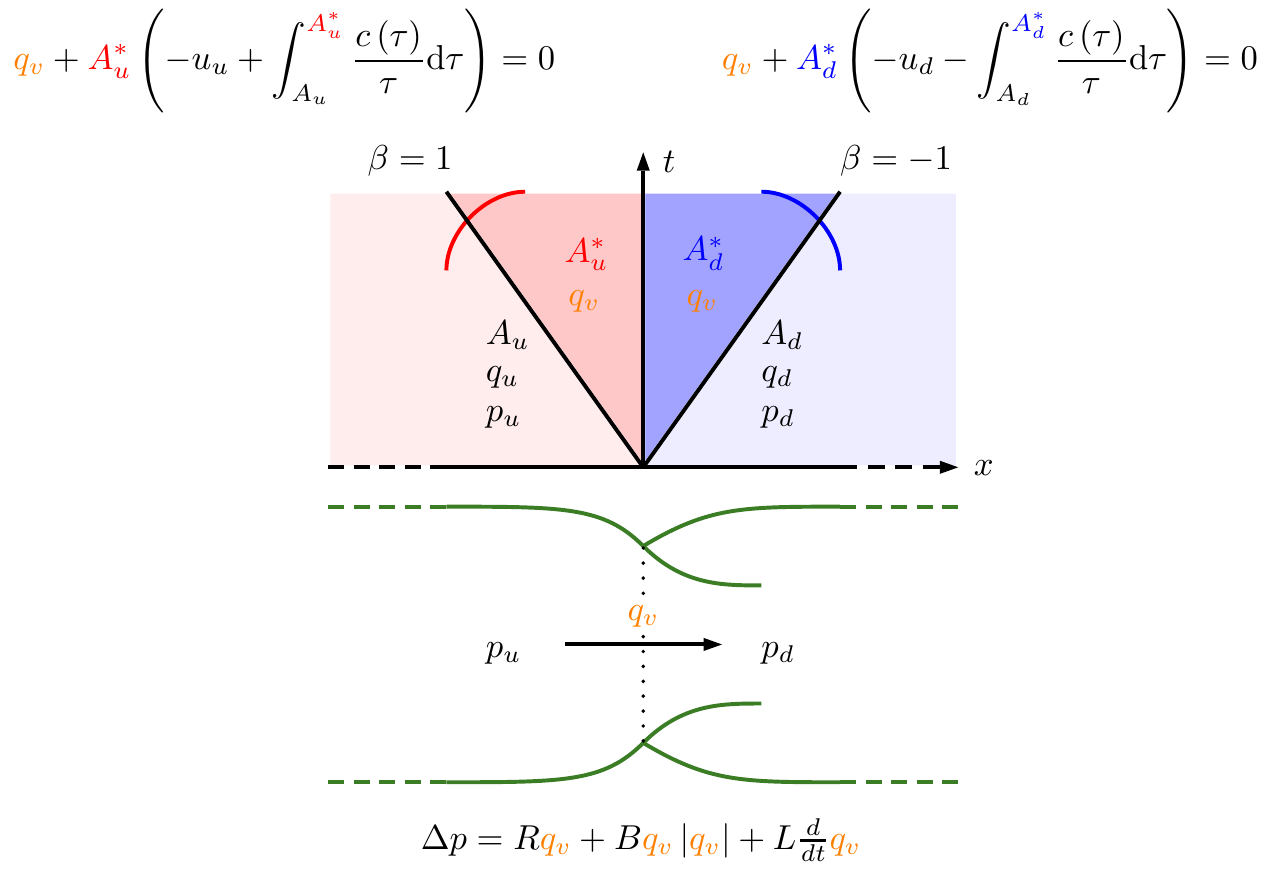} 
\caption{\scriptsize {\bf Illustration of the coupling method between two lymphangions and one valve.} Riemann invariants are used to couple the flow through the valve and the two lymphangions.} \label{fig:Coupling}
\end{figure}

Here we aim to couple valves and lymphangions. For each lymphangion, we need to calculate the numerical flux at the interface in which the valve is located, which can be either $\F_{\f{1}{2}}$ or $\F_{M+\f{1}{2}}$ according to Fig. \ref{fig:FiniteVolume}. There are three possible configurations for a lymphatic valve. It can be the leftmost or rightmost valve of a collector, or it can be interposed between two lymphangions.
In every case, the flow across the lymphatic valve is calculated from \eqref{eq:systemValve}, where the pressure gradient $\Delta p$ in \eqref{eq:ValveModel} is evaluated at the current time $t^n$ using either the two lymphangions, or one of the lymphangions and a prescribed time-varying pressure. Specifically, at $t=t^n$ the pressure gradient $\Delta p(t^n)$ is 
\begin{equation}
\Delta p(t^n) = p_u(t^n) - p_d(t^n)\;,
\end{equation}
where values $p_u(t^n)$ and $p_d(t^n)$ are
\begin{equation}
p_u(t^n) := \begin{dcases} 
p_{M}^n \;, & {\text{lymphatic pressure at $i=M$, $t=t^n$}},\\
P_{in}(t^n)\;, & {\text{prescribed upstream pressure at $t=t^n$}},
\end{dcases}
\end{equation}
and 
\begin{equation}
p_d(t^n) := \begin{dcases} 
p_{1}^n \;, & {\text{lymphatic pressure at $i=1$, $t=t^n$}}, \\
P_{out}(t^n)\;, & {\text{prescribed downstream pressure at $t=t^n$}},
\end{dcases}
\end{equation}
where pressures $p_{M}^n$ and $p_1^n$ refers to the upstream and downstream lymphangions, respectively, and $P_{in}$ and $P_{out}$ are prescribed functions of time. 
The three possible configurations are summarized here
\begin{equation}
\Delta p(t^n) := \begin{dcases} 
P_{in}(t^n) - p_1^n \;, & {\text{leftmost valve}},\\
p_M^n - p_1^n \;, & {\text{valve between two lymphangions}},\\
p_M^n - P_{out}(t^n)\;, & {\text{rightmost valve}}.
\end{dcases}
\end{equation}
Once we numerically solve system \eqref{eq:systemValve}, the flow across the valve at the future time $q_v^{n+1}$ is determined.  

In the present paper, to find $A^*$ and calculate the numerical flux at the boundary we follow the numerical methodology proposed by Alastruey et al. \cite{Alastruey:2008b}. This method has already been used in \cite{Mueller:2014a, Mueller:2014b, Contarino:2016a, Sherwin:2003b}. 
To impose the valve flow at the boundaries of a lymphatic vessel, we use the Riemann invariant
\begin{equation}
\df{q_v^{n+1}}{A^*}-u^n+\beta\int_{A^n}^{A^*}\df{c\left(\tau\right)}{\tau}\mathrm{d}\tau=0\;,
\end{equation}
where $A^n$ and $u^n$ are the cross-sectional area and the velocity at the cell adjacent to the boundary at current time $t=t^n$, $q_v^{n+1}$ is the known flow rate across the valve and 
\begin{align}\label{eq:Auxillary}
\beta=\begin{dcases}
-1     & \text{downstream lymphangion}\;, \\
1 & \text{upstream lymphangion}\;.
\end{dcases}
\end{align}
To find $A^*$, we solve the following non-linear algebraic equation
\begin{equation}\label{eq:Non-linearF}
\mathcal{F}(A^*):=q_v^{n+1}+A^*\left(-u^n+\beta\int_{A^n}^{A^*}\df{c\left(\tau\right)}{\tau}\mathrm{d}\tau\right)=0\;,
\end{equation}
using the Newton-Raphson iterative method. Then the numerical flux at the boundary is
\begin{equation}\label{eq:FluxBoundary}
\F_{\f{1}{2} \text{ or } M+\f{1}{2}}=\F\left(\Q^*\right)\;,
\end{equation}
where
\begin{equation}
\Q^*=\colvec{A^*}{q_v^{n+1}}\;.
\end{equation}
When a valve is interposed between two lymphangions, then the non-linear problem \eqref{eq:Non-linearF} has to be solved twice: one for the upstream lymphangion ($\beta=1$) and one for the downstream lymphangion ($\beta=-1$).
The above numerical flux is a Godunov-type flux. Indeed, the numerical flux is evaluated at the solution along the $t$-axis of a classical Riemann problem where one of the two unknowns is given by an ODE, which in this case is the valve model, see \cite{Borsche:2016a} for similar works. 

Fig. \ref{fig:Coupling} illustrates the two non-linear problems to be solved to couple two lymphangions connected by a valve. First, we need to solve system \eqref{eq:systemValve}, and then the resulting flow rate $q_v$ is imposed at both lymphangions. Once we solve the two numerical problems, then the numerical fluxes at the boundaries are provided by \eqref{eq:FluxBoundary}.

\subsection{Imposed pressure at boundaries} \label{sec:imposedPressure}
In various experiments reported in the literature \cite{Davis:2011a, Davis:2012a, Scallan:2012a, Scallan:2013a}, lymphatic collectors containing two or more valves were isolated from the animal, cannulated at each end with a glass micropipette and pressurised. Here we present a numerical method to simulate this kind of experiments. The procedure is similar to the coupling method for valves and lymphangions: both cross-sectional area and flow rate need to be found at the ending interface where the pressure needs to be imposed. 

Consider a time-varying pressure $p_{I}\left(t\right)$ at a terminal interface. 
From pressure $p_I\left(t\right)$, cross-sectional area $A_{I}\left(t\right)$ can be calculated by using the inverse of the tube law \eqref{eq:pressure}. The flow rate $q^*$ can be found by applying the Riemann invariants as described in \ref{sec:couplingvalves}, and in this case it can be explicitly calculated as
\begin{equation}
q^*=A_{I}\left(t^n\right)\left(u^n-\beta\int_{A^n}^{A_I\left(t^n\right)}\df{c\left(\tau\right)}{\tau}\mathrm{d}\tau\right)\;,
\end{equation}
where $A^n$, $u^n$ and $\beta$ are the cross-sectional area and the velocity at the cell adjacent to the boundary at current time $t=t^n$ and $\beta$ is given by Eq. \eqref{eq:Auxillary}.
As before, the numerical flux at the boundary is
\begin{equation}
\F_{\f{1}{2} \text{ or } M+\f{1}{2}}=\F\left(\Q^*\right)\;,
\end{equation}
where
\begin{equation}
\Q^*=\colvec{A_I\left(t^n\right)}{q^*}\;.
\end{equation}

\subsection{Numerical methods for the systems of ODEs}
The systems of ODEs \eqref{eq:systemValve} and \eqref{eq:FHNSystem} were numerically solved with a second-order implicit Runge-Kutta method using the Lobatto IIIC method. The Butcher tableau is
\begin{center}
\begin{tabular}{c|cc}
0 & 1/2 & -1/2 \\
1 & 1/2 & 1/2 \\
\hline
&1/2&1/2
\end{tabular}
\end{center}
In the next section, we present the coupling of the systems of PDEs and ODEs, through an algorithm description.


\subsection{Complete algorithm}

Here we provide the complete algorithm to update the solution from time $t^n$ to time $t^{n+1}=t^n+\Delta t$. When not specified, the initial conditions are: $p(x,0)=P_{in}(0)$, $u(x,0)=0$, $v(0)=0.1$, $w(0)=s(0)=I(0)=0$ and $q_v(0)=\xi(0)=0$.

\begin{enumerate}
\item Assume data for all variables at $t=t^n$.
\item Calculate the time step $\Delta t$ as explained in Section \ref{sec:finitevolume}.
\item Evolve the valve flow $q$ and valve state $\xi$ of each lymphatic valve from time $t^n$ to $t^{n+1}$ by applying a second-order implicit Runge-Kutta method to the system of ODEs \eqref{eq:systemValve} and assuming the pressure difference $\Delta p$ at time $t^n$.
\item Calculate the numerical fluxes at the boundaries $\F_{\f{1}{2}}$ and $\F_{M+\f{1}{2}}$ of each lymphangion, as described in Sections \ref{sec:couplingvalves} and \ref{sec:imposedPressure}, using the lymphatic valve flow rates at time $t^{n+1}$.
\item Using the contraction state $s$ at the current time $t^n$, calculate the numerical fluxes $\F_{i+\f{1}{2}}$ within each domain of the lymphangions using the SLIC method (Section 14.5.3 of \cite{Toro:2009a}).
\item Using the contraction state $s$ at the current time $t^n$, calculate the numerical sources $\S_{i}$ within each domain of the lymphangions using a second-order method in space and time (Chapter 19 of \cite{Toro:2009a}).
\item Update the conserved variables $\Q$ of the PDEs of each lymphangion from time $t^n$ to $t^{n+1}$ according with finite volume formula \eqref{eq:FiniteVolumeFormula}.
\item Evolve the variables of the EFMC model of each lymphangion from time $t^n$ to $t^{n+1}$ by applying a second-order implicit Runge-Kutta method to the system of ODEs \eqref{eq:FHNSystem} and using the space-time averaged cross-sectional area and WSS at time $t^{n+1}$.
\end{enumerate}
The EFMC model and the system of PDEs are coupled through the contraction state $s$. The variable $s$ gives the actual value of the Young modulus in Eq. \eqref{eq:ContractionYoungModulus} to be used to calculate the physical flux in Eq. \eqref{eq:PhysicalFlux}.
Observe that even though we use second-order methods for every model, the accuracy of the global algorithm is formally of only first order. This is caused by the coupling methods. As a matter of fact, we couple the set of ODEs and PDEs using only use a first-order method. There are more sophisticated high-order coupling methods in the literature, see for instance \cite{Borsche:2016a}.

}
\section{Results and discussion}\label{sec:results}
In this section, we assemble the components of the model, implement the explained numerical methods and show numerical results. First, we numerically solve a Riemann problem for the PDEs and compare the numerical solution with the exact solution; then we show results for a single, three and ten lymphangions. We analyse lymphodynamical indexes for a wide range of boundary pressures, perform two sensitivity analyses for positive and negative pressure gradients and finally investigate the effect of defective lymphatic valves. Table \ref{table:parameters} gives parameters used in the numerical simulations.

\subsection{Test problem with piecewise initial condition: a Riemann problem} \label{sec:RP}
\begin{figure*}[htb]
\begin{center}
\includegraphics[width=0.8\textwidth]{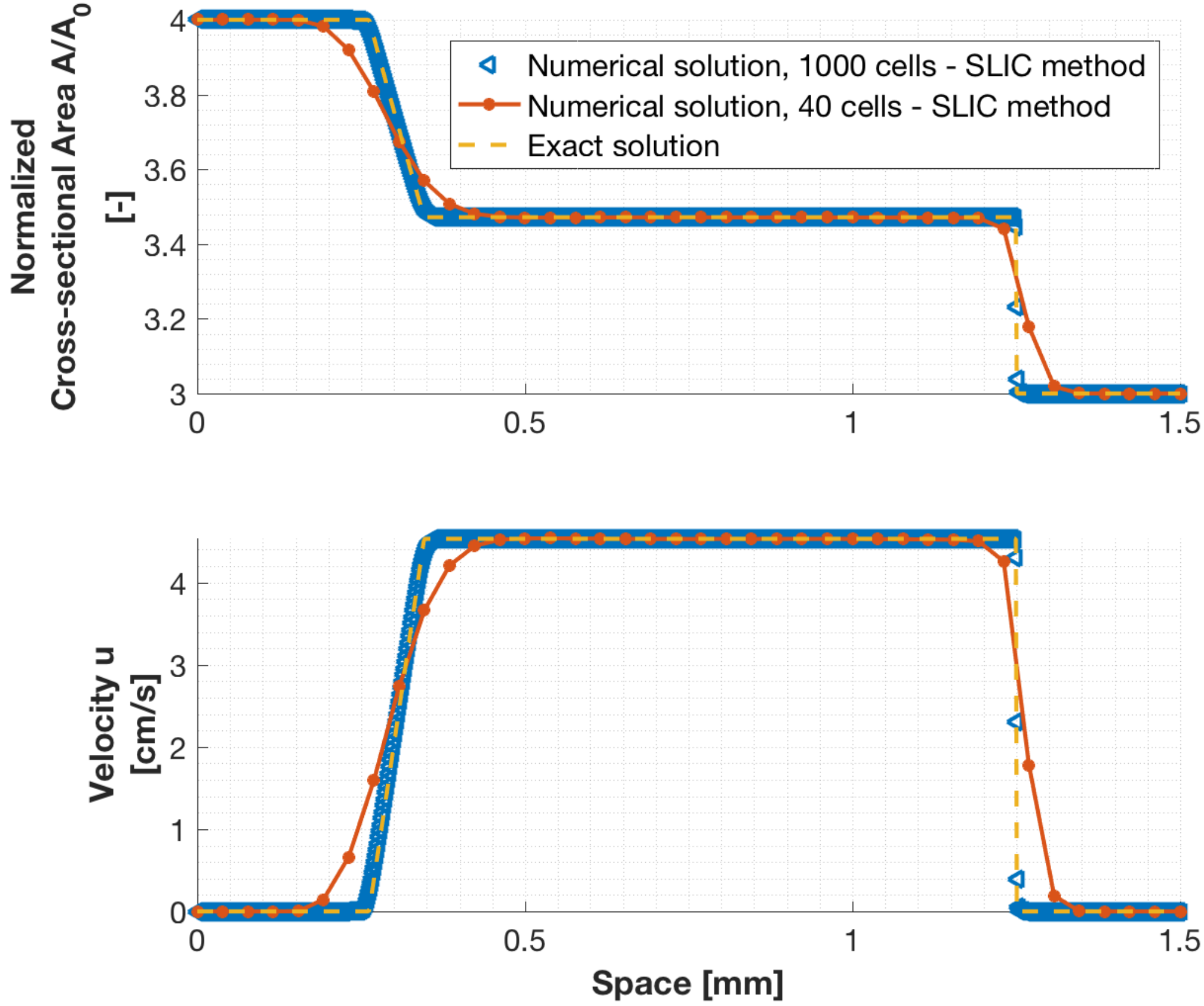}
\end{center}
\caption{\scriptsize {\bf Riemann problem for a single lymphangion without contractions.} Top and bottom frames show the following: normalised cross-sectional area and velocity at a fixed prescribed time. We compare the numerical results with the exact solution. The initial conditions are given in Section \ref{sec:RP} and the output time is $t_{output}=0.0015$ $s$. Here we used $M=40$ and $M=1000$ computational cells to discretize the lymphangion.} \label{fig:Test2}
\end{figure*}

Here we solve numerically a Riemann problem for the one-dimensional model for a single lymphangion without valves. The Riemann problem is a particular Cauchy problem where the initial conditions are piecewise constant with a single discontinuity. The exact solution in subcritical regime is available for this PDE system but is not reported here. For the exact solution of the Riemann problem for arteries and veins see \cite{Toro:2012a,Toro:2013a,Spiller:2017a}. Here, the chosen initial condition for $A$ and $u$ is
\begin{equation}
A(x,0)=\begin{cases}
A_L=4A_0\;, & x<\frac{L}{2}\;, \\
A_R=3A_0\;, & x>\frac{L}{2}\;, \\
\end{cases}
\end{equation}
where $L$ is the lymphangion length and $A_0$ is the cross-sectional area at equilibrium, and $u(x,0)=0$. We assumed no contractions, dynamic viscosity $\mu=0$ and transmissive boundary conditions. Fig. \ref{fig:Test2} shows the numerical results using $M=40$ and $M=1000$ cells and the exact solution at the output time time $t_{output}=0.0015$ $s$. The numerical solution with $M=40$ is comparable with the exact solution, which is composed of a left rarefaction and a right shock. The numerical solution with $M=1000$ confirms that the numerical solution converges to the exact solution. This result is typical of a $2\times 2$ non-linear system of hyperbolic differential equations and is comparable with the Riemann problem for the Euler equations, shallow water equations and the blood flow equations \cite{Toro:2009a}. 

\subsection{Four representative test cases of a collecting lymphatic}
Here we show four representative test cases where we simulated from one to ten lymphangions. In the first test case, we highlight the diastolic and systolic phases of the lymphatic cycle. The second test clearly shows the frequency increases as the intraluminal pressure increases. The third case shows the negative chronotropic effect given by the WSS. The four case shows an example of a bigger collecting lymphatic composed of ten lymphangions and eleven valves. Figs. \ref{fig:Test1}, \ref{fig:Test3} and \ref{fig:Test4} show the numerical results of cases from one to three. From top to bottom frames we show: an illustration of lymphangions and lymphatic valves, time-varying valve states (open $\xi=1$ and closed $\xi=0$), flow rates across the valves, and pressures, velocities, diameters, WSS, Young's modulus at the centre of the lymphangions, excitable variable and stimulus. The line colours shown from the second to the last panels refer the colour configuration shown in collector illustrated in the first panel. In the last panel, the {\em unstable spiral-node} region is represented by a blue-shaded area, while the {\em unstable node} region is represented by a red-shaded area. 

\subsubsection{Test case 1: representative case of a single lymphangion}

\begin{figure}[H]
\begin{center}
\includegraphics[width=0.45\textwidth,height=0.8\textheight,keepaspectratio]{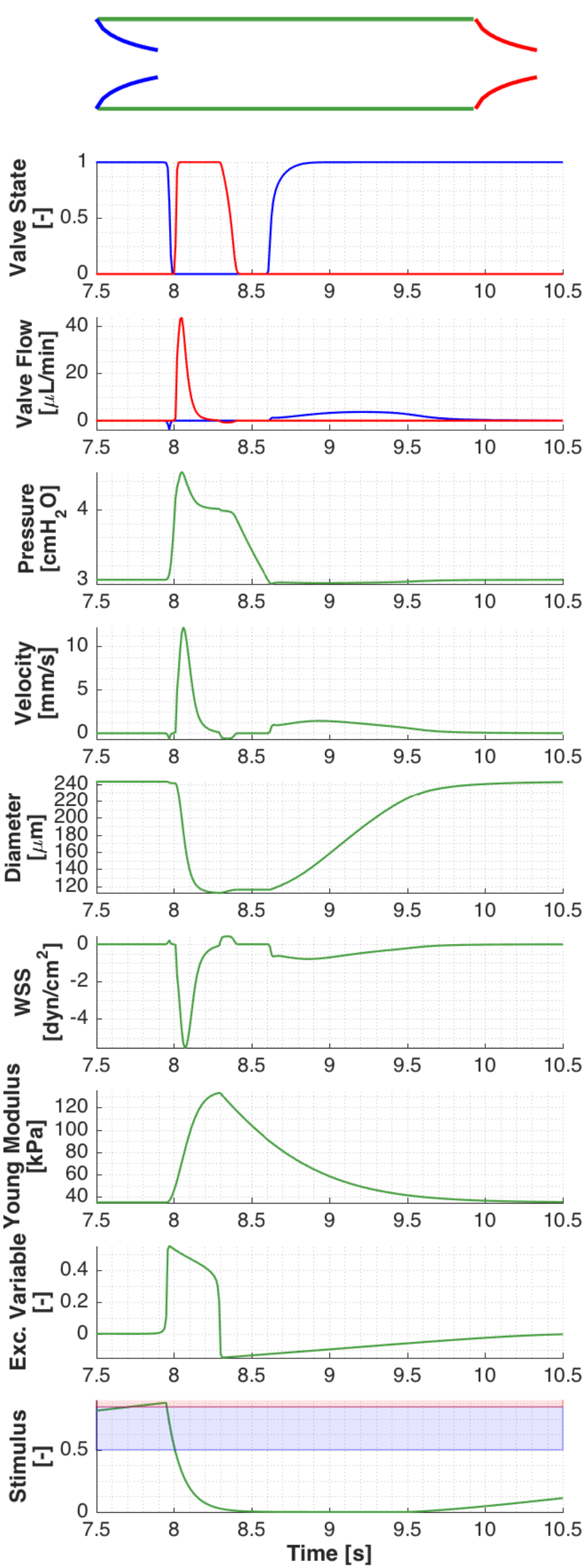}
\end{center}
\caption{\scriptsize 
{\bf Test case 1: representative case of a single lymphangion.} Boundary pressures: $P_{in}$ = 3 cmH$_2$O and $P_{out}$ = 4 cmH$_2$O. From the top to the bottom frames we show the following: illustration of lymphangions and lymphatic valves, time-varying valve states (open $\xi = 1$ and closed $\xi = 0$), flow rates across the valves, and pressures, velocities, diameters, WSS, Young modulus at the centre of the lymphangion, excitable variable excitable variables and stimulus. The colors shown from the second to the last panels refer the colour configuration shown in the first panel. In the last panel, blue and red shaded area illustrate the {\em unstable spiral-node} and the {\em unstable node} region, respectively. In this numerical test we used $M=20$ computational cells to discretize the lymphangion.} \label{fig:Test1}
\end{figure}

Here we show a representative test where the EFMC model, the valve model and the one-dimensional model for lymph flow are coupled together. We show results for a single lymphangion with two terminal valves. Upstream (left) and downstream (right) boundary pressures are $P_{in}=3$ cmH$_2$O and $P_{out}=4$ cmH$_2$O, respectively. We applied the numerical method at the boundaries explained in Section \ref{sec:couplingvalves}.
Fig. \ref{fig:Test1} shows the numerical results. As soon as the stimulus goes beyond the unstable spiral-node region and falls into the unstable node region, an action potential occurs followed by a contraction. The lymphatic pressure increases and the downstream valve opens: both flow rate across the downstream valve and the lymphatic velocity increase, the diameter decreases and the absolute value of the WSS increases. When the lymphatic pressure decreases below the downstream boundary pressure $P_{out}$, the downstream valve closes, and then the flow rate and the velocity return to zero. Then the lymphatic pressure decreases below the upstream boundary pressure $P_{in}$ and the upstream valve opens. Subsequently, the lymphangion is filled with lymph. 

\begin{figure*}
\subfloat{\label{fig:a}\includegraphics[width=0.5\textwidth]{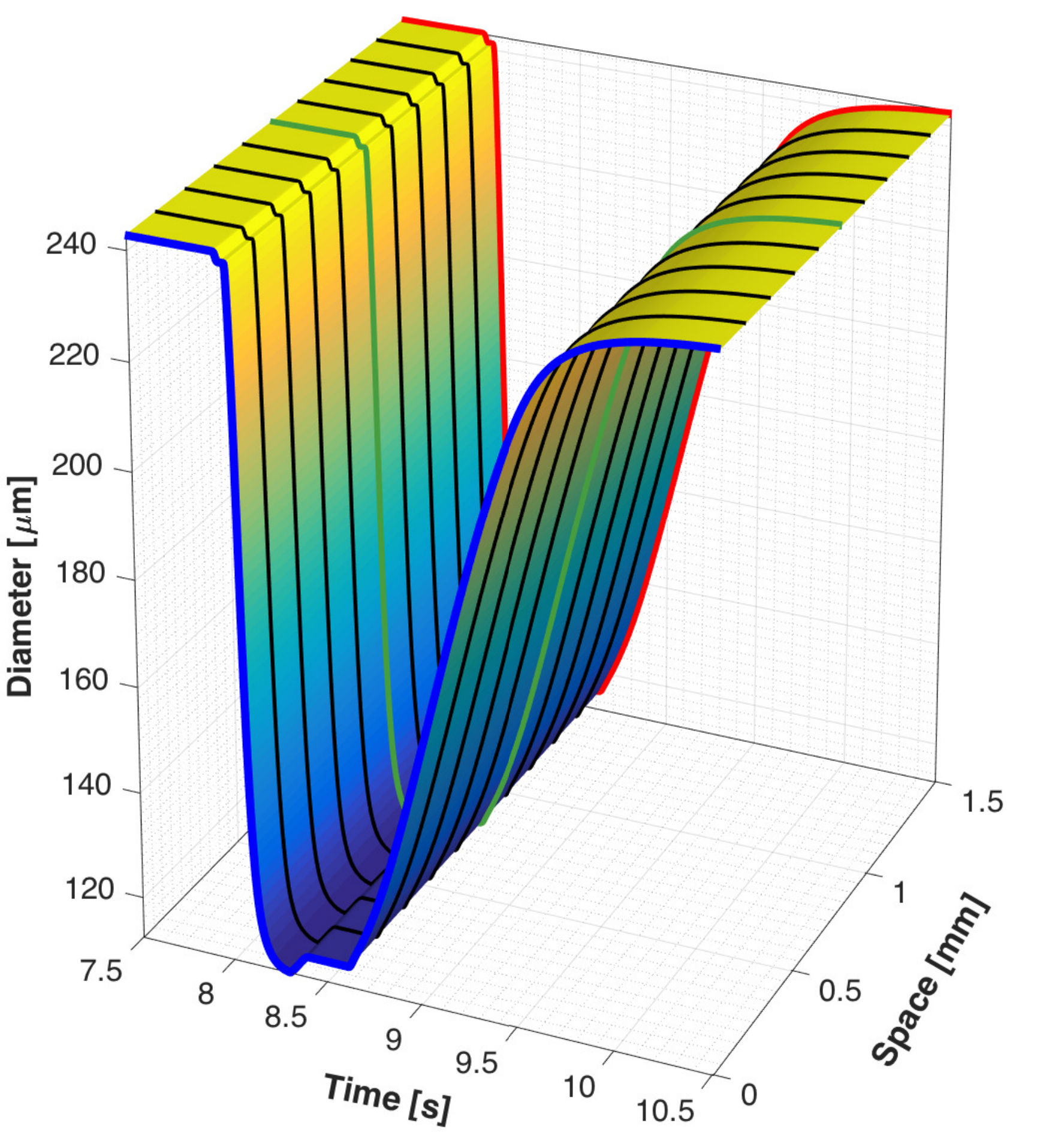}}
\subfloat{\label{fig:a}\includegraphics[width=0.5\textwidth]{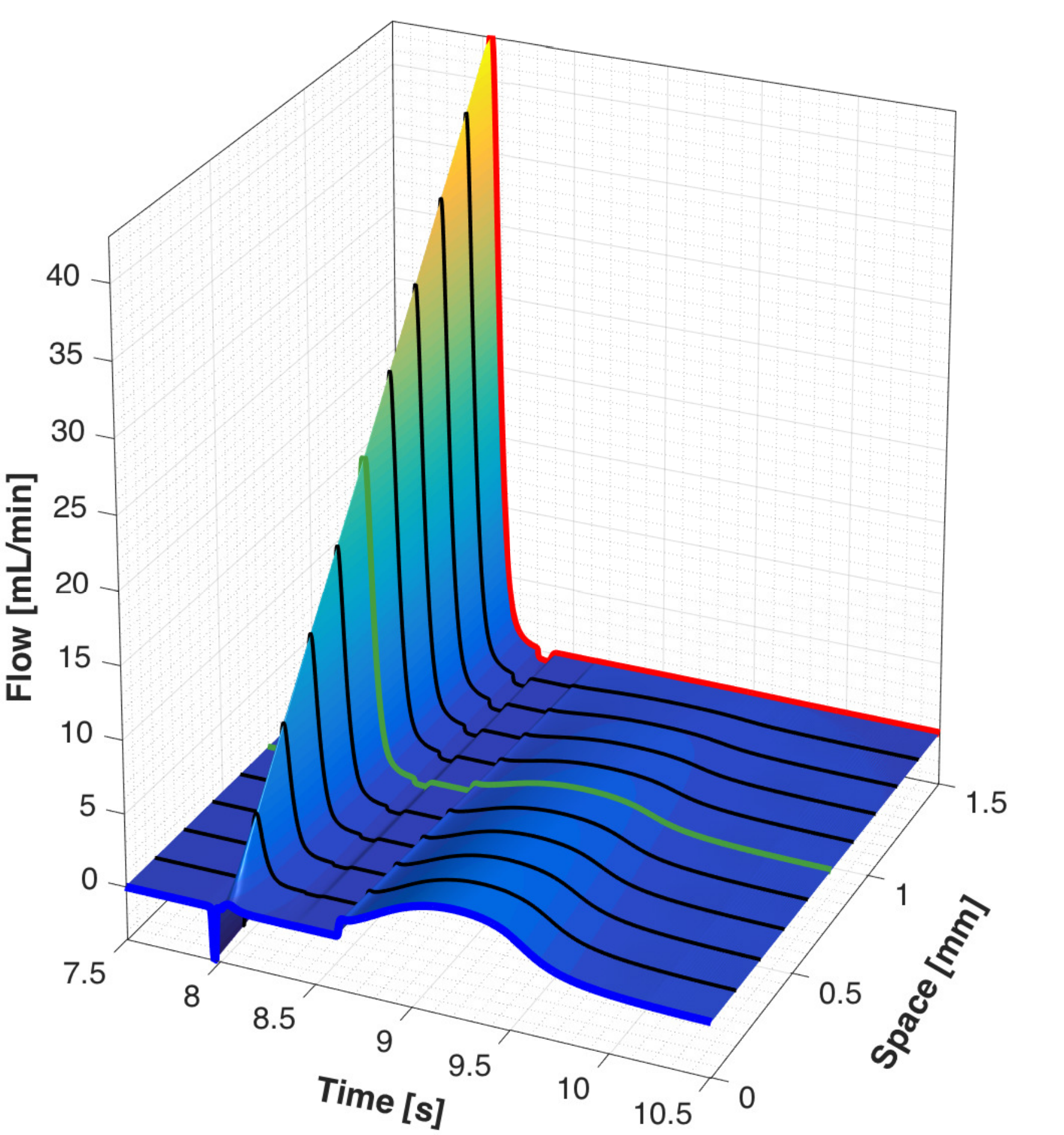}} \\
\subfloat{\label{fig:a}\includegraphics[width=0.5\textwidth]{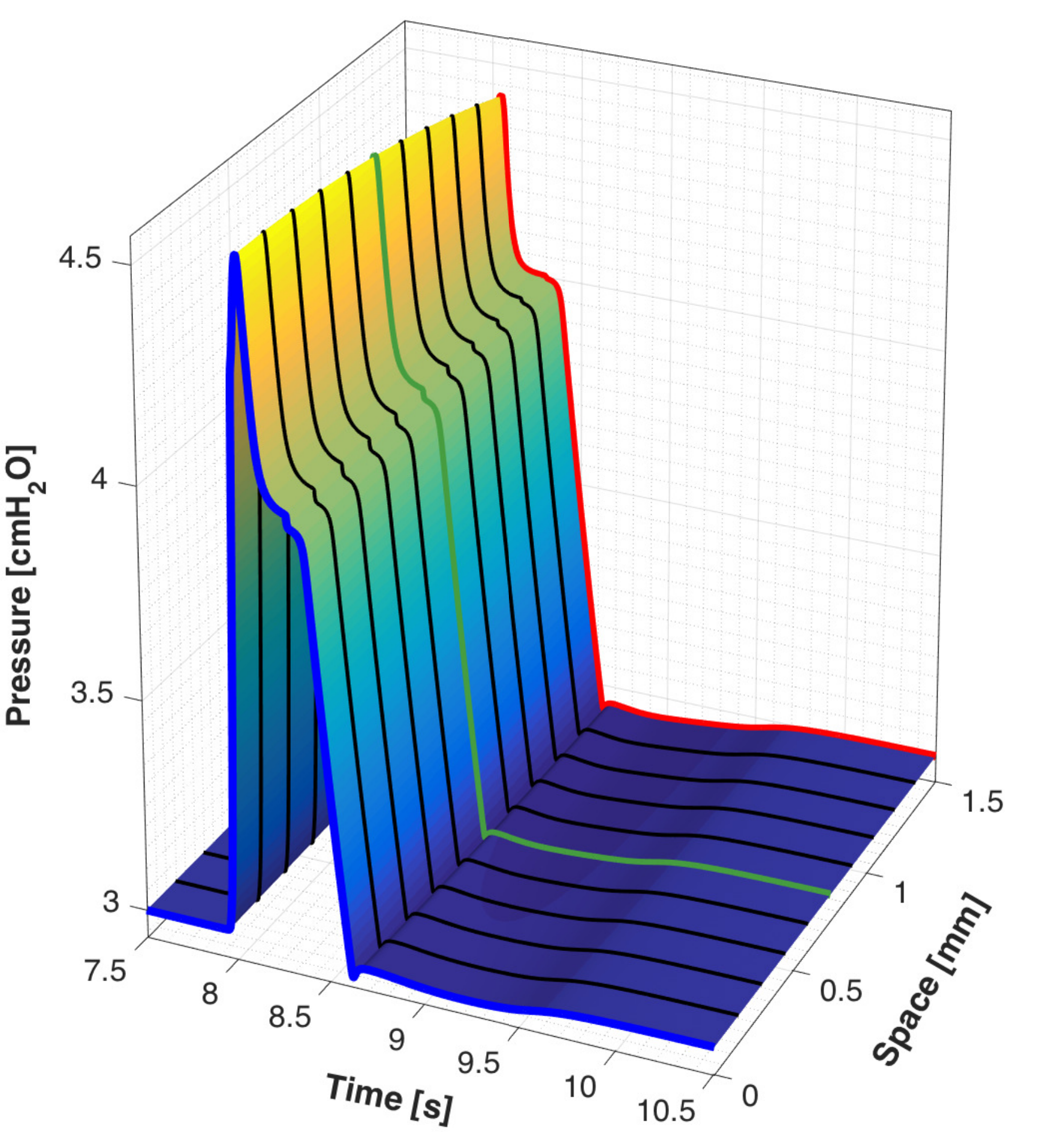}}
\subfloat{\label{fig:a}\includegraphics[width=0.5\textwidth]{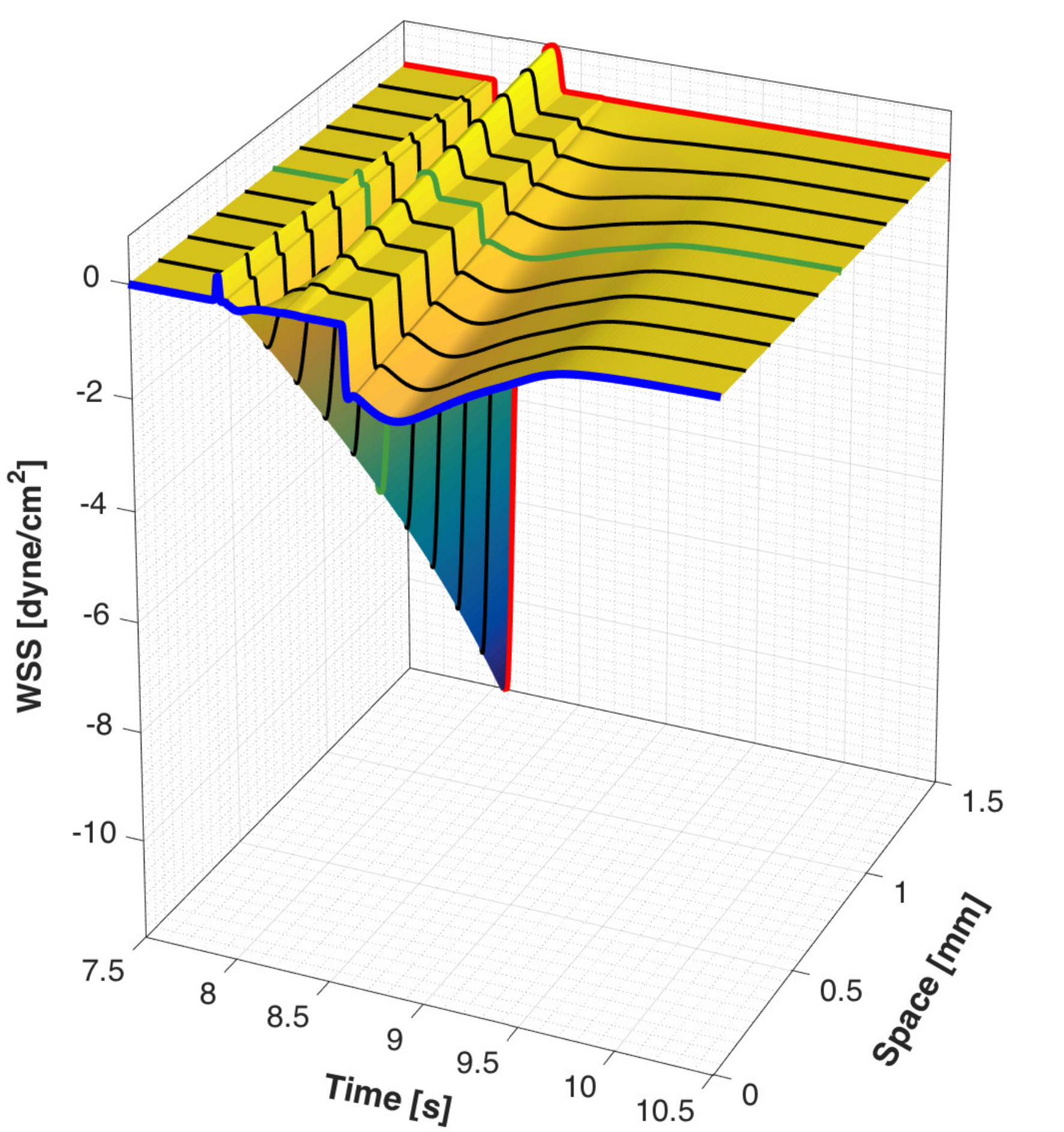}}
\caption{\scriptsize {\bf Test case 1: representative case of a single lymphangion (space-time).} Boundary pressures: $P_{in}$ = 3 cmH$_2$O and $P_{out}$ = 4 cmH$_2$O. Here we show numerical results in space and time of diameter, flow, pressure and WSS. Blue and red lines represent the numerical solutions close to the upstream and downstream valve, respectively. The green line represents the numerical solution at the centre of the lymphangion. In this numerical test we used $M=51$ computational cells to discretize the lymphangion.}  \label{fig:Test1_spacetime}
\end{figure*}

Fig. \ref{fig:Test1_spacetime} shows the space-time numerical solution of diameter, flow, pressure and WSS. The blue and red line represent the numerical solutions close to the upstream and downstream valve, respectively. The green line represents the numerical solution at the centre of the lymphangion. 
The diameter decreases almost homogeneously throughout the lymphangion. The same happens to the pressure, but it can be noticed that the systolic pressure varies in space, as the maximum is reached at the upstream side, while the minimum is reached at the downstream one. The flow rate and the WSS share a similar behaviour with opposite signs. During the systolic phase, the flow rate (WSS) reaches its maximum (minimum) at the downstream side, while it reaches its minimum (maximum) at the upstream one. The red and blue lines in the flow rate are similar to the valve flow rates shown in Fig. \ref{fig:Test1}. This result highlights that the mathematical model can give quantitative information throughout the length of the lymphangion.

\subsubsection{Test case 2: contraction frequency increases as the intraluminal pressure increases} 
The numerical test shown here was inspired by the experiments performed in several works \cite{Davis:2011a, Davis:2012a, Scallan:2012a, Scallan:2013a} where time-varying pressures were imposed at the boundaries of the collector. More specifically, this test imitates the {\em ramp-wise $P_{out}$ elevation} shown in \cite{Davis:2012a}. We simulated a collector  
\begin{figure*}[h]
\begin{center}
\includegraphics[width=0.95\textwidth]{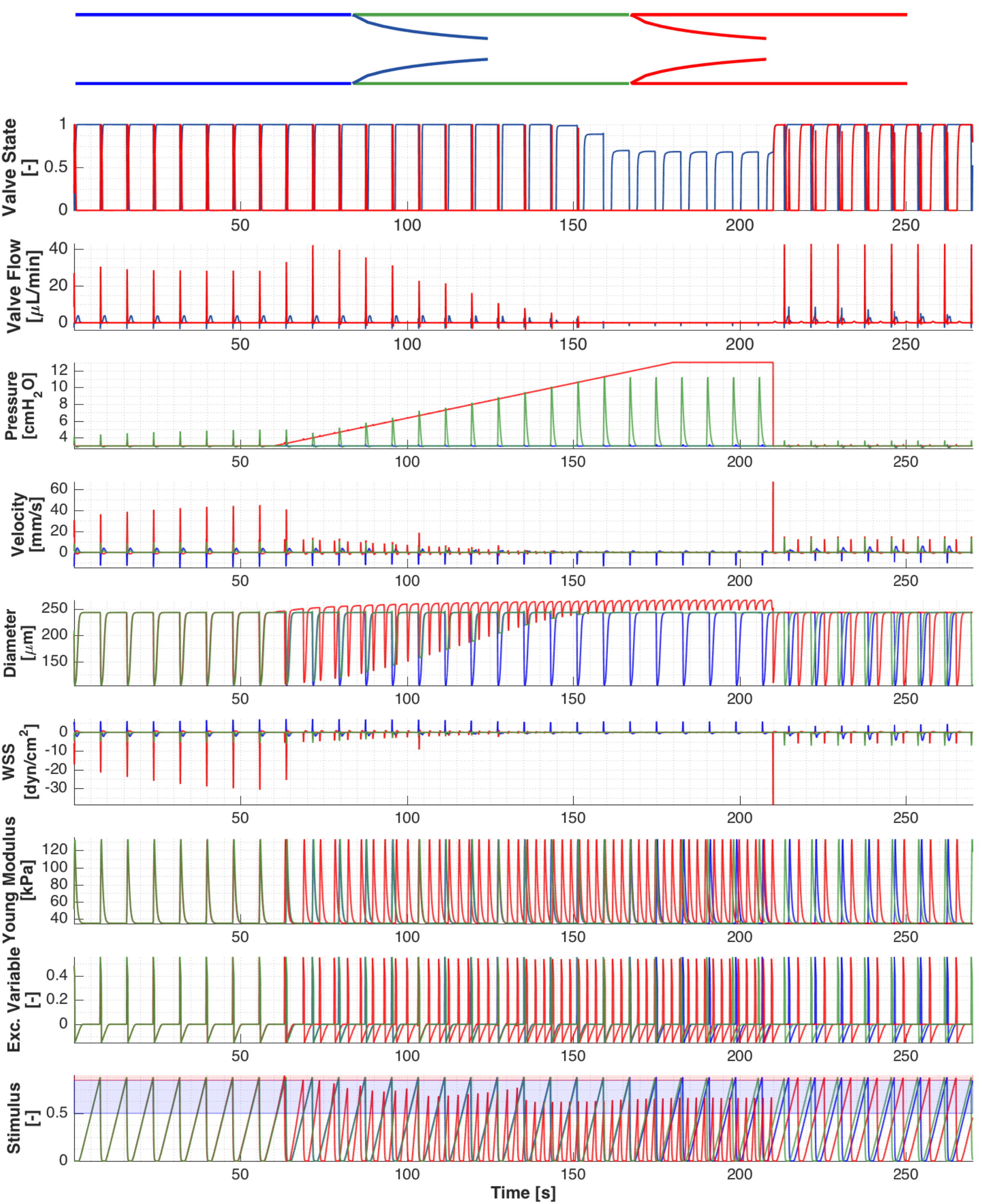}
\end{center}
\caption{\scriptsize 
{\bf Test case 2: the contraction frequency increases as the intraluminal pressure increases.}
Time-varying boundary pressures can be found in Eq. \eqref{eq:Test3boundary}. From the top to the bottom frames we show the following: illustration of lymphangions and lymphatic valves, time-varying valve states (open $\xi = 1$ and closed $\xi = 0$ ), flow rates across the valves, and pressures, velocities, diameters, WSS, Young modulus at the centre of the lymphangions, excitable variables and stimulus. The colors shown from the second to the last panels refer the colour configuration shown in the first panel. In the last panel, blue and red shaded area illustrate the {\em unstable spiral-node} and the {\em unstable node} region, respectively. In this numerical test we used $M=20$ computational cells to discretize each lymphangions.} \label{fig:Test3}
\end{figure*}
\afterpage{\clearpage}
\begin{figure*}[h]
\begin{center}
\includegraphics[width=0.95\textwidth]{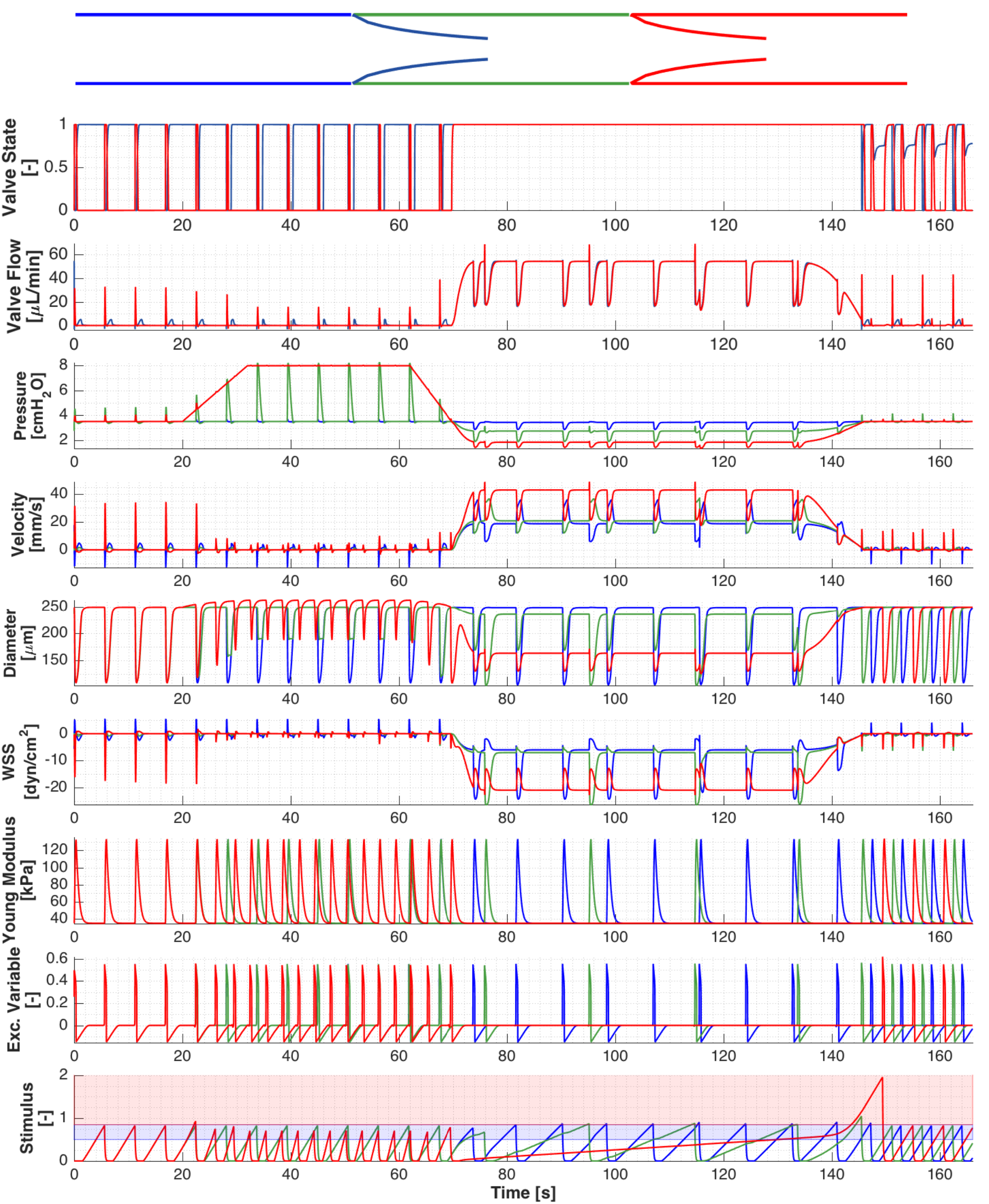}
\end{center}
\caption{\scriptsize {\bf Test case 3: contraction frequency decreases with increasing WSS.} Time-varying boundary pressures can be found in Eq. \eqref{eq:Test4boundary}. From the top to the bottom frames we show the following: illustration of lymphangions and lymphatic valves, time-varying valve states (open $\xi = 1$ and closed $\xi = 0$ ), flow rates across the valves, and pressures, velocities, diameters, WSS, Young modulus at the centre of the lymphangions, excitable variables and stimulus. The colors shown from the second to the last panels refer the colour configuration shown in the first panel. In the last panel, blue and red shaded area illustrate the {\em unstable spiral-node} and the {\em unstable node} region, respectively. In this numerical test we used $M=20$ computational cells to discretize each lymphangions.} \label{fig:Test4}
\end{figure*}
\afterpage{\clearpage} 
of three lymphangions and we imposed the following time-varying pressures
\begin{align}\label{eq:Test3boundary}
P_{out}(t)&=\begin{cases}
p_1\;, & t<t_1\;,\\
\f{p_2-p_1}{t_2-t_1}\left(t-t_2\right)+p_2\;, & t_1<t<t_2\;,\\
p_2\;, & t_2<t<t_3\;,\\
p_3\;, & t_3<t<t_{output}\;,\\
\end{cases} \\
P_{in}(t)&=p_1\;,
\end{align}
where $p_1=p_3=3$ cmH$_2$O, $p_2=13$ cmH$_2$O, $t_1=60$ s, $t_2=t_1+120$ s, $t_3=t_2+30$ s and $t_{output}=t_3+60$ s. 
Applying the numerical methods explained in \ref{sec:imposedPressure}, the inlet pressure $P_{in}$ was imposed at the leftmost interface of the upstream lymphangion, while the output pressure was imposed at the rightmost interface of the downstream lymphangion. Only a negative transaxial-pressure gradient is taken into account. 

Fig. \ref{fig:Test3} shows the results of the numerical simulation.
Even though the upstream and downstream lymphangions contract, their pressures are controlled. The downstream pressure (red line) follows the behaviour of the imposed output pressure $P_{out}$, while the upstream pressure (blue line) is almost constantly $P_{in}$. The lymphangion at the centre responds to these changes of boundary pressures (green line). Initially, both valves close and open, but when the downstream pressure $P_{out}$ reaches a certain value ($\approx$ 11 cmH$_2$O), the centred lymphangion cannot open the downstream valve anymore, and this can be seen by the valve state and the valve flow of the downstream lymphangion. Since the pressure of the downstream lymphangion increases during the numerical simulation, its frequency of contraction increases as well. As soon as $P_{out}$ falls back to 3 cmH$_2$O, then the contraction frequency suddenly decreases to the initial value. The frequency of centred lymphangion is not affected by the increase in the frequency of the downstream lymphangion. This comes from the fact that the current model does not take into account the interaction between adjacent lymphangions. From this numerical results, we can see that: 1) the frequency of contractions depends on the intraluminal pressure and 2) the centred lymphangion tries to overcome the downstream time-varying pressure by increasing the end-systolic pressure, but this is possible up to a certain threshold.

\begin{figure*}[t]
\begin{center}
\includegraphics[width=1\textwidth]{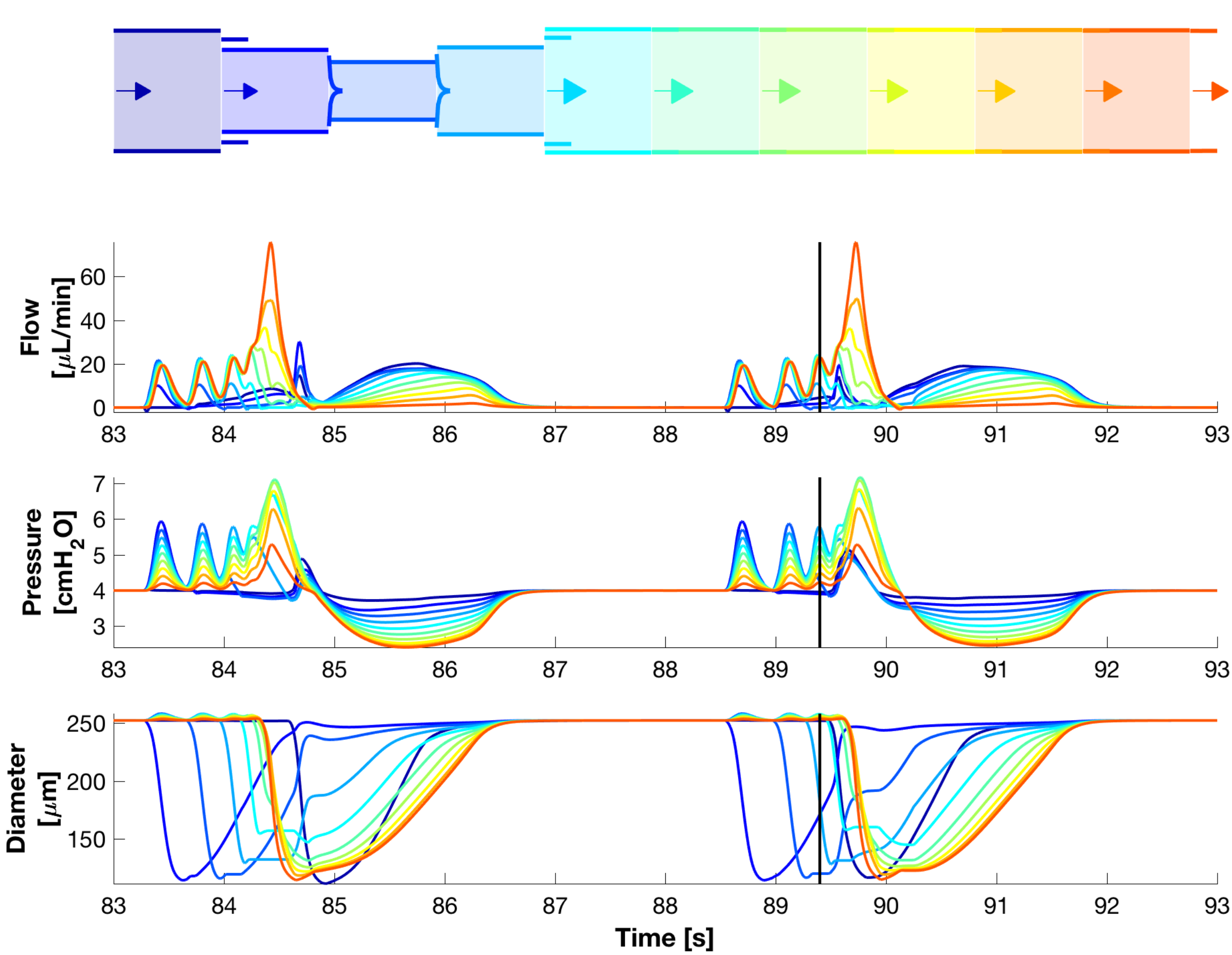}
\end{center}
\caption{\scriptsize {\bf Test case 4: representative case of a collector of ten lymphangions.} Boundary pressures: $P_{in}$ = 4 cmH$_2$O and $P_{out}$ = 4 cmH$_2$O. From the top to the bottom frames we show the following: an illustration of the collector at time $t=89.4$ (shown by the black vertical line), time-varying flow rate, pressure and diameter at the centre of each lymphangion. The colors shown from the second to the last panels refer the colour configuration shown in the first panel. In this numerical test we used $M=20$ computational cells to discretize each lymphangions.} \label{fig:TestTenLymphangion}
\end{figure*}

\subsubsection{Test case 3: contraction frequency decreases with increasing WSS}
The test proposed here simulates a collector of three lymphangions with two valves and highlights the effect of the WSS on the frequency of contractions.
As done for the test case 2, we imposed the following time-varying pressures at the terminal interfaces of the collector
\begin{align}\label{eq:Test4boundary}
P_{out}(t)&=\begin{cases}
p_1\;, & t<t_1\;,\\
\f{p_2-p_1}{t_2-t_1}\left(t-t_2\right)+p_2\;, & t_1<t<t_2\;,\\
p_2\;, & t_2<t<t_3\;,\\
\f{p_3-p_2}{t_4-t_3}\left(t-t_4\right)+p_3\;, & t_3<t<t_4\;,\\
p_3\;, & t_4<t<t_5\;,\\
\f{p_4-p_3}{t_6-t_5}\left(t-t_6\right)+p_4\;, & t_5<t<t_6\;,\\
p_4\;, & t_6<t<t_{output}\;,\\
\end{cases} \\
P_{in}(t)&=p_1\;,
\end{align}
where $p_1=p_4=3.5$ cmH$_2$O, $p_2=8$ cmH$_2$O, $p_3=1$ cmH$_2$O, $t_1=20$ s, $t_2=t_1+12$ s, $t_3=t_2+30$ s, $t_4=t_3+12$ s, $t_5=t_4+60$ s, $t_6=t_5+12$ s and $t_{output}=t_6+20$ s. 
In this test, both negative and positive transaxial pressure gradients are taken into account. The inlet pressure $P_{in}$ is fixed to $3.5$ cmH$_2$O, while the output pressure $P_{out}$ is initially equal to the inlet pressure, it increases up to $8$ cmH$_2$O, it decreases to $1$ cmH$_2$O, and finally returns back to the initial value. 

Fig. \ref{fig:Test4} shows the numerical results. The downstream lymphangion initially contracts at the same frequency of the other lymphangions, then as the output pressure $P_{out}$ increases, it contracts faster. The end-systolic pressure of the centre lymphangion follows the behaviour of the downstream pressure $P_{out}$. When a favourable pressure gradient occurs ($P_{in}>P_{out}$), then the absolute value of the WSS of each lymphangion increases, and thus the frequencies drastically decrease. Interestingly, when one of the lymphangions contracts, then the flow at the downstream valve decreases from $\approx$ 50 $\mu$L min$^{-1}$ to $\approx$ 18 $\mu$L min$^{-1}$. This means that the contractions of the lymphangions decrease the outflow, and therefore the negative chronotropic effect given by the increment of the WSS gives an overall positive effect on the averaged outflow. 

\subsubsection{Test case 4: representative case of a collecting lymphatic composed of ten lymphangions}
Here we illustrate a representative example of a collecting lymphatic composed of ten lymphangions and eleven valves. We applied the numerical method at the boundaries explained in Section \ref{sec:couplingvalves} with prescribed upstream and downstream pressures $P_{in}=P_{out}=4$ cmH$_2$O. We used $M=20$ computational cells to discretize each one-dimensional lymph vessels and we set the output time $t_{output}=93$ $s$.

Fig. \ref{fig:TestTenLymphangion} shows the numerical results of the simulation. Figure shows an illustration of the collector at time $t=89.4$ (shown by the black vertical line), and time-variation of flow rate, pressure and diameter at the centre of each lymphangion. Lymphatic contractions occur without prescribed delays between lymphangions. When the leftmost lymphangion contracts, the adjacent lymphangion becomes stimulated as its diameter increases and will thus contract as well. This leads to a chain of consecutive contractions, which occur in a non-linear manner. 

\subsection{Pressure versus normalised cross-sectional area (PA) plots for a single lymphangion }

\begin{figure}[tb]
\begin{center}
\includegraphics[width=0.5\textwidth]{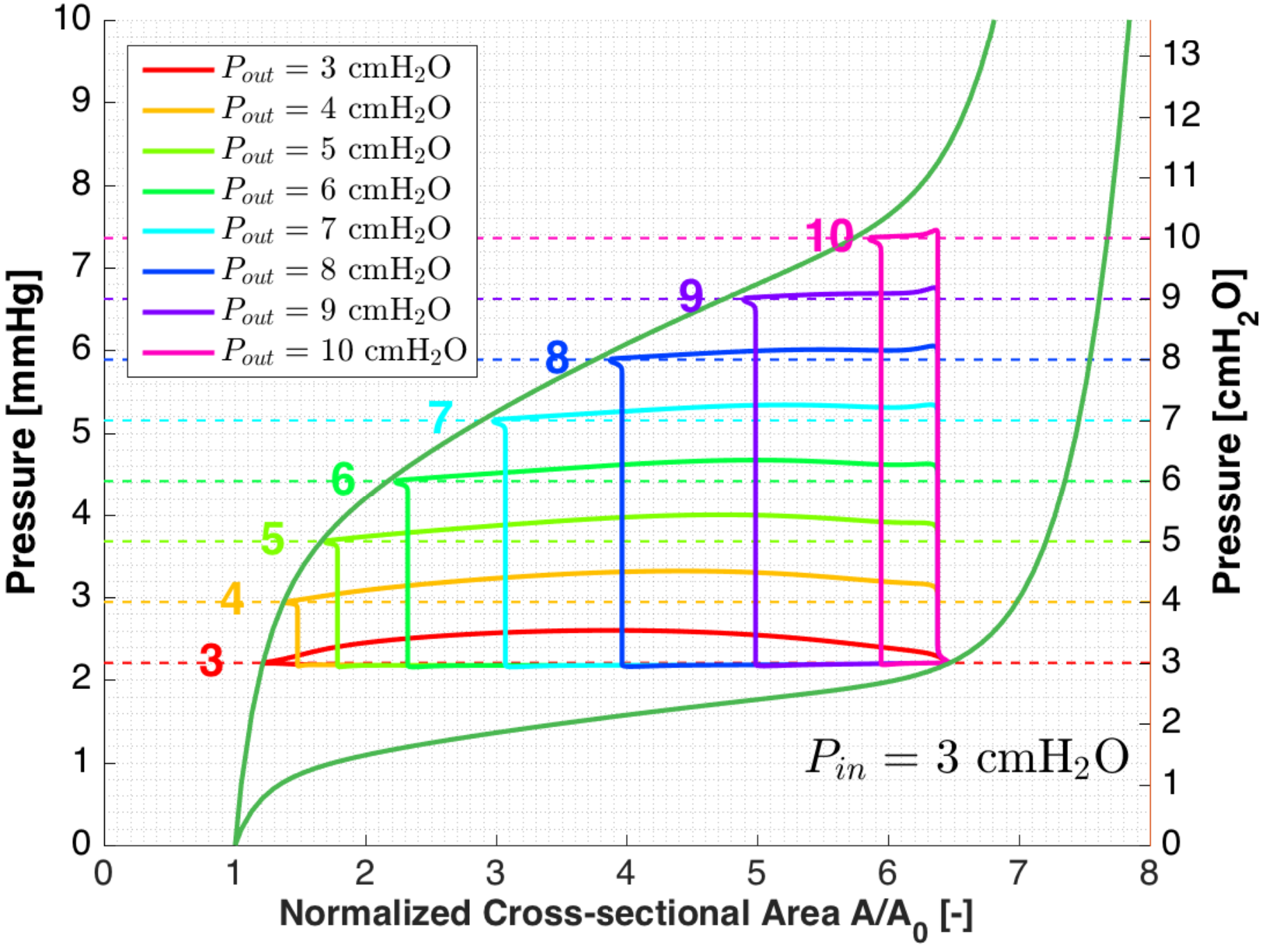}
\end{center}
\caption{\scriptsize {\bf Pressure against normalised cross-sectional area (PA) plots of during lymphatic contractions}. Here we simulated a single lymphangion with two lymphatic valves and different downstream pressures $P_{out}$ from 3 to 10 cmH$_2$O, while keeping fixed the upstream pressure $P_{in}$ to 3 cmH$_2$O. Pressure results are shown in mmHg and cmH$_2$O. The figure also shows the tube laws with and without contraction. Pressures and diameters were calculated at the centre of the lymphangion. \label{fig:pIn_pOut}}
\end{figure}
In this section we show plots of pressure against normalised cross-sectional area (PA) of a single lymphangion with fixed upstream pressure $P_{in}=3$ cmH$_2$O and different downstream pressures $P_{out}$, from $3$ to $10$ cmH$_2$O. The aim of this exercise is to show that the numerical results of the mathematical model imitate the experimental measurements of the pressure-volume relationship \cite{Davis:2012a}. Fig. \ref{fig:pIn_pOut} shows the numerical results and also the contracted and relaxed tube laws (green lines). Results show a qualitatively good agreement with \cite{Davis:2012a,Scallan:2012a}. As the downstream pressure $P_{out}$ increases, the PA plots shrink and the systolic pressure increases. As a consequence, the stroke work and the ejection fraction decrease towards zero. 
The systolic pressure can increase up to a certain level, depending on the baseline pressure. As a matter of fact, in the current case $P_{in}=3$ cmH$_2$O, the maximum systolic pressure is $\approx 11$ cmH$_2$O and decreases as $P_{in}$ decreases. For instance, for $A/A_0=4$ e $P_{in}\approx 2$ cmH$_2$O, the maximum systolic reachable pressure is $\approx$ 8 cmH$_2$O. From the literature, we know that lymphangions can increase the strength of contraction under stresses \cite{Davis:2012a}, which means that the slope of the {\em Ending Systolic Pressure-Volume Relationship} (ESPVR) increases \cite{Scallan:2012a} or, in analogy to our terminology, that the maximum Young's modulus $E_{max}$ increases somewhat. For an example of a mathematical model with adaptation of contractility, see Caulk et al. \cite{Caulk:2016a}. This represents a limitation of our current mathematical model as the tube laws do not undergo modification under continuous stresses.

\subsection{Analysis of lymphodynamical indexes by varying $P_{in}$ and $P_{out}$} \label{sec:pIn_pOut}

\begin{figure*}
\subfloat{\label{fig:a}\includegraphics[width=0.33\textwidth]{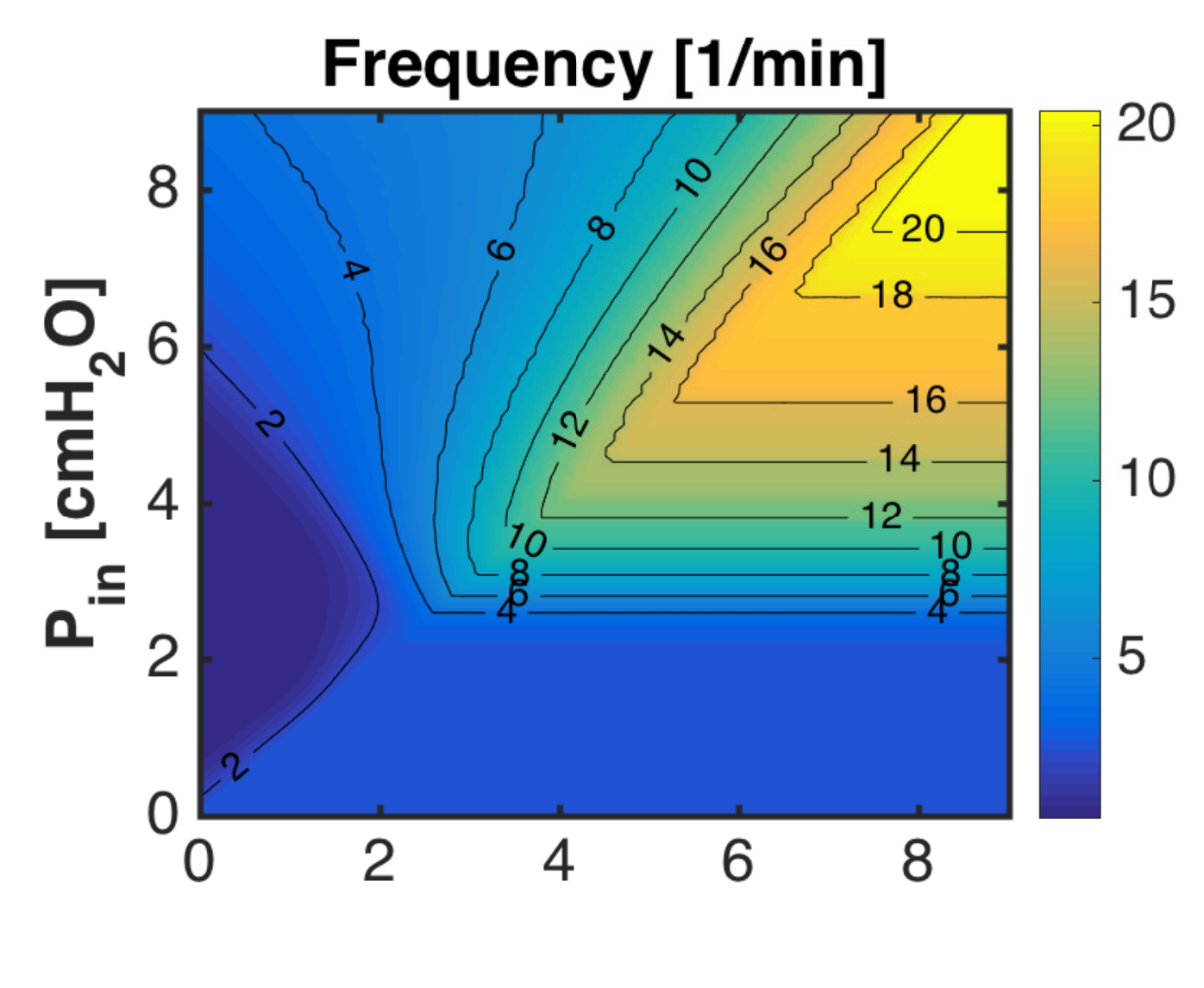}}
\subfloat{\label{fig:a}\includegraphics[width=0.33\textwidth]{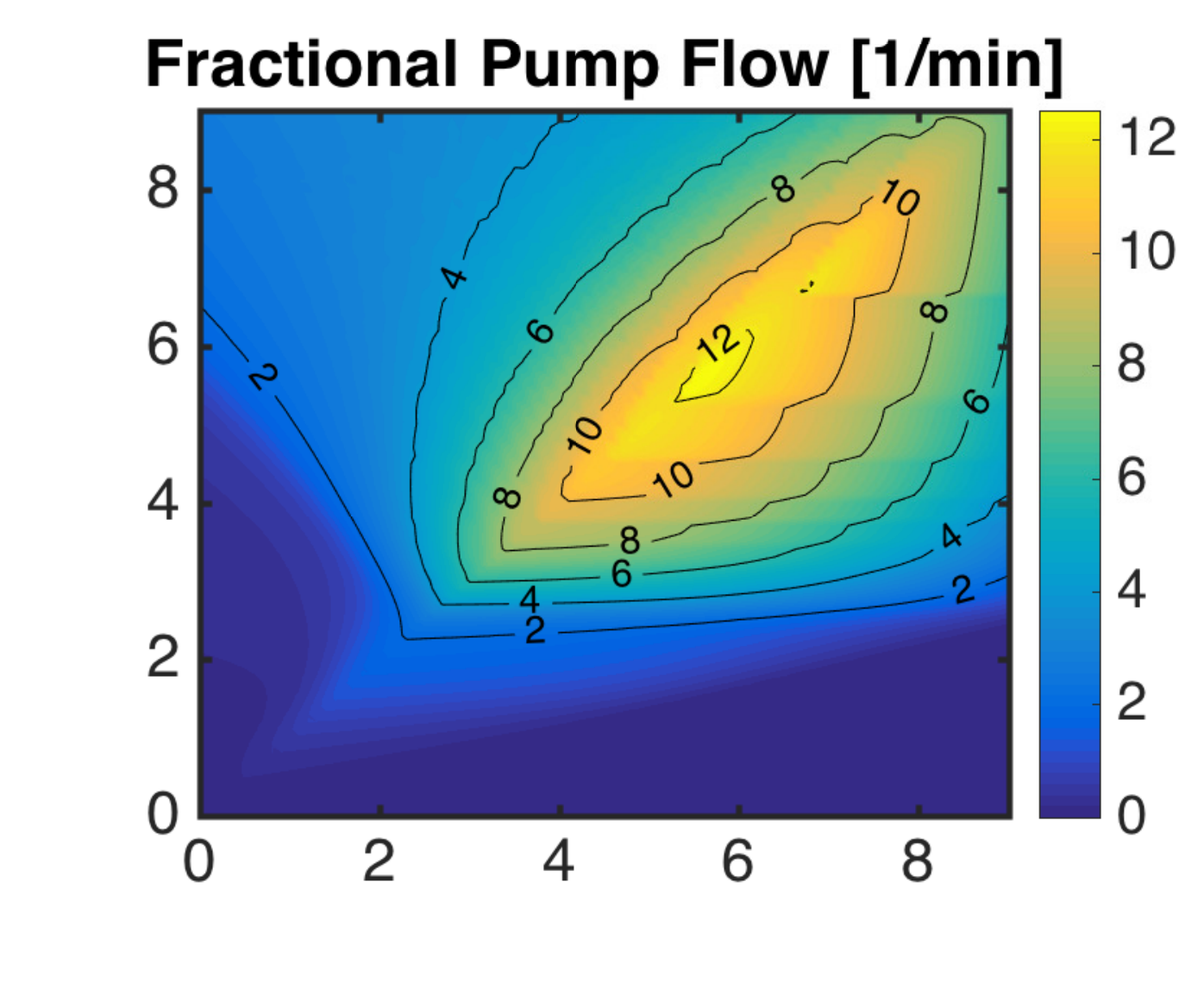}}
\subfloat{\label{fig:a}\includegraphics[width=0.33\textwidth]{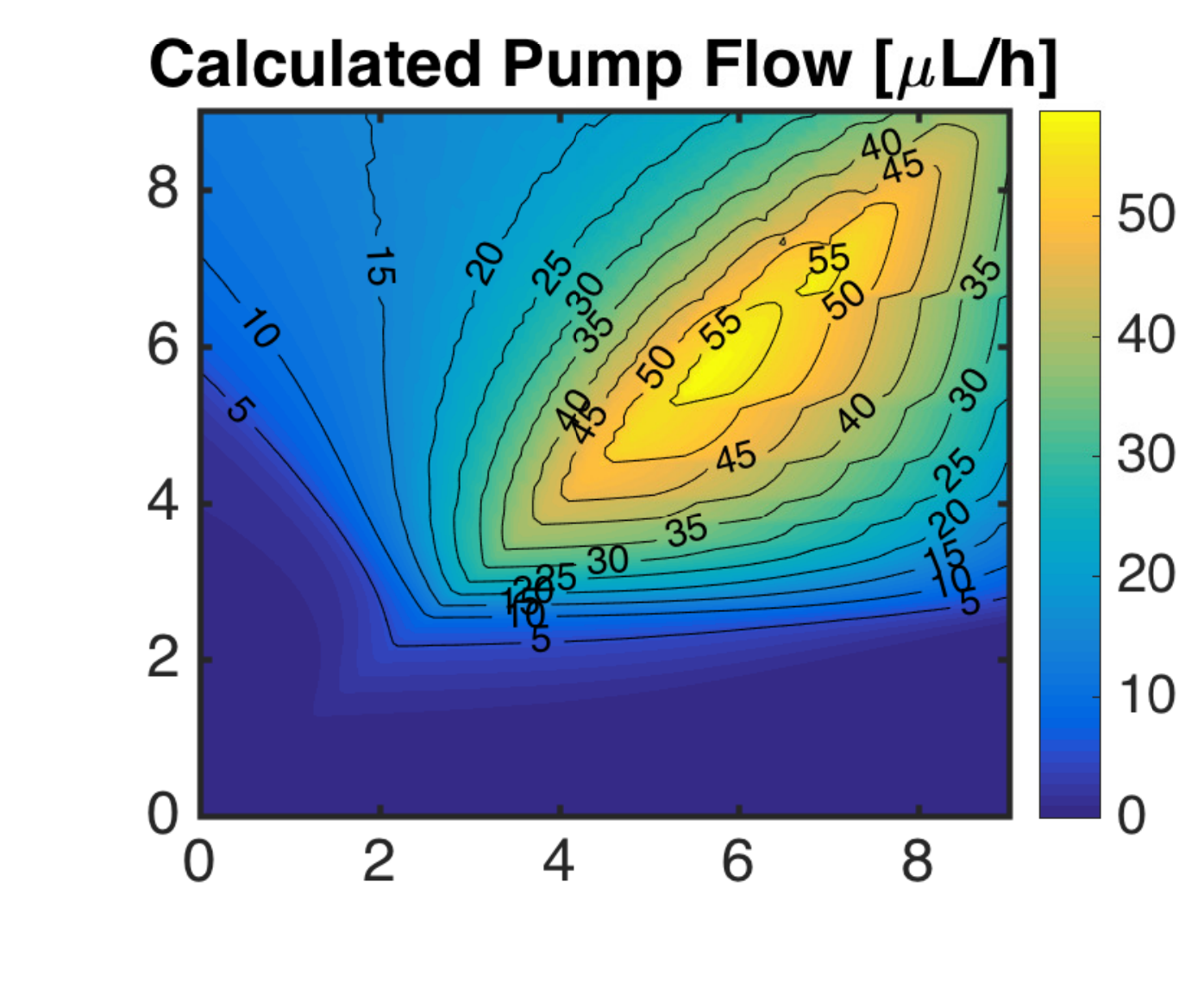}} \\
\subfloat{\label{fig:a}\includegraphics[width=0.33\textwidth]{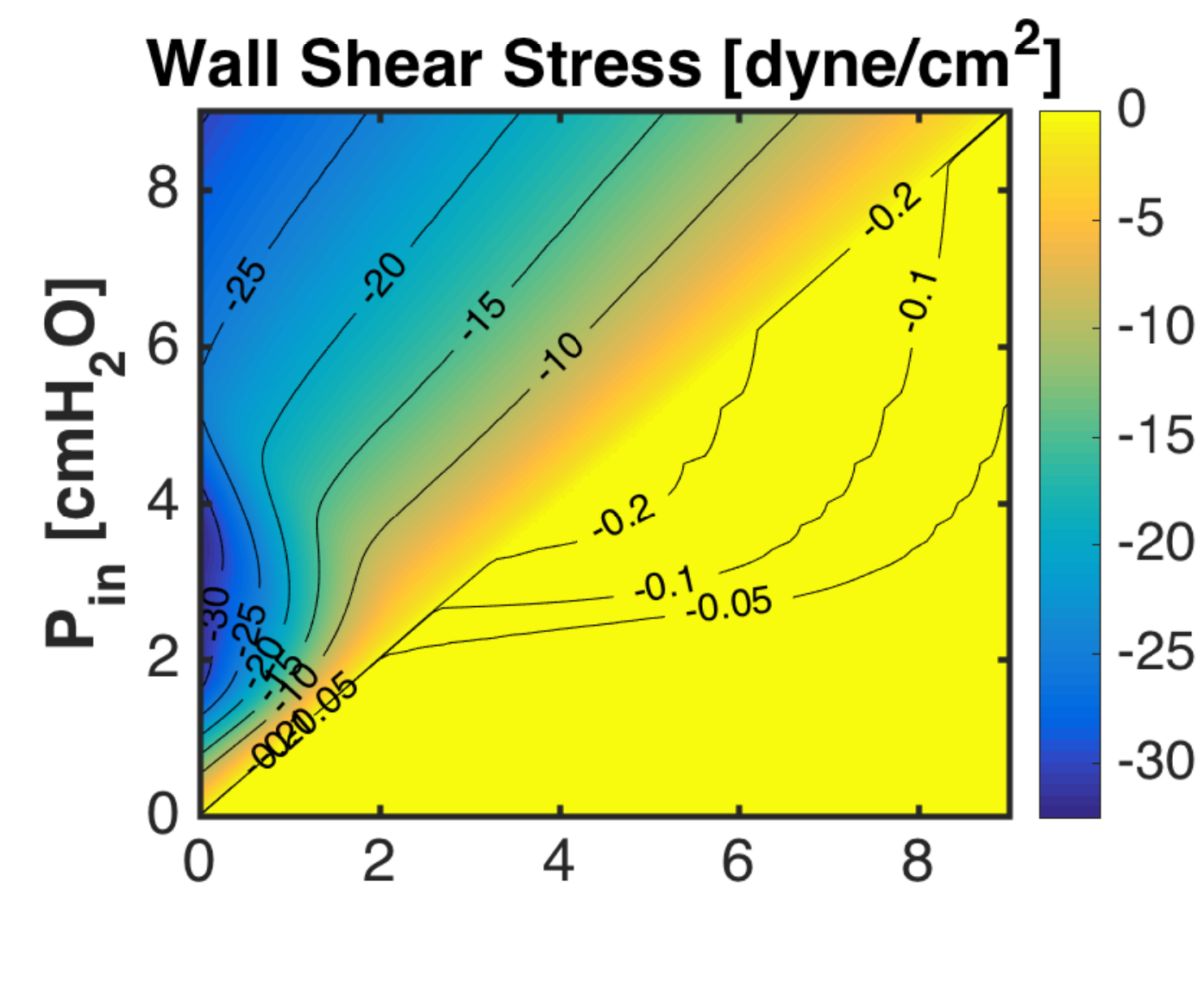}}
\subfloat{\label{fig:a}\includegraphics[width=0.33\textwidth]{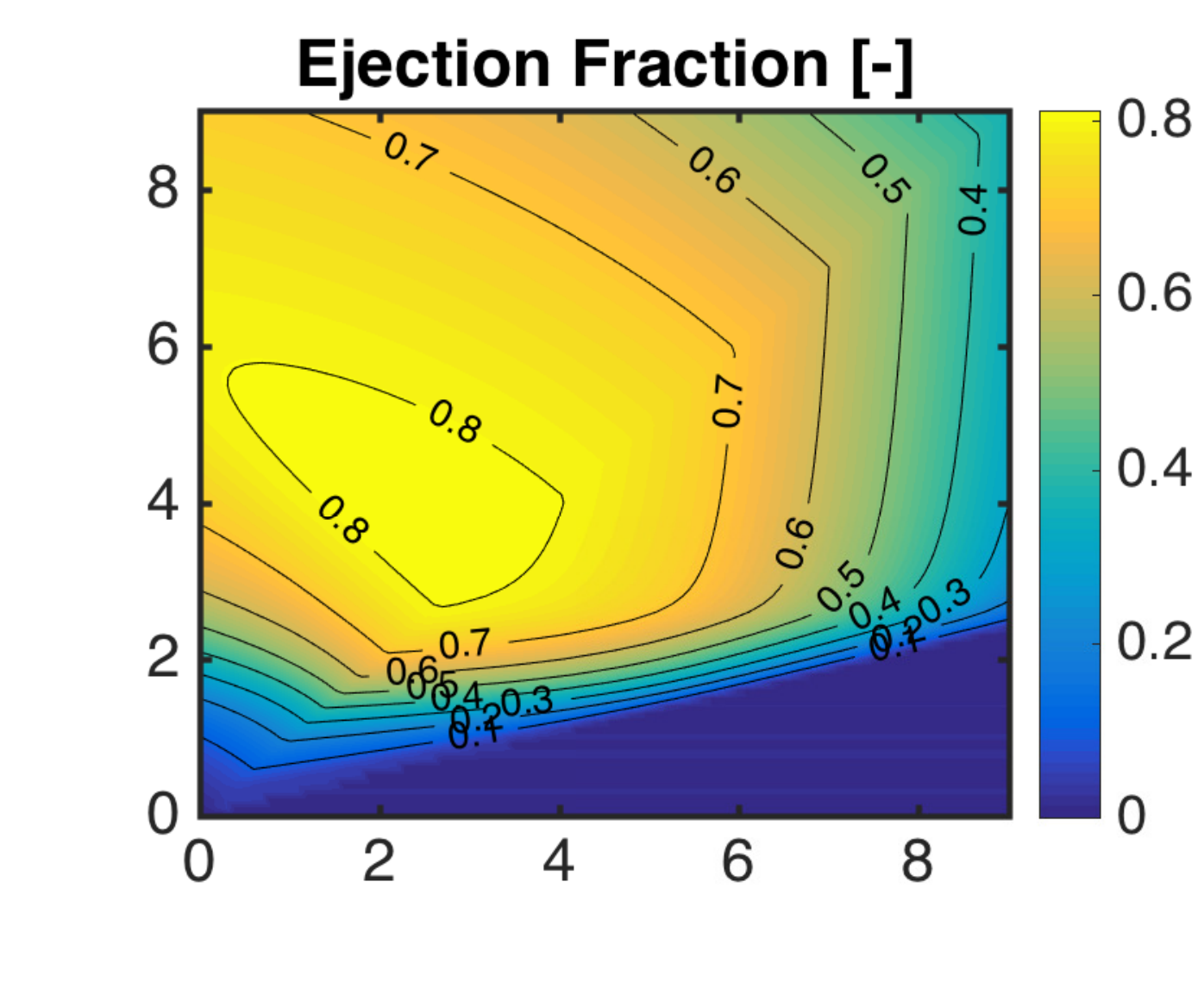}}
\subfloat{\label{fig:a}\includegraphics[width=0.33\textwidth]{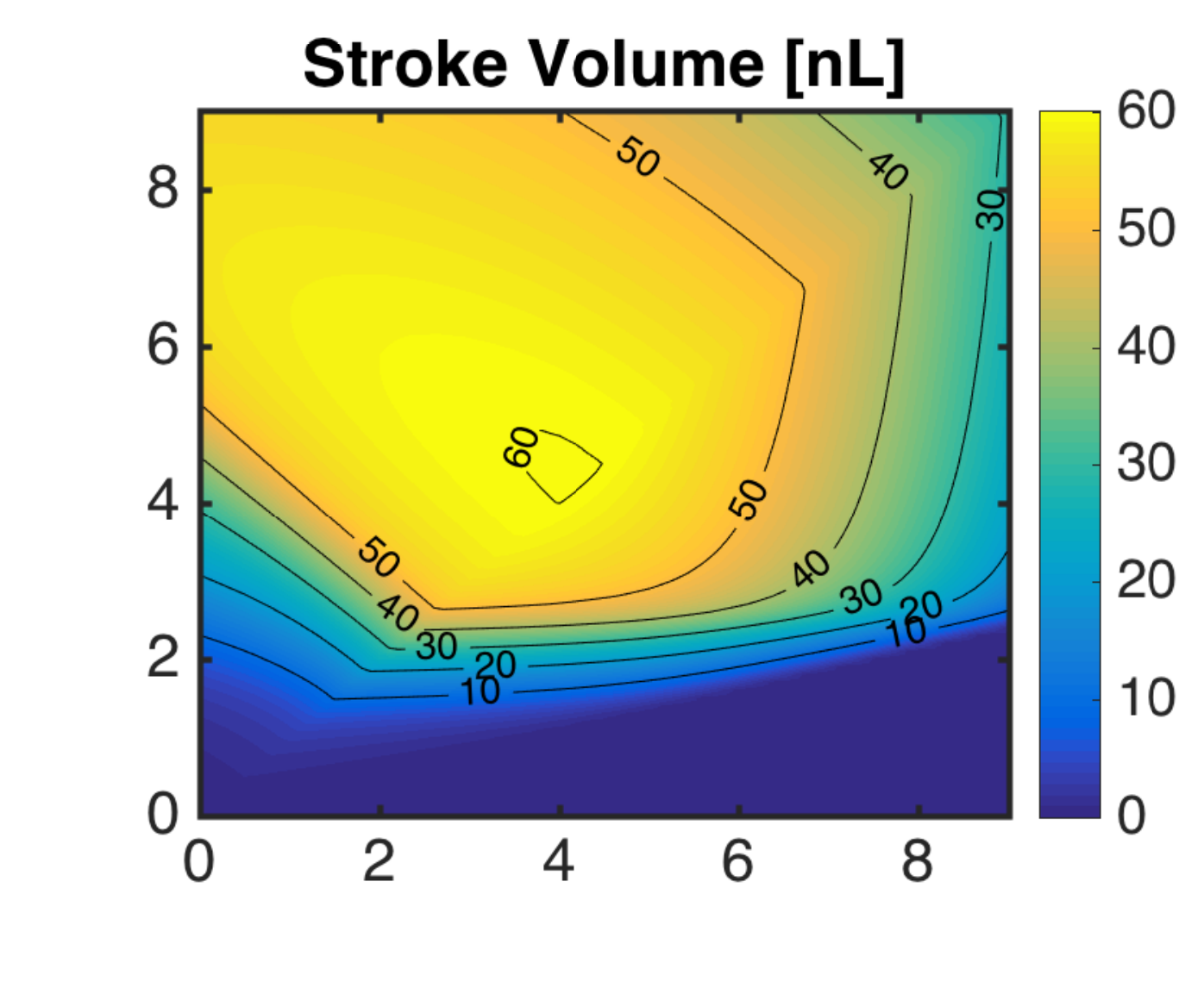}} \\
\subfloat{\label{fig:a}\includegraphics[width=0.33\textwidth]{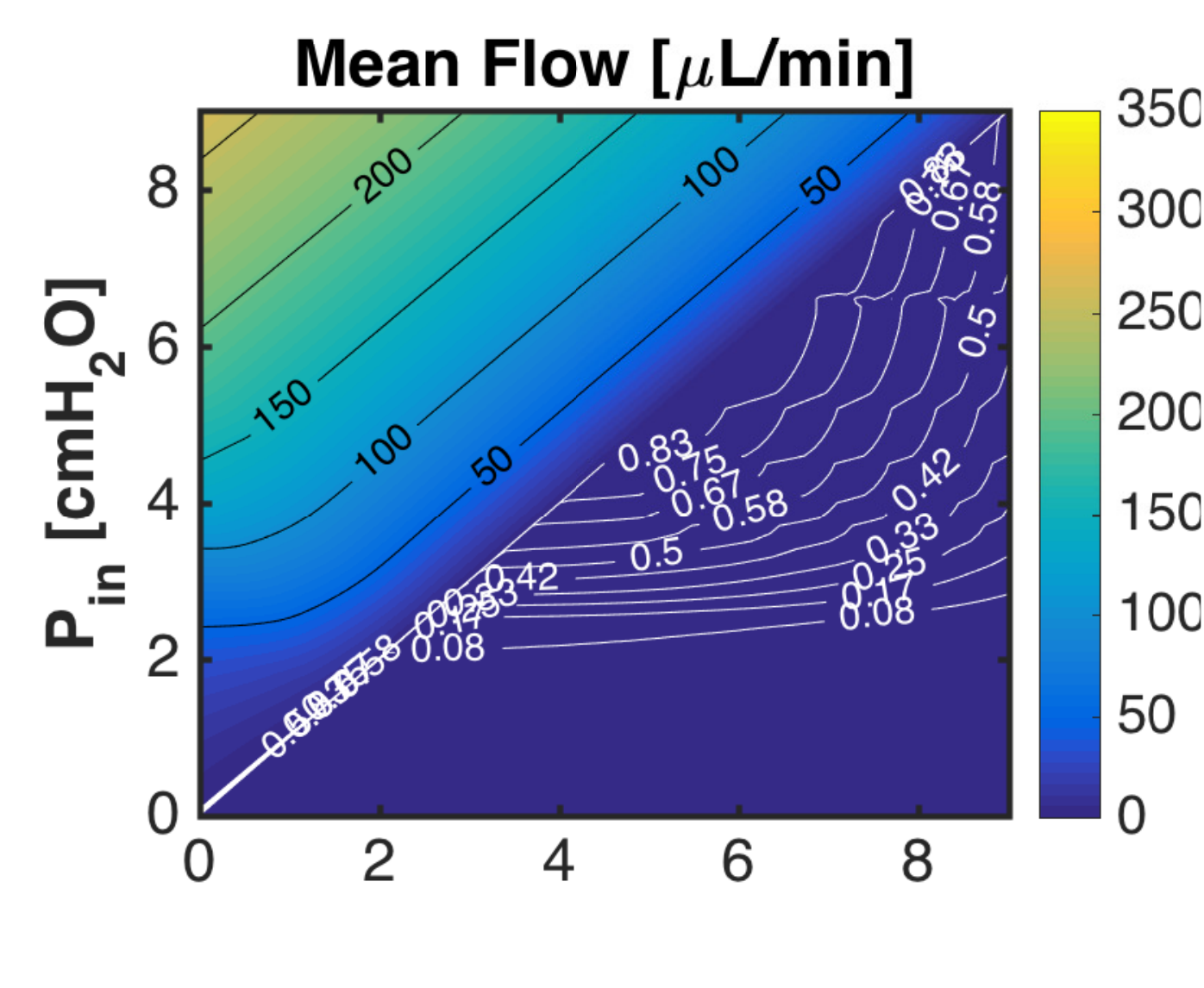}}
\subfloat{\label{fig:a}\includegraphics[width=0.33\textwidth]{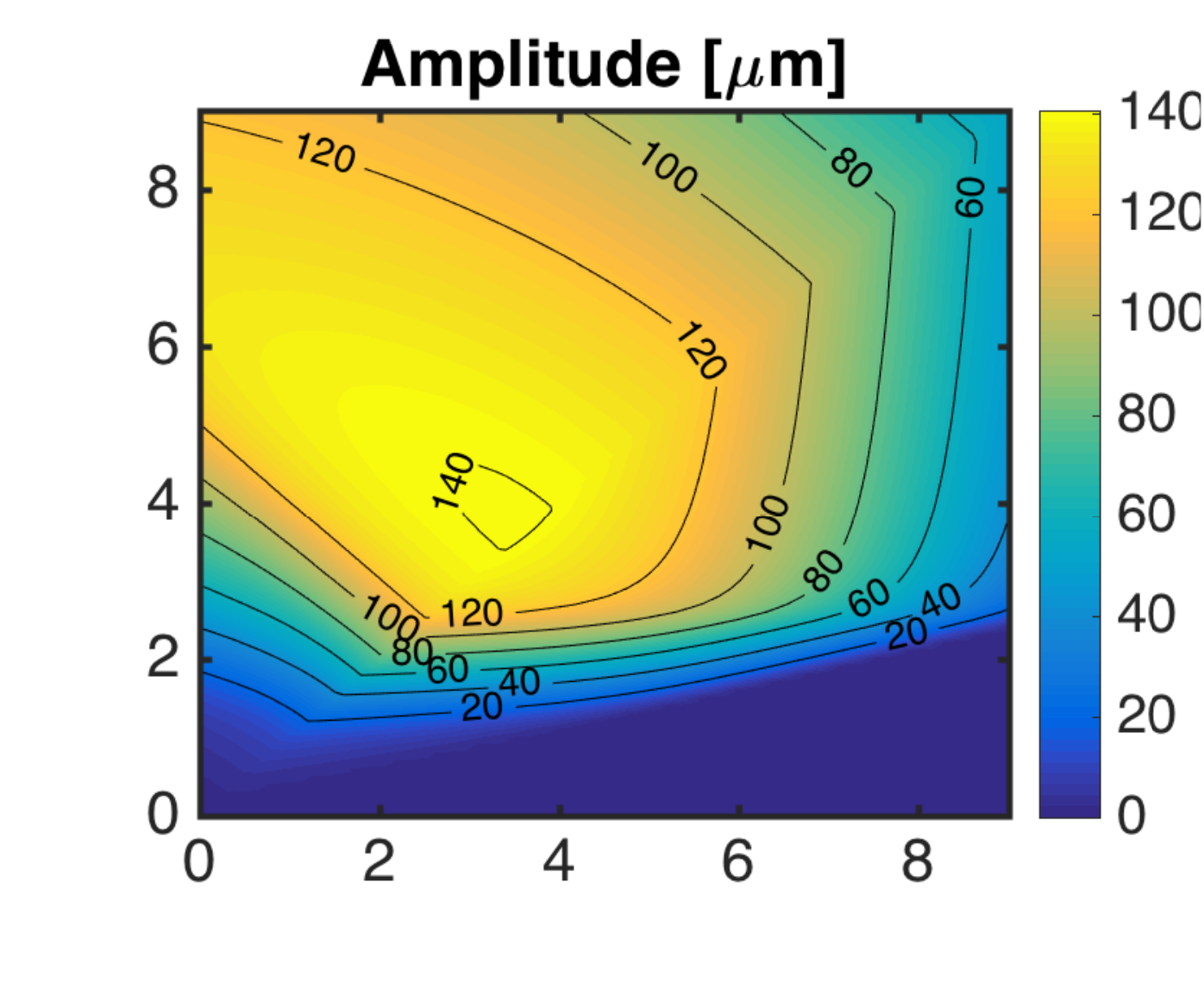}}
\subfloat{\label{fig:a}\includegraphics[width=0.33\textwidth]{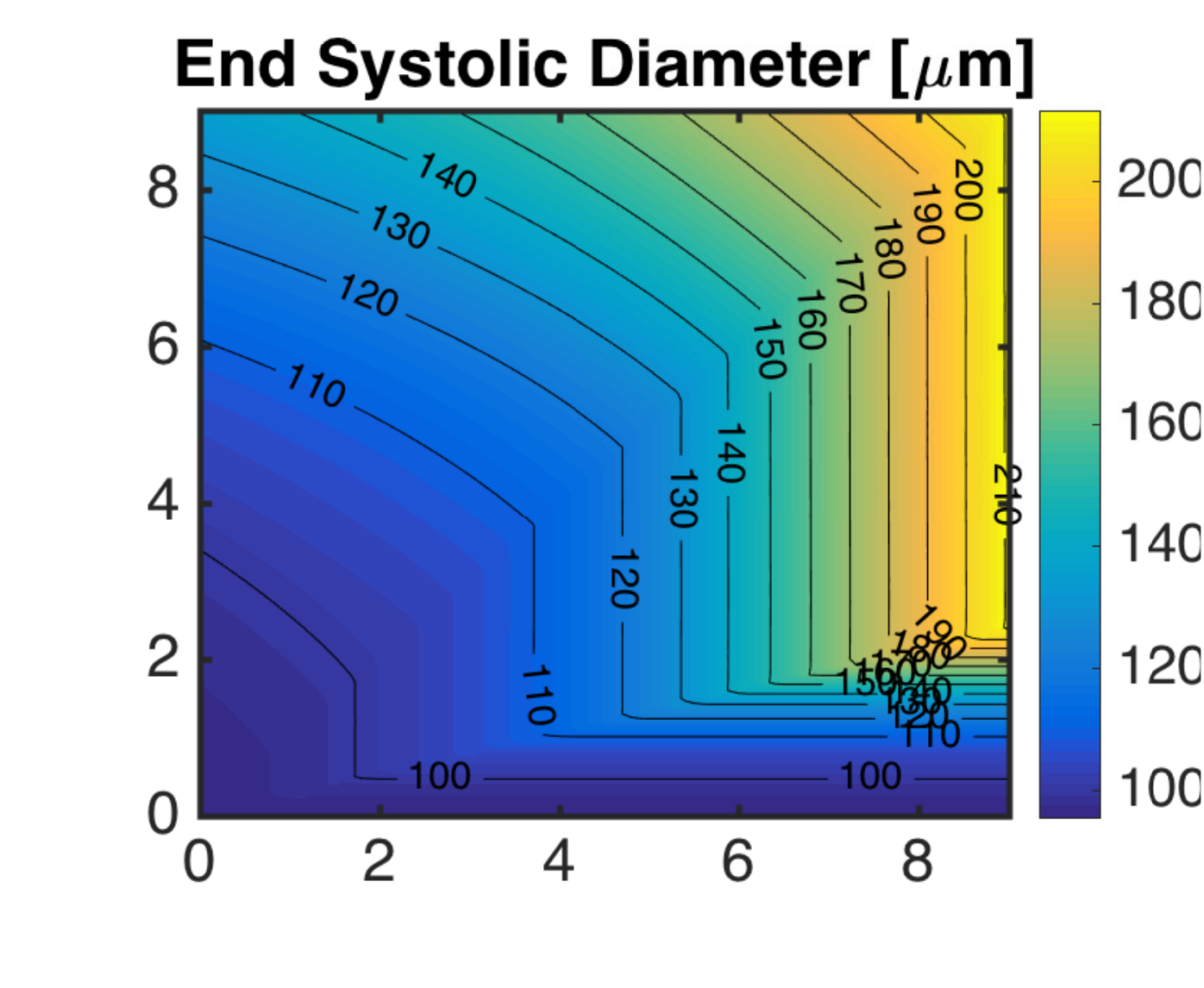}} \\
\subfloat{\label{fig:a}\includegraphics[width=0.33\textwidth]{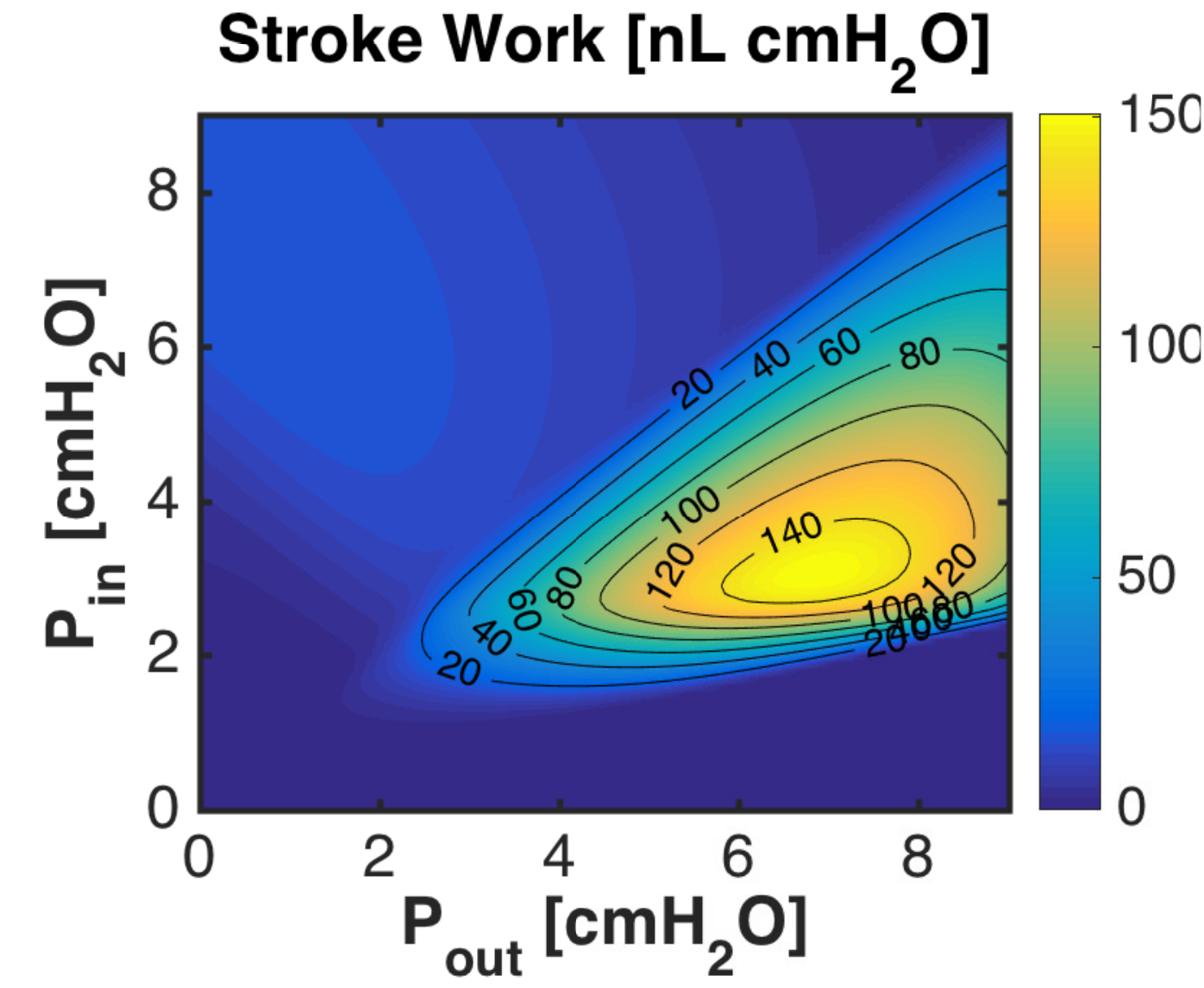}}
\subfloat{\label{fig:a}\includegraphics[width=0.33\textwidth]{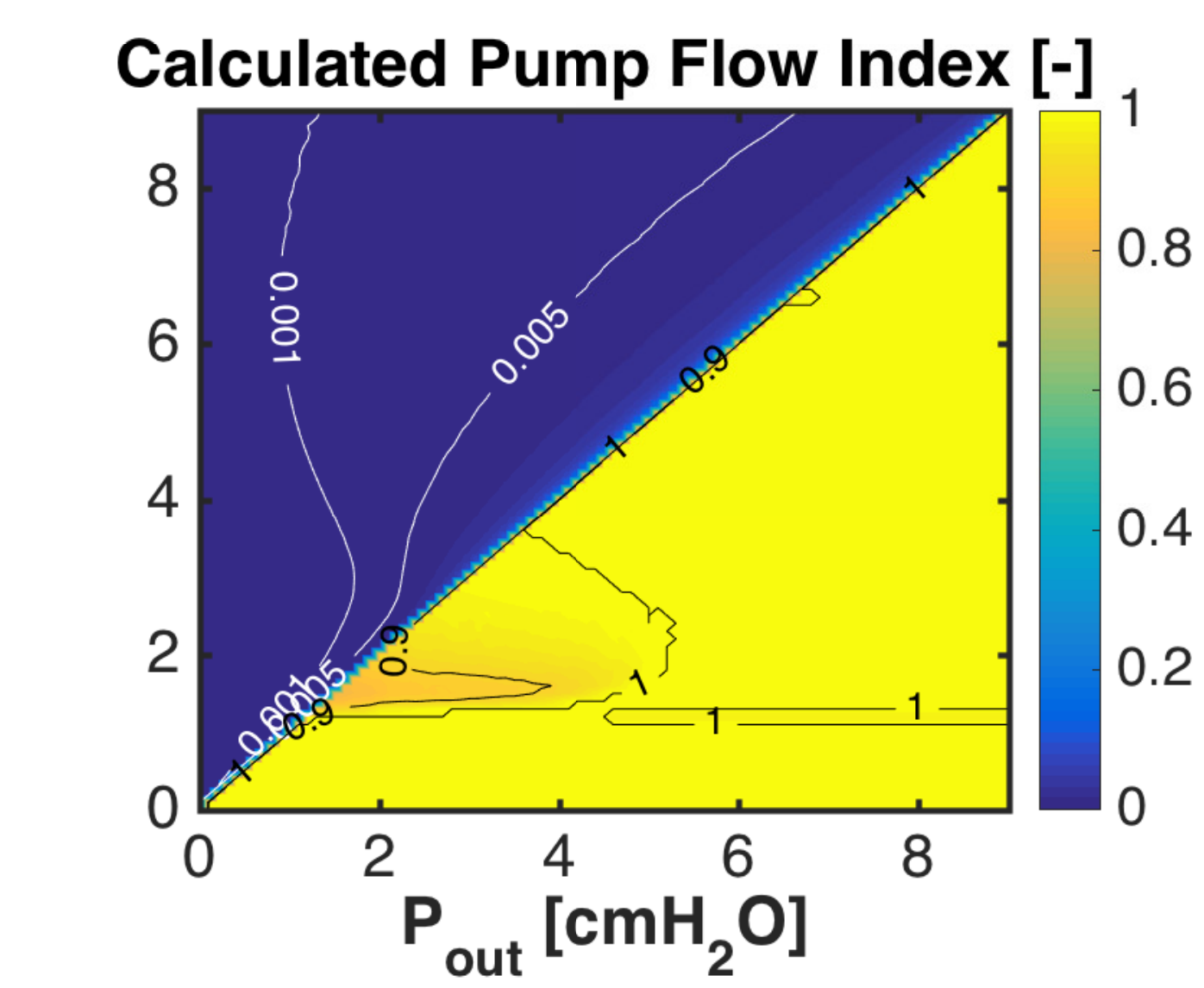}}
\subfloat{\label{fig:a}\includegraphics[width=0.33\textwidth]{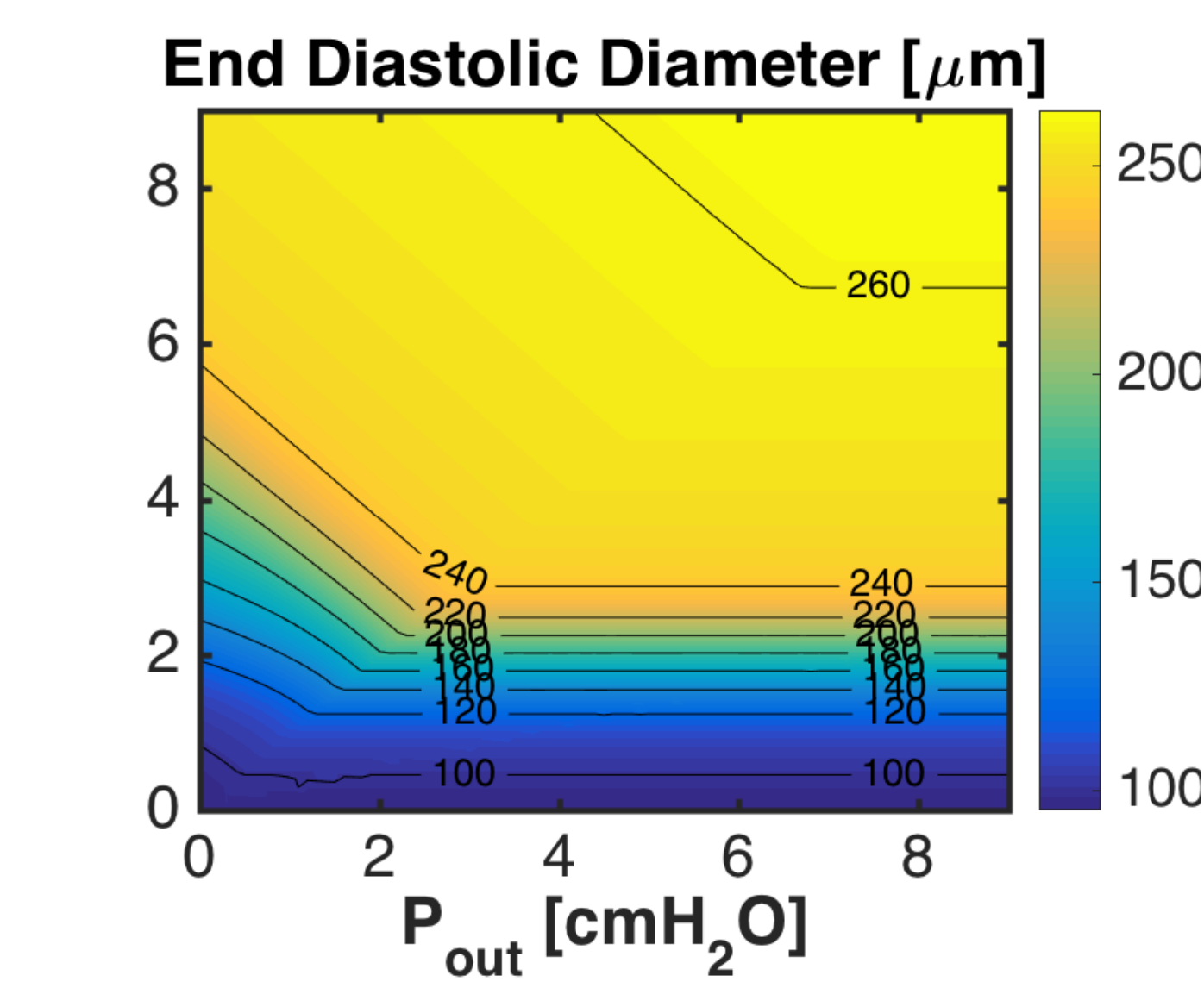}}
\caption{\scriptsize {\bf Analysis of lymphodynamical indexes by varying $P_{in}$ and $P_{out}$.} The frames show contour plot and isolines for: frequency, FPF, CPF, WSS, EF, SV, mean flow, mean WSS, AMP, ESD, SW, CPFI and EDD. We constructed a grid of points $\left(P_{out},P_{in}\right)$ with all the possible combinations of $P_{in}=\left(0,0.1,\dots,8.9,9\right)$ and $P_{out}=\left(0,0.1,\dots,8.9,9\right)$. For each combination of $P_{in}$ and $P_{out}$, we simulated a single lymphangion with two terminal valves, $t_{output}=160$ $s$ and $M=20$ computational cells to discretize the lymphangion.
} \label{fig:Indexes}
\end{figure*}
The aim of this study is to show lymphodynamical indexes for several combinations of $P_{in}$ and $P_{out}$ from 0 to 9 cmH$_2$O. We constructed a grid of points $\left(P_{out},P_{in}\right)$ with all the possible combinations of $P_{in}=\left(0,0.1,\dots,8.9,9\right)$ and $P_{out}=\left(0,0.1,\dots,8.9,9\right)$. For each combination of $P_{in}$ and $P_{out}$, we simulated a single lymphangion with two terminal valves, $t_{output}=160$ $s$ and $M=20$ computational cells. The total number of simulations was $91\times 91=8281$. We applied the numerical method at the boundaries as explained in Section \ref{sec:couplingvalves}. We analyzed the indexes reported in Table \ref{table:Indexes}.
\begin{table*}[t!] 
\centering 
\resizebox{\linewidth}{!}{ 
{\def\arraystretch{1.5}\tabcolsep=8pt
\begin{tabular}{l l c l c} 
\toprule[1.5pt] 
& {\bf Index} & {\bf Formula} & {\bf  Description} & {\bf Units} \\ \hline 
ESD & End-Systolic Diameter & - & Diameter at the end of lymphatic contraction & $\mu$m\\ 
EDD & End-Diastolic Diameter & - & Diameter at the beginning of filling & $\mu$m\\ 
ESP & End-Systolic Pressure & - & Pressure at the end of lymphatic contraction & cmH$_2$O\\ 
EDP & End-Diastolic Pressure & - & Pressure at the beginning of filling & cmH$_2$O\\ 
EF & Ejection Fraction & $1-\f{\text{EDV}}{\text{ESV}}$ & Fractional amount of ejected lymph & -\\ 
SV & Stroke Volume & $\text{EDV}-\text{ESV}$ & Ejected volume amount & nL \\ 
FPF & Fractional Pump Flow & $\text{EF}\times \text{FREQ}$ & Fractional change in lymphatic volume per minute & min$^{-1}$\\
CPF & Calculated Pump Flow & $\text{SV}\times \text{FREQ}$ & Flow produced by lymphatic contraction & $\mu$L h$^{-1}$\\
AMP & Amplitude & $\text{EDD}-\text{ESD}$&  Difference between diastolic and systolic diameter & $\mu$m\\
SW & Stroke Volume &$ \int P \mathrm{d}V$ & Area inside the pressure-volume loop & nL cmH$_2$O\\
CPFI & Calculated Pump Flow Index & $\left| \f{\text{CPF}}{{q}_{\text{mean}}}\right|$ & Ratio between propelled and passive flows & - \\
$\tau_{mean}$ & Time-averaged Wall Shear Stress & $\f{1}{t_2-t_1}\int_{t_1}^{t_2} \tau\left(\f{L}{2},t\right) \mathrm{d}t$ & Averaged WSS during a lymphatic cycle & dyne cm$^{-2}$ \\
$q_{mean}$ & Time-averaged Flow & $\f{1}{t_2-t_1}\int_{t_1}^{t_2} q\left(\f{L}{2},t\right) \mathrm{d}t$ & Averaged flow during a lymphatic cycle & $\mu$L h$^{-1}$ \\
\bottomrule[1.5pt] 
\end{tabular}}} \caption{\scriptsize {\bf Lymphodynamical indexes.} $t_1$ and $t_2$ corresponds to the initial and ending time of the lymphatic contractions and $L$ is the lymphangion length. } \label{table:Indexes} 
\end{table*} 
CPF is the flow produced only by the contraction and does not take into account the passive flow induced by a positive transaxial-pressure gradient. To understand the influence of the lymphatic pump on the total flow during a lymphatic cycle, we introduce the {\em Calculated Pump Flow Index}, which indicates if the flow is produced by contractions (CPFI $\approx$ 1) or by a positive transaxial-pressure gradient (CPFI $\approx$ 0). We observe that assuming constant frequency of contraction, the CPFI can also be calculated as 
\begin{equation}
\text{CPFI}=\left|\f{\text{SV}}{\displaystyle \int_{t_1}^{t_2} q\left(\f{L}{2},t\right) \mathrm{d}t}\right|\;,
\end{equation}
 which is the ratio between the ejected lymph volume due to the lymphatic contraction and the total lymph volume difference. Davis et al. \cite{Davis:2012a} performed an analysis of the lymphodynamical indexes by keeping fixed the upstream pressure $P_{in}$ and varying the downstream pressure $P_{out}$ from 0 to 18 cmH$_2$O. Similar experiments were performed by Scallan et al. \cite{Scallan:2012a} where the effect of increasing the upstream pressure $P_{in}$ while keeping fixed the downstream pressure $P_{out}$ was studied. In the present paper, both experimental setups are comprised in the numerical simulations.

Fig. \ref{fig:Indexes} shows the numerical results of the analysis through contour plots and isolines.
To describe the figure, we divide the $P_{in}$-$P_{out}$ space in two regions: 1) the negative transaxial-pressure gradient $\Delta P=P_{in}-P_{out}<0$ region (lower triangle) and 2) the positive transaxial-pressure gradient $\Delta P=P_{in}-P_{out}>0$ region (upper triangle). The two regions are therefore divided by the bisector $P_{in}=P_{out}$.

{\em \bf Negative pressure gradient $\Delta P=P_{in}-P_{out}<0$}. Here the upstream and downstream valves open non-linearly and close during the lymphatic cycle (results not shown). The frequency increases as $P_{in}$ rises, and this is in agreement with Fig. \ref{fig:EFMC}. This comes from muscle-stretch feedback from the EFMC model in Eq. \eqref{eq:Trigger}. The frequency does not increase when $P_{out}$ rises. This might be surprising because it is well-known that contraction-waves propagate between lymphangions through gap-junctional communications \cite{Zawieja:1993a, Shields:2009a}. Therefore, if the downstream pressure rises, then it would be natural to expect a rise in frequency among the neighbouring lymphangions. However, we modelled only a single lymphangion with variable downstream and upstream pressures, and we did not take into account the interaction between adjacent lymphangions. 
Moreover, the gap-junctional communications between lymphangions have not been modelled in this paper for the sake of simplicity, though we speculate that it is mathematically possible. 
EF tends to decrease as the $P_{out}$ increases, while increases when $P_{in}$ increases. FPF combines both frequency and EF: it increases when $P_{in}$ rises, and it decreases when $P_{out}$ increases. Intriguingly, the maximum of FPF is when $P_{in}\approx P_{out}$ and near 6 cmH$_2$O. For higher pressures, FPF decreases. The results for FPF are not comparable with Davis et al. \cite{Davis:2012a}, as the frequency remain constant when $P_{out}$ rises. The SV and AMP follow the same behaviour of EF. ESD increases only when $P_{out}$ rises, and it remains constant when $P_{in}$ increases. On the contrary, EDD remains constant when $P_{out}$ increases, and this is in agreement with \cite{Davis:2012a}. SW is maximum for $P_{out}\approx 7$ cmH$_2$O and $P_{in}\approx 3$ cmH$_2$O, and it tends to decrease elsewhere. The mean flow and CPF are comparable almost everywhere. This means that lymph flow is produced only by contractions, and not by a passive flow induced by pressure gradients. This is summarised in CPFI, which is almost 1 everywhere in this region. Finally, note that WSS is generally low in this region. 

{\em \bf Positive pressure gradient $\Delta P=P_{in}-P_{out}>0$}. Here the valves remain open for most of the time during the lymphatic cycle (results not shown).
The mean flow increases as $P_{in}$ rises, and it decreases as $P_{out}$ increases. Subsequently, the absolute value of WSS rises when $P_{in}$ increases, and this gives a negative chronotropic effect on the frequency. In other words, the frequency decreases as the absolute value of WSS increases. This comes from function $f_{NO}$ in Eq. \eqref{eq:Trigger}: the greater the absolute value of the WSS, the greater the negative chronotropic effect on the lymphatic contractions. The decrease in frequency is controlled by parameters $k_{NO}$, which in these simulations is set to $k_{NO}=0.8$. Lower values of $k_{NO}$ would have weakened the negative chronotropic effects, leading to higher rates of lymphatic contractions. Observe that the CPF differ from the mean flow insofar as the CPF only takes into account the flow given by contractions. This is summarised in CPFI which is almost zero in this region. The ESD and EDD increase when $P_{in}$ and $P_{out}$ rise. EF, SV and AMP share a similar behaviour and reach their maximum value at $P_{in}=P_{out}\approx 4$ cmH$_2$O.

\subsection{Sensitivity analyses}


 \begin{sidewaystable}
 
\resizebox{\linewidth}{!}{ 
\centering 
\begin{tabular}{c | c l| c c c c c c c c c c c c c c } 
\toprule[1.5pt] 
$P_{in}$ = 3 - $P_{out}$ = 4& \multicolumn{2}{ c |}{$\mathbf{P}(\mathbf{X})$} &7.80 $\pm$ 5.43&0.74 $\pm$ 0.05&47.19 $\pm$ 14.78&5.75 $\pm$ 4.03&22.87 $\pm$ 16.60&-0.12 $\pm$ 0.09&118.93 $\pm$ 18.05&232.78 $\pm$ 34.08&4.49 $\pm$ 0.15&3.00 $\pm$ 0.00&3.07 $\pm$ 0.05&21.73 $\pm$ 15.79\\ 
\bottomrule[1.5pt] 
\multirow{2}{*}{$\mathbf{X}$} & \multicolumn{2}{ c |}{\multirow{2}{*}{$\bar{\mathbf{S}}$ [\%]}} &Frequency&EF&SV&FPF&CPF&WSS&ESD&EDD&ESP&EDP&Mean Pressure & Mean Flow  \\  
& & & [min$^{-1}$]& [-]& [nL]& [min$^{-1}$]& [$\mu$L h$^{-1}$]& [dyne cm$^{-2}$]& [$\mu$m]& [$\mu$m]& [cmH$_2$O]& [cmH$_2$O]& [cmH$_2$O]& [$\mu$L h$^{-1}$]\\ \hline 
47.52 $\pm$ 8.21 &$r_0$ & [$\mu$m] &-&-&{\color{red}207.7} {\color{red} $\pm$ 61.9}&0.8 $\pm$ 1.2&{\color{red}249.1} {\color{red} $\pm$ 170.3}&{\color{red}121.6} {\color{red} $\pm$ 82.9}&{\color{red}119.2} {\color{red} $\pm$ 18.2}&{\color{red}105.6} {\color{red} $\pm$ 15.1}&-&2.0 $\pm$ 1.3&-&{\color{red}235.7} {\color{red} $\pm$ 161.6}\\  
100.23 $\pm$ 16.92 &$a_1$ & [s$^{-1}$] &5.7 $\pm$ 6.6&-0.9 $\pm$ 1.3&-1.1 $\pm$ 1.6&3.2 $\pm$ 5.5&3.4 $\pm$ 6.2&-1.4 $\pm$ 6.5&1.5 $\pm$ 2.2&-&-&-&-&3.0 $\pm$ 7.1\\  
0.49 $\pm$ 0.09 &$a_2$ & [-] &{\color{blue}-76.9} {\color{blue} $\pm$ 108.6}&22.7 $\pm$ 29.5&{\color{darkgreen}26.3} {\color{darkgreen} $\pm$ 34.0}&{\color{darkgreen}-30.6} {\color{darkgreen} $\pm$ 73.6}&{\color{darkgreen}-33.4} {\color{darkgreen} $\pm$ 80.8}&-4.3 $\pm$ 83.0&{\color{darkgreen}-40.1} {\color{darkgreen} $\pm$ 50.7}&-&-&1.5 $\pm$ 0.7&3.6 $\pm$ 2.3&{\color{darkgreen}-33.4} {\color{darkgreen} $\pm$ 81.0}\\  
25.09 $\pm$ 4.31 &$a_3$ & [-] &{\color{darkgreen}-30.5} {\color{darkgreen} $\pm$ 46.4}&6.5 $\pm$ 8.2&7.5 $\pm$ 9.4&-16.1 $\pm$ 34.7&-17.0 $\pm$ 36.2&4.6 $\pm$ 32.7&-11.5 $\pm$ 14.2&-&-&0.6 $\pm$ 0.3&1.5 $\pm$ 0.9&-17.1 $\pm$ 36.4\\  
3.00 $\pm$ 0.52 &$b_1$ & [s$^{-1}$] &{\color{darkgreen}28.5} {\color{darkgreen} $\pm$ 49.1}&-2.5 $\pm$ 3.0&-3.5 $\pm$ 4.3&19.8 $\pm$ 39.1&20.4 $\pm$ 40.1&-14.4 $\pm$ 33.1&4.0 $\pm$ 4.9&-&-&-&-0.6 $\pm$ 0.5&20.6 $\pm$ 40.4\\  
109.79 $\pm$ 19.28 &$c_1$ & [s$^{-1}$] &-0.9 $\pm$ 1.7&4.6 $\pm$ 6.4&5.4 $\pm$ 7.4&4.0 $\pm$ 5.8&4.9 $\pm$ 7.4&-10.8 $\pm$ 13.7&-8.2 $\pm$ 11.0&-&-&-&-&5.5 $\pm$ 8.4\\  
3.01 $\pm$ 0.53 &$c_2$ & [s$^{-1}$] &15.3 $\pm$ 20.5&-&-&14.0 $\pm$ 18.7&17.5 $\pm$ 23.8&-13.5 $\pm$ 19.0&-&-&-&-1.3 $\pm$ 0.4&-&17.3 $\pm$ 23.3\\  
0.10 $\pm$ 0.02 &$R_I$ & [s$^{-1}$] &4.6 $\pm$ 4.4&-&-&4.6 $\pm$ 4.5&5.6 $\pm$ 5.8&-4.8 $\pm$ 5.4&-&-&-&-&-&5.3 $\pm$ 6.6\\  
10.10 $\pm$ 1.70 &$k_{rel}$ & [s$^{-1}$] &2.0 $\pm$ 2.9&-1.4 $\pm$ 2.2&-1.6 $\pm$ 2.5&-&-&1.8 $\pm$ 4.0&2.5 $\pm$ 3.8&-&-&-&-&-0.5 $\pm$ 3.8\\  
9.92 $\pm$ 1.72 &$n_{Ca}$ & [-] &{\color{blue}-52.0} {\color{blue} $\pm$ 34.5}&-&-&{\color{darkgreen}-47.1} {\color{darkgreen} $\pm$ 32.4}&{\color{blue}-51.4} {\color{blue} $\pm$ 40.1}&{\color{darkgreen}45.6} {\color{darkgreen} $\pm$ 33.2}&-&-&-&-&-1.1 $\pm$ 0.8&{\color{blue}-52.3} {\color{blue} $\pm$ 40.0}\\  
7.78 $\pm$ 1.31 &$A_{Ca}/A_0$ & [-] &{\color{red}-377.7} {\color{red} $\pm$ 281.3}&-&-&{\color{red}-349.8} {\color{red} $\pm$ 257.3}&{\color{red}-405.2} {\color{red} $\pm$ 289.7}&{\color{red}339.3} {\color{red} $\pm$ 255.3}&-&-&-&-&-8.2 $\pm$ 5.8&{\color{red}-406.0} {\color{red} $\pm$ 287.1}\\  
0.50 $\pm$ 0.09 &$k_{No}$ & [-] &-&-&-&-&-&-&-&-&-&-&-&-\\  
5.96 $\pm$ 1.03 &$t_{No}$ & [dyne cm$^{-2}$] &-&-&-&-&-&-&-&-&-&-&-&-\\  
1.20 $\pm$ 0.20 &$n_{No}$ & [-] &-&-&-&-&-&-&-&-&-&-&-&-\\  
9.89 $\pm$ 1.76 &$K_{vo}$ & [Pa$^{-1}$ s$^{-1}$] &-&-&-&-&-&-&-&-&-&-&-&-\\  
10.01 $\pm$ 1.74 &$K_{vc}$ & [Pa$^{-1}$ s$^{-1}$] &-&-&-&-&-&-7.2 $\pm$ 4.5&-&-&-&-&-&5.5 $\pm$ 4.2\\  
99.49 $\pm$ 17.43 &$L_{eff}$ & [$\mu$m] &-&-&-&-&-&-6.4 $\pm$ 4.4&-&-&-&-0.6 $\pm$ 0.4&-&5.0 $\pm$ 4.9\\  
1007.99 $\pm$ 176.63 &$\rho$ & [kg m$^{-3}$] &-&-&-&-&-&-&-&-&-&-&-&-\\  
1.00 $\pm$ 0.17 &$\mu$ & [cP] &-&-&-&-&-&{\color{red}-114.7} {\color{red} $\pm$ 79.6}&-&-&-&-1.0 $\pm$ 0.6&-&4.9 $\pm$ 4.8\\  
2.01 $\pm$ 0.35 &$\gamma$ & [-] &-&-&-&-&-&{\color{blue}-56.7} {\color{blue} $\pm$ 39.8}&-&-&-&-0.5 $\pm$ 0.4&-&2.6 $\pm$ 3.5\\  
134280.92 $\pm$ 22699.60 &$E_{max}$ & [Pa] &-5.7 $\pm$ 7.9&{\color{darkgreen}29.1} {\color{darkgreen} $\pm$ 15.7}&{\color{darkgreen}33.1} {\color{darkgreen} $\pm$ 19.2}&24.6 $\pm$ 20.4&{\color{darkgreen}28.5} {\color{darkgreen} $\pm$ 24.9}&{\color{blue}-78.6} {\color{blue} $\pm$ 55.8}&{\color{blue}-54.3} {\color{blue} $\pm$ 25.8}&-&-&-&-&{\color{darkgreen}29.7} {\color{darkgreen} $\pm$ 26.5}\\  
34591.78 $\pm$ 5844.81 &$E_{min}$ & [Pa] &{\color{red}-195.7} {\color{red} $\pm$ 171.4}&-19.6 $\pm$ 17.3&{\color{blue}-80.6} {\color{blue} $\pm$ 52.3}&{\color{red}-193.2} {\color{red} $\pm$ 159.6}&{\color{red}-270.3} {\color{red} $\pm$ 191.3}&{\color{red}199.8} {\color{red} $\pm$ 166.2}&-1.9 $\pm$ 3.5&{\color{darkgreen}-32.9} {\color{darkgreen} $\pm$ 22.6}&-&1.5 $\pm$ 0.8&-3.8 $\pm$ 3.5&{\color{red}-275.6} {\color{red} $\pm$ 192.3}\\  
\bottomrule[1.5pt] 
\end{tabular}} 
\caption{\scriptsize {\bf Sensitivity analysis of the one-dimensional lymph flow equations coupled to the EFMC model and valve dynamics. Negative pressure gradient}. Symbol $"-"$ indicates that the absolute value of the sensitivity index $\bar{S}_{i,j}$ was less than 0.5 and therefore parameter $x_i$ did not influence index $P_j$. Likewise, green-coloured parameters show an influence between 25 and 50, blue-coloured parameters show an influence between 50 and 100, and red-coloured parameters show an influence greater than 100. Results are shown using as mean $\pm$ SD. } \label{table:SA<0}

\end{sidewaystable}

The mathematical model for lymphatic collectors proposed here depends on several parameters shown in Table \ref{table:parameters}, and the lymphodynamical indexes studied in Section \ref{sec:pIn_pOut} are strongly related to these parameters. 
To investigate the influence of each parameter on the indexes, we performed a sensitivity analysis following the idea of \cite{Griensven:2006a}. This procedure was already applied, for instance, to find the most relevant model parameters of a zero-dimensional model of the cardiovascular system \cite{Liang:2014a}. 

The method is divided into a local and global analysis.
In the local analysis we calculated $N$ local sensitivity matrixes $S_{i,j}^k$, for $k=1,\dots,N$, as follows: starting from the reference value in Table \ref{table:parameters}, we randomly varied each parameter from $70$ to $130\% $ and obtained a new set of parameters. Here, the baseline value for $k_{NO}$ was set to 0.5. With this varied set of parameters, we calculated the local sensitivity matrix as follows

\begin{equation}\label{eq:sensitivity}
S_{i,j}^k=\left|\f{x_i}{P_j(\mathbf{X})}\right|\f{\partial P_j\left(\mathbf{X}\right)}{\partial x_i}\times 100\;, 
\end{equation}
where $\mathbf{X}=\left(x_1,x_2,\dots,x_m\right)$ is the vector of the varied model parameters, $\mathbf{P}=\left(P_1,P_2,\dots,P_n\right)$ is the vector of the lymphodynamical indexes and $\mathbf{S}^k=\left(S_{i,j}^k\right)_{i,j}$ is local sensitivity matrix. The value $S_{i,j}^k$ represents the non-dimensional percentage change in index $P_j$ induced by a small change in parameter $x_i$. For instance, if the model parameter $x_i$ does not influence index $P_j$, then $S_{i,j}^k$ will be almost zero. Viceversa, if there is a significant influence of $x_i$ on $P_j$, then the absolute value of $S_{i,j}^k$ will be greater than zero. A positive sign of $S_{i,j}^k$ indicates that an increase of parameter $x_i$ induces an increase of index $P_j$. Viceversa, a negative sign of $S_{i,j}^k$ indicates that an increase of parameter $x_i$ induces a decrease of index $P_j$.

Subsequently, in the global analysis we performed a statistical analysis of $S_{i,j}^k$ by calculating its mean $\bar{S}_{i,j}$ and its standard deviation $\sigma_{i,j}$, namely
\begin{equation}\label{eq:SensitivityMean}
\bar{S}_{i,j}=\f{1}{N}\sum_{k=1}^{N}S^k_{i,j}\;,
\end{equation}

\begin{equation}\label{eq:SensitivitySD}
\sigma_{i,j}=\sqrt{\f{1}{N}\sum_{k=1}^{N}\left(S^k_{i,j}-\bar{S}_{i,j}\right)^2}\;.
\end{equation}
A large standard deviation $\sigma_{i,j}$ indicates a strong correlation of the studied parameter with the remaining parameters in determining the sensitivity index. 
To calculate \eqref{eq:SensitivityMean} and \eqref{eq:SensitivitySD}, we removed possible outliers by discarding the data below the $3$rd percentile and above the $97$th percentile.
 
The partial derivative in Eq. \eqref{eq:sensitivity} were approximated using a second-order finite difference method based on a percentage change of the parameter as follows:
\begin{equation}
S_{i,j}^k\approx\f{\text{sgn}(x_i)}{\left|P_j(\mathbf{X})\right|}\f{P_j\left(\mathbf{X}^{i,\epsilon_+}\right)-P_j\left(\mathbf{X}^{i,\epsilon_-}\right)}{2\epsilon}\times 100\;,
\end{equation}
where
\begin{equation}
\mathbf{X}^{i,\epsilon_{\pm}}=\left(x_1,\dots,x_i\left(1\pm \epsilon\right),\dots,x_m\right)\;.
\end{equation}
The parameter $\epsilon$ was chosen as $\epsilon=0.05$. 
Compared to \cite{Griensven:2006a} and \cite{Liang:2014a}, we did not constructed a stratified sampling space, but rather a simple random variation in the considered range. Here, we studied parameter ${A_{Ca}}/{A_0}$ with baseline value $7.75$.

Based on the results of the lymphodynamical indexes shown in Section \ref{sec:pIn_pOut}, we performed two sensitivity analyses: one for a positive pressure gradient $\Delta P=P_{in}-P_{out}>0$ and one for a negative pressure gradient $\Delta P=P_{in}-P_{out}<0$. For each analysis, we calculated $N=500$ local sensitivity matrixes. We used $M=20$ computational cells to discretize the one-dimensional lymph vessel and output time $t_{output}=60$ $s$. Results are shown in Tables \ref{table:SA<0} and \ref{table:SA>0}. In both tables, results are shown in the following manner: the second column from the left shows parameters $\mathbf{X}$ of the model while the first column shows their means $\pm$ SDs. The second row from the top shows the studied lymphodynamical indexes $\mathbf{P}$ while the first row shows the resulting means $\pm$ SDs. Symbol $"-"$ indicates that the absolute value of the sensitivity index $\bar{S}_{i,j}$ was less than 0.5 and therefore parameter $x_i$ did not influence index $P_j$. Likewise, green-coloured parameters show an influence between 25 and 50, blue-coloured parameters show an influence between 50 and 100, and red-coloured parameters show an influence greater than 100. 

 \begin{sidewaystable}
 
\resizebox{\linewidth}{!}{ 
\centering 
\begin{tabular}{c | c l| c c c c c c c c c c c c c c } 
\toprule[1.5pt] 

$P_{in}$ = 4 - $P_{out}$ = 2& \multicolumn{2}{ c |}{$\mathbf{P}(\mathbf{X})$} &5.75 $\pm$ 4.62&0.78 $\pm$ 0.03&50.36 $\pm$ 13.91&4.54 $\pm$ 3.72&17.44 $\pm$ 14.75&-12.98 $\pm$ 2.05&108.40 $\pm$ 15.53&231.99 $\pm$ 31.88&3.37 $\pm$ 0.08&3.01 $\pm$ 0.01&3.02 $\pm$ 0.01&4821.88 $\pm$ 2221.97\\ 
\bottomrule[1.5pt] 
\multirow{2}{*}{$\mathbf{X}$} & \multicolumn{2}{ c |}{\multirow{2}{*}{$\bar{\mathbf{S}}$ [\%]}} &Frequency&EF&SV&FPF&CPF&WSS&ESD&EDD&ESP&EDP&Mean Pressure & Mean Flow  \\  
& & & [min$^{-1}$]& [-]& [nL]& [min$^{-1}$]& [$\mu$L h$^{-1}$]& [dyne cm$^{-2}$]& [$\mu$m]& [$\mu$m]& [cmH$_2$O]& [cmH$_2$O]& [cmH$_2$O]& [$\mu$L h$^{-1}$]\\ \hline 
47.59 $\pm$ 8.44 &$r_0$ & [$\mu$m] &-18.5 $\pm$ 18.4&-&{\color{red}301.1} {\color{red} $\pm$ 82.2}&-19.1 $\pm$ 19.3&{\color{red}965.1} {\color{red} $\pm$ 786.6}&{\color{darkgreen}-47.5} {\color{darkgreen} $\pm$ 18.8}&{\color{red}108.6} {\color{red} $\pm$ 15.5}&{\color{red}118.5} {\color{red} $\pm$ 16.3}&-&-&-&{\color{red}328.7} {\color{red} $\pm$ 134.1}\\  
99.68 $\pm$ 17.26 &$a_1$ & [s$^{-1}$] &14.1 $\pm$ 19.9&-&-0.6 $\pm$ 0.8&13.1 $\pm$ 20.2&17.1 $\pm$ 28.3&-&0.6 $\pm$ 0.8&-&-&-&-&-\\  
0.50 $\pm$ 0.09 &$a_2$ & [-] &{\color{red}-182.1} {\color{red} $\pm$ 287.9}&8.5 $\pm$ 10.4&12.6 $\pm$ 15.6&{\color{red}-152.6} {\color{red} $\pm$ 273.5}&{\color{red}-200.6} {\color{red} $\pm$ 342.7}&{\color{darkgreen}-27.7} {\color{darkgreen} $\pm$ 17.0}&-15.8 $\pm$ 19.5&-&-&-&-&-13.7 $\pm$ 8.8\\  
24.87 $\pm$ 4.24 &$a_3$ & [-] &{\color{blue}-67.3} {\color{blue} $\pm$ 115.0}&2.4 $\pm$ 2.9&3.6 $\pm$ 4.4&{\color{blue}-60.6} {\color{blue} $\pm$ 116.9}&{\color{blue}-78.7} {\color{blue} $\pm$ 147.9}&-9.8 $\pm$ 6.0&-4.5 $\pm$ 5.4&-&-&-&-&-4.9 $\pm$ 3.2\\  
2.99 $\pm$ 0.52 &$b_1$ & [s$^{-1}$] &{\color{blue}55.8} {\color{blue} $\pm$ 109.9}&-0.9 $\pm$ 1.1&-1.7 $\pm$ 2.1&{\color{blue}54.3} {\color{blue} $\pm$ 111.5}&{\color{blue}69.0} {\color{blue} $\pm$ 132.9}&2.8 $\pm$ 3.5&1.5 $\pm$ 1.9&-&-&-&-&1.4 $\pm$ 1.6\\  
110.24 $\pm$ 18.67 &$c_1$ & [s$^{-1}$] &-2.4 $\pm$ 5.1&1.7 $\pm$ 2.2&2.5 $\pm$ 3.3&2.8 $\pm$ 5.5&4.4 $\pm$ 8.3&-1.5 $\pm$ 1.5&-3.2 $\pm$ 4.1&-&-&-&-&-0.8 $\pm$ 0.8\\  
2.99 $\pm$ 0.52 &$c_2$ & [s$^{-1}$] &{\color{darkgreen}36.6} {\color{darkgreen} $\pm$ 59.2}&-&-&{\color{darkgreen}39.4} {\color{darkgreen} $\pm$ 64.4}&{\color{blue}57.8} {\color{blue} $\pm$ 90.5}&7.6 $\pm$ 4.5&-&-&-&-&-&3.0 $\pm$ 1.9\\  
0.10 $\pm$ 0.02 &$R_I$ & [s$^{-1}$] &10.6 $\pm$ 13.1&-&-&12.3 $\pm$ 15.0&17.4 $\pm$ 21.0&-0.7 $\pm$ 1.0&-&-&-&-&-&-\\  
9.92 $\pm$ 1.66 &$k_{rel}$ & [s$^{-1}$] &3.8 $\pm$ 6.4&-&-0.8 $\pm$ 1.1&1.9 $\pm$ 5.3&2.6 $\pm$ 7.4&0.7 $\pm$ 0.9&0.9 $\pm$ 1.4&-&-&-&-&-\\  
10.08 $\pm$ 1.69 &$n_{Ca}$ & [-] &{\color{red}-132.8} {\color{red} $\pm$ 89.8}&-&-&{\color{red}-139.2} {\color{red} $\pm$ 95.4}&{\color{red}-194.1} {\color{red} $\pm$ 148.8}&7.2 $\pm$ 5.4&-&-&-&-&-&3.0 $\pm$ 2.5\\  
7.75 $\pm$ 1.31 &$A_{Ca}/A_0$ & [-] &{\color{red}-1272.7} {\color{red} $\pm$ 1070.8}&-&-&{\color{red}-1350.4} {\color{red} $\pm$ 1136.1}&{\color{red}-1938.7} {\color{red} $\pm$ 1626.2}&{\color{blue}68.0} {\color{blue} $\pm$ 54.9}&-&-&-&-&-2.1 $\pm$ 1.7&{\color{darkgreen}26.5} {\color{darkgreen} $\pm$ 22.1}\\  
0.51 $\pm$ 0.09 &$k_{No}$ & [-] &{\color{blue}-94.3} {\color{blue} $\pm$ 53.1}&-&-&{\color{blue}-99.0} {\color{blue} $\pm$ 55.8}&{\color{red}-145.3} {\color{red} $\pm$ 92.1}&5.1 $\pm$ 3.1&-&-&-&-&-&2.3 $\pm$ 1.6\\  
6.00 $\pm$ 1.04 &$t_{No}$ & [dyne cm$^{-2}$] &{\color{blue}56.5} {\color{blue} $\pm$ 34.2}&-&-&{\color{blue}59.6} {\color{blue} $\pm$ 36.3}&{\color{blue}83.8} {\color{blue} $\pm$ 52.1}&-3.1 $\pm$ 2.0&-&-&-&-&-&-1.3 $\pm$ 1.0\\  
1.22 $\pm$ 0.20 &$n_{No}$ & [-] &{\color{darkgreen}-25.6} {\color{darkgreen} $\pm$ 16.8}&-&-&{\color{darkgreen}-26.9} {\color{darkgreen} $\pm$ 17.6}&{\color{darkgreen}-39.7} {\color{darkgreen} $\pm$ 29.8}&1.5 $\pm$ 1.3&-&-&-&-&-&0.6 $\pm$ 0.6\\  
9.98 $\pm$ 1.68 &$K_{vo}$ & [Pa$^{-1}$ s$^{-1}$] &-&-&-&-&-&-&-&-&-&-&-&-\\  
10.06 $\pm$ 1.72 &$K_{vc}$ & [Pa$^{-1}$ s$^{-1}$] &&&&&&&&&&&&\\  
99.81 $\pm$ 17.06 &$L_{eff}$ & [$\mu$m] &{\color{darkgreen}31.0} {\color{darkgreen} $\pm$ 20.5}&-&-&{\color{darkgreen}33.4} {\color{darkgreen} $\pm$ 22.5}&{\color{darkgreen}46.5} {\color{darkgreen} $\pm$ 31.9}&{\color{darkgreen}43.0} {\color{darkgreen} $\pm$ 6.1}&-&-&-&-0.8 $\pm$ 0.3&-&{\color{blue}-50.3} {\color{blue} $\pm$ 16.4}\\  
999.59 $\pm$ 172.23 &$\rho$ & [kg m$^{-3}$] &8.5 $\pm$ 7.4&-&-&9.0 $\pm$ 7.8&13.7 $\pm$ 13.5&11.9 $\pm$ 5.0&-&-&-&-&-&-17.7 $\pm$ 12.7\\  
0.99 $\pm$ 0.17 &$\mu$ & [cP] &-17.1 $\pm$ 13.4&-&-&-18.2 $\pm$ 14.3&{\color{darkgreen}-28.3} {\color{darkgreen} $\pm$ 26.0}&-23.6 $\pm$ 10.0&-&-&-&-&-&{\color{blue}-63.0} {\color{blue} $\pm$ 20.6}\\  
1.99 $\pm$ 0.34 &$\gamma$ & [-] &-8.4 $\pm$ 7.4&-&-&-8.9 $\pm$ 8.0&-13.5 $\pm$ 13.5&-11.7 $\pm$ 4.8&-&-&-&-&-&{\color{darkgreen}-30.8} {\color{darkgreen} $\pm$ 10.1}\\  
136201.26 $\pm$ 22028.42 &$E_{max}$ & [Pa] &-14.2 $\pm$ 21.7&12.5 $\pm$ 5.8&17.9 $\pm$ 8.9&21.2 $\pm$ 19.6&{\color{darkgreen}30.8} {\color{darkgreen} $\pm$ 30.5}&-9.3 $\pm$ 6.7&-23.9 $\pm$ 10.2&-&-&-&-&-5.4 $\pm$ 3.9\\  
35051.20 $\pm$ 5948.93 &$E_{min}$ & [Pa] &{\color{red}-716.3} {\color{red} $\pm$ 643.2}&-19.2 $\pm$ 15.9&{\color{red}-118.9} {\color{red} $\pm$ 82.6}&{\color{red}-796.9} {\color{red} $\pm$ 692.8}&{\color{red}-1368.9} {\color{red} $\pm$ 1072.0}&{\color{darkgreen}-28.1} {\color{darkgreen} $\pm$ 48.5}&-0.5 $\pm$ 1.0&{\color{darkgreen}-38.7} {\color{darkgreen} $\pm$ 26.7}&1.3 $\pm$ 1.4&3.1 $\pm$ 1.0&-&-1.7 $\pm$ 13.4\\ 

\bottomrule[1.5pt] 
\end{tabular}} 
\caption{\scriptsize {\bf Sensitivity analysis of the one-dimensional lymph flow equations coupled to the EFMC model and valve dynamics. Positive pressure gradient}. Symbol $"-"$ indicates that the absolute value of the sensitivity index $\bar{S}_{i,j}$ was less than 0.5 and therefore parameter $x_i$ did not influence index $P_j$. Likewise, green-coloured parameters show an influence between 25 and 50, blue-coloured parameters show an influence between 50 and 100, and red-coloured parameters show an influence greater than 100. Results are shown using as mean $\pm$ SD.} \label{table:SA>0}

\end{sidewaystable}

{\em \bf Negative pressure gradient $\Delta P=P_{in}-P_{out}<0$}. 
Here we used $P_{in}$ = 3 cmH$_2$O and $P_{out}$ = 4 cmH$_2$O.  
The cross-sectional radius at equilibrium $r_0$ positively influence indexes SV, CPF, ESD, EDD, WSS and mean flow. Among the parameters of the EFMC model, $a_2$ is the most influential parameter. Indeed, it is the threshold to change the nature of the stationary point described in Section \ref {sec:equilibriumFHN} from stable to unstable. 
Parameter $R_I$ does not have a remarkable influence on the indexes. Likewise for parameters $k_{rel}$. The frequency, and thus FPF and CPF, are strongly influenced by $A_{Ca}/A_0$ and $n_{Ca}$. Indeed, the greater these parameters, the lower the frequency. Results for $n_{Ca}$ are in agreement with Frame \ref{fig:c}. Moreover, these parameters also influence the mean flow. Then parameters $k_{NO}$, $\tau_{NO}$, $n_{NO}$, which are related to WSS and passive flows, do not affect the lymphodynamical indexes since $\Delta P<0$. The parameter of the valve model $K_{vo}$, $K_{vc}$ and $L_{eff}$ do not affect the indexes. The fluid properties of lymph $\mu$ and $\gamma$ only affect the WSS, while the density $\rho$ has no effects on the indexes. The maximum and minimum Young modulus $E_{max}$ and $E_{min}$ affect the ESD and EDD, respectively, and also influence most of the parameters, such as the frequency, EF, SV, FPF, CPF and mean flow. To conclude, among the analysed parameters, nothing appears to influence the ESP. 

{\em \bf Positive pressure gradient $\Delta P=P_{in}-P_{out}>0$}.
Here we used $P_{in}$ = 4 cmH$_2$O and $P_{out}$ = 2 cmH$_2$O. 
The cross-sectional radius at equilibrium $r_0$ influences the CPF, SV, WSS, mean flow, ESD, EDD and WSS.
The most influential parameter of the EFMC model is $a_2$, followed by $a_3$ and $b_1$, which influences the frequency, CPF, FPF and WSS. Radius $R_I$ and parameter $k_{rel}$ do not affect the lymphodynamical indexes. As before, an increase of parameters $n_{Ca}$ and $A_{Ca}/A_0$ influence the WSS and leads to a negative chronotropic effect, which in turn influences the FPF and the CPF. An increase of parameter $k_{NO}$ decreases the frequency, indeed the greater this parameter, the greater the influence of the contraction inhibition given by the WSS. On the contrary, the increase of parameter $\tau_{NO}$ increases the frequency. Results for $k_{NO}$ and $\tau_{NO}$ are in agreement with Frames \ref{fig:f} and \ref{fig:e}, respectively, of Fig. \ref{fig:EFMC}. 
To a lesser degree, parameters $L_{eff}$, $\rho$, $\mu$ and $\gamma$ influence the frequency, FPF, CPF, WSS and the mean flow. Among these parameters, the greatest effect is seen on the mean flow given by the dynamic viscosity $\mu$. Between the minimum and maximum Young's modulus, the most influential one is $E_{min}$, as it affects the frequency, SV, FPF and CPF. It also affects WSS and EDD, though to a lesser extent.

\subsection{A quantitative study on the lymphodynamical effect of stenotic and regurgitant lymphatic valves}

\begin{table*}[t!] 
\centering 
\resizebox{\linewidth}{!}{ 
\begin{tabular}{c | c c c c | c c c c | c } 
\toprule[1.5pt] 
 & \multicolumn{4}{ c |}{Stenotic valve} & \multicolumn{4}{ c |}{Regurgitant valve} & Healthy valves  \\ 
& \multicolumn{2}{ c }{Left} & \multicolumn{2}{ c |}{Right}  & \multicolumn{2}{ c }{Left}  & \multicolumn{2}{ c |}{Right}  &  \\  \hline 
\\[-1em] 
& $M^L_{st}=$ 0.5 & $M^L_{st}$= 0.1 & $M^R_{st}=$ 0.5 & $M^R_{st}=$ 0.1 & $M^L_{rg}=$ 0.1 & $M^L_{rg}=$ 0.8 & $M^R_{rg}=$ 0.1 & $M^R_{rg}=$ 0.8 & $M_{st}=1$ $M_{rg}=0$\\  \hline 
\\[-1em] 
Frequency [min$^{-1}$] & 7.51 $\approx$&{\color{darkgreen}5.87 $\downarrow$}&7.54 $\approx$&7.54 $\approx$&7.54 $\approx$&7.54 $\approx$&{\color{blue}17.19 $\uparrow$}&{\color{blue}17.54 $\uparrow$}&7.54 \\ 
SW [nL cmH$_2$O]  & 160.74 $\approx$&{\color{darkgreen}174.27 $\uparrow$}&167.63 $\approx$&{\color{darkgreen}122.43 $\downarrow$}&162.01 $\approx$&{\color{blue}54.15 $\downarrow$}&126.10 $\approx$&{\color{blue}24.65 $\downarrow$}&148.10 \\ 
EF [-] & 0.66 $\approx$&0.65 $\approx$&0.65 $\approx$&{\color{blue}0.32 $\downarrow$}&{\color{darkgreen}0.77 $\uparrow$}&{\color{darkgreen}0.81 $\uparrow$}&0.70 $\approx$&0.69 $\approx$&0.66 \\ 
SV [nL] & 46.81 $\approx$&45.30 $\approx$&45.58 $\approx$&{\color{blue}22.07 $\downarrow$}&{\color{darkgreen}53.78 $\uparrow$}&{\color{darkgreen}56.66 $\uparrow$}&{\color{darkgreen}54.97 $\uparrow$}&{\color{darkgreen}54.76 $\uparrow$}&45.71 \\ 
FPF [min$^{-1}$] & 4.97 $\approx$&{\color{darkgreen}3.83 $\downarrow$}&4.94 $\approx$&{\color{blue}2.39 $\downarrow$}&{\color{darkgreen}5.83 $\uparrow$}&{\color{darkgreen}6.14 $\uparrow$}&{\color{blue}11.98 $\uparrow$}&{\color{blue}12.19 $\uparrow$}&4.95 \\ 
CPF [$\mu$L h$^{-1}$] & 21.09 $\approx$&{\color{darkgreen}15.97 $\downarrow$}&20.62 $\approx$&{\color{blue}9.99 $\downarrow$}&{\color{darkgreen}24.34 $\uparrow$}&{\color{darkgreen}25.64 $\uparrow$}&{\color{blue}56.70 $\uparrow$}&{\color{blue}57.65 $\uparrow$}&20.68 \\ 
Mean Flow [$\mu$L h$^{-1}$] & 19.09 $\approx$&{\color{darkgreen}12.90 $\downarrow$}&18.04 $\approx$&{\color{darkgreen}9.45 $\downarrow$}&{\color{darkgreen}11.38 $\downarrow$}&{\color{blue}1.42 $\downarrow$}&{\color{blue}4.36 $\downarrow$}&{\color{blue}-0.55 $\downarrow$}&17.70 \\ 
CPFI [-] & 1.10 $\approx$&1.24 $\approx$&1.14 $\approx$&1.06 $\approx$&{\color{blue}2.14 $\uparrow$}&{\color{red}18.06 $\uparrow$}&{\color{red}13.01 $\uparrow$}&{\color{red}103.90 $\uparrow$}&1.17 \\ 
WSS [mdyne cm$^{-2}$] & -75.35 $\approx$&{\color{darkgreen}-51.90 $\uparrow$}&-72.99 $\approx$&{\color{blue}-24.75 $\uparrow$}&{\color{blue}-35.04 $\uparrow$}&{\color{blue}-20.86 $\uparrow$}&{\color{blue}-19.41 $\uparrow$}&{\color{blue}1.23 $\uparrow$}&-70.81 \\ 
Peak Velocity [mm s$^{-1}$] & 6.56 $\approx$&6.88 $\approx$&{\color{darkgreen}5.74 $\downarrow$}&{\color{blue}0.94 $\downarrow$}&{\color{darkgreen}5.29 $\downarrow$}&{\color{darkgreen}9.49 $\uparrow$}&7.39 $\approx$&6.66 $\approx$&6.77 \\ 
ESD [$\mu$m] & 142.36 $\approx$&143.04 $\approx$&142.86 $\approx$&{\color{darkgreen}200.90 $\uparrow$}&{\color{darkgreen}115.98 $\downarrow$}&{\color{darkgreen}104.91 $\downarrow$}&142.40 $\approx$&142.84 $\approx$&142.45 \\ 
EDD [$\mu$m] & 244.95 $\approx$&242.72 $\approx$&243.09 $\approx$&243.10 $\approx$&243.10 $\approx$&243.11 $\approx$&258.72 $\approx$&258.63 $\approx$&243.08 \\ 
ESP [cmH$_2$O] & 6.56 $\approx$&6.41 $\approx$&6.95 $\approx$&{\color{darkgreen}9.48 $\uparrow$}&6.37 $\approx$&{\color{darkgreen}4.21 $\downarrow$}&6.43 $\approx$&6.52 $\approx$&6.40 \\ 
EDP [cmH$_2$O] & 3.00 $\approx$&2.93 $\approx$&3.00 $\approx$&3.00 $\approx$&3.00 $\approx$&3.00 $\approx$&{\color{blue}5.93 $\uparrow$}&{\color{blue}6.00 $\uparrow$}&3.00 \\ 

\bottomrule[1.5pt] 
\end{tabular}} \caption{\scriptsize {\bf Analysis of the effect of lymphatic valve deficits on lymphodynamical indexes}. Here we compare indexes for healthy and defective valves. In the first column, we show lymphodynamical indexes, while in the second and third we show results for the stenotic and regurgitant valve, respectively. In the last column results for healthy valves are shown. A green-coloured result indicates a normalised, percentage change in absolute compared to the healthy case value between 15 and 50 $\%$. Likewise, a blue-coloured result indicates a change between 50 and 200 $\%$ and a red-coloured result indicates a change above 200 $\%$. The arrows indicate a positive or a negative change. } \label{table:Valves} 
\end{table*} 

\begin{figure}[!h]
\begin{center}
\includegraphics[width=0.435\textwidth]{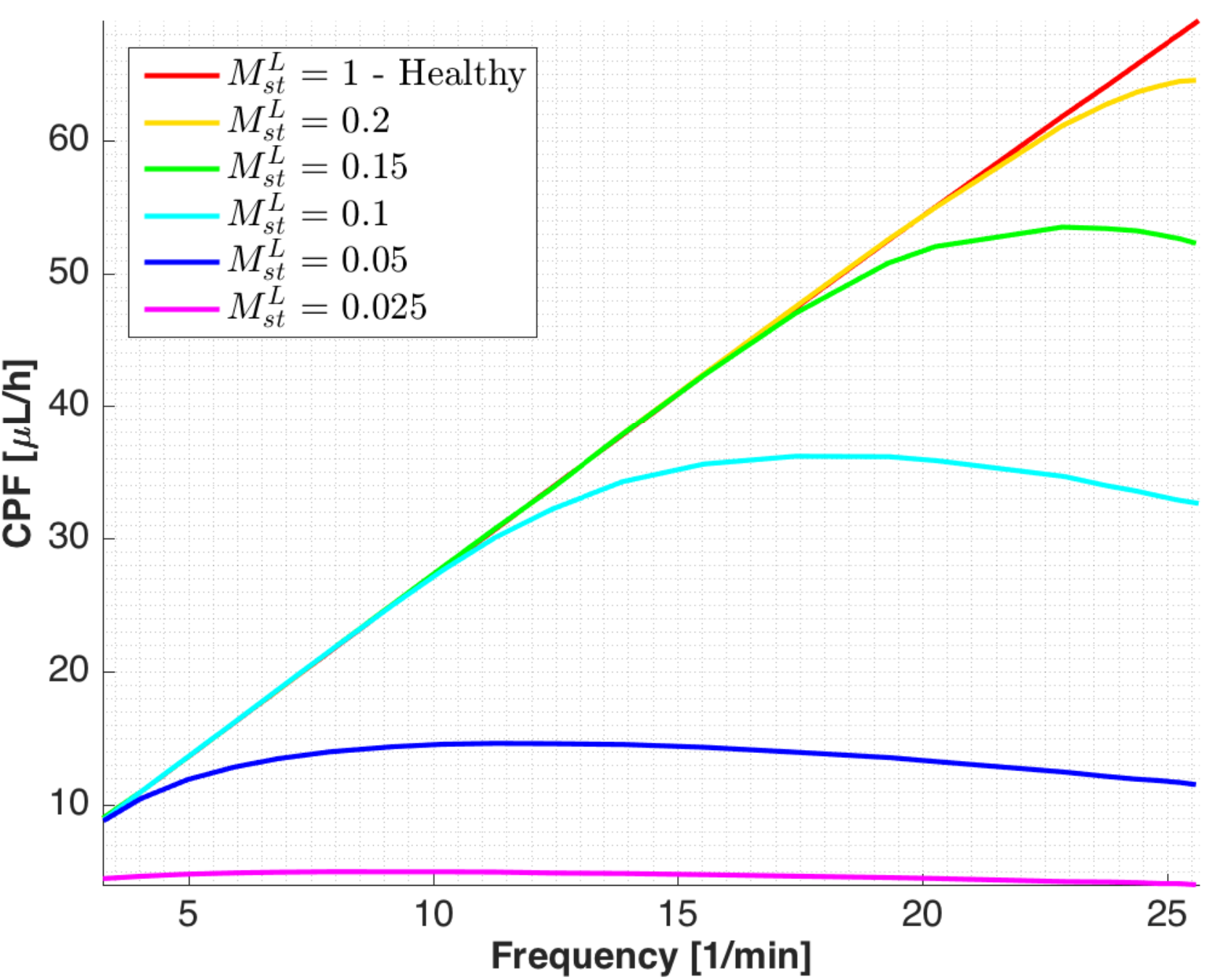}
\includegraphics[width=0.435\textwidth]{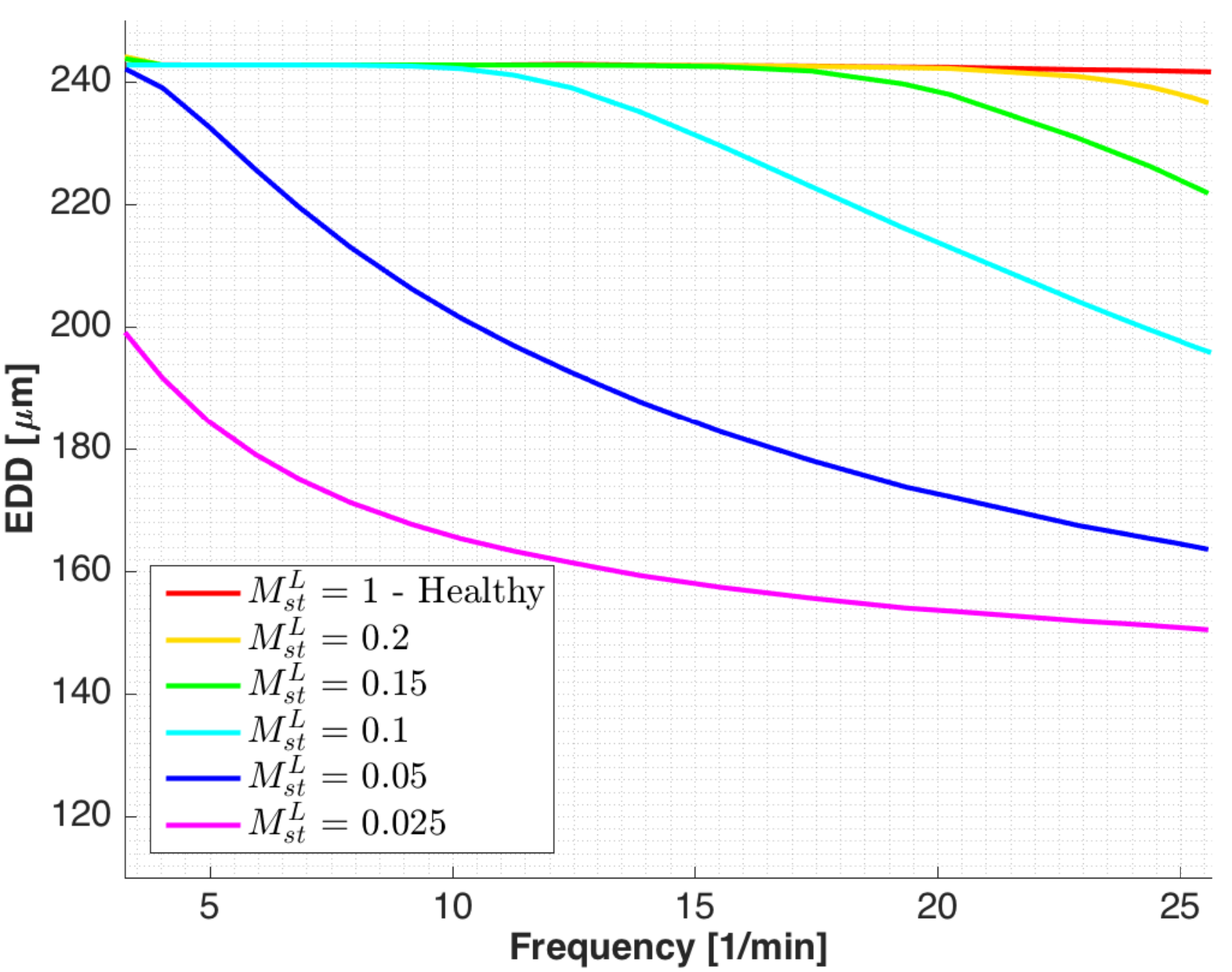}
\includegraphics[width=0.435\textwidth]{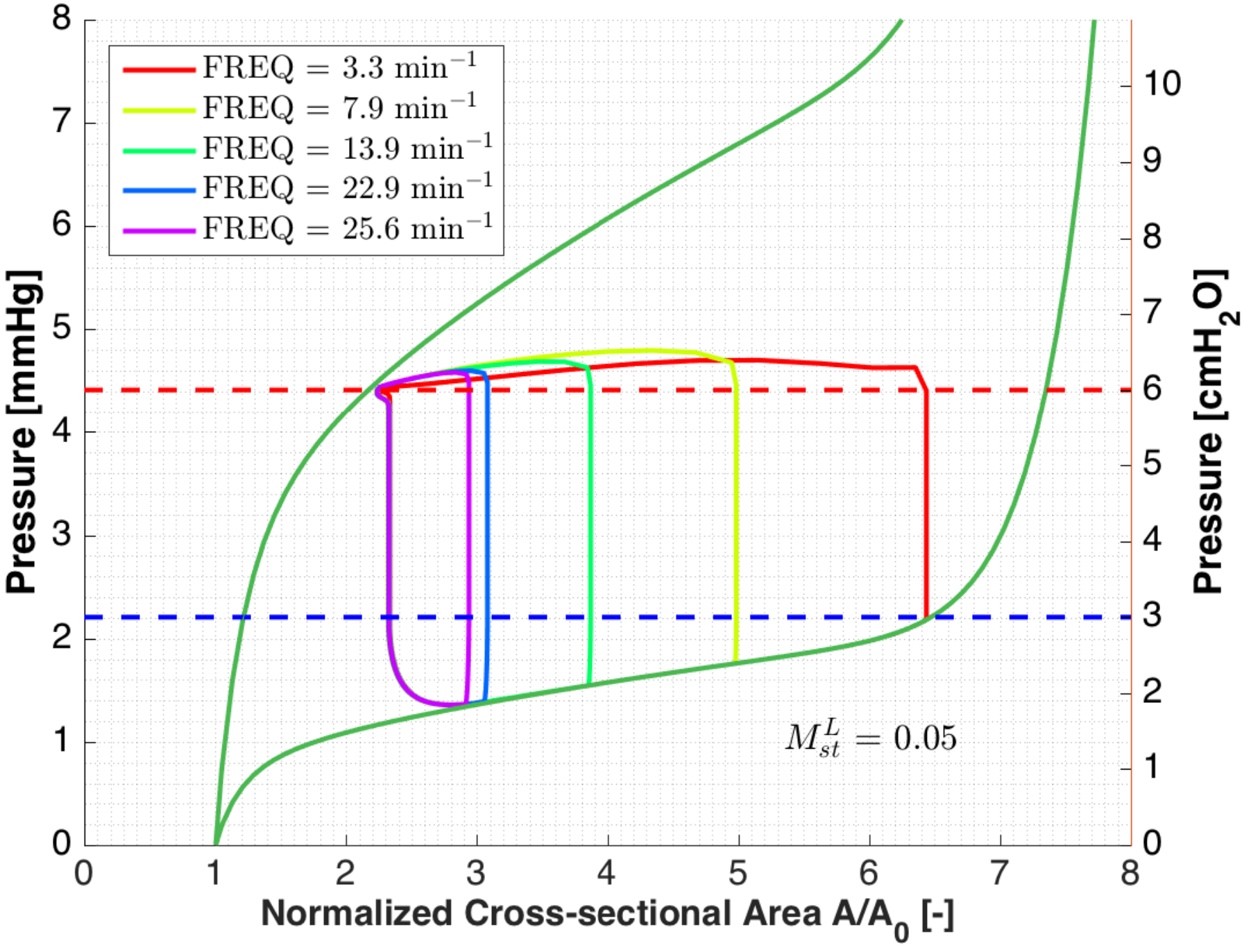}
\end{center}
\caption{\scriptsize {\bf High frequencies of contractions with a left stenotic valve diminish the CPF}. We simulated three lymphangions with two valves and two imposed boundary pressures at the ending interfaces $P_{in}$ = 3 cmH$_2$O and $P_{out}$ = 6 cmH$_2$O. The left valve is assumed stenotic. Results for the centre lymphangion are shown. Here we show frequency against CPF (top frame) and EDD (centre frame) for different severity of the stenosis, and pressure against normalised area (bottom frame) for different frequencies and with a severe stenosis. The higher the frequency, the greater the negative effect on the CPF caused by a severe left stenotic valve.} \label{fig:frequencyMstL}
\end{figure}

The mathematical model for collectors proposed in the present paper includes a well-established model for valves proposed by Mynard et al. \cite{Mynard:2012a}. It has already been used for the heart valve modelling \cite{Mynard:2015a}, as well as for the venous valves \cite{Toro:2015aa}. More interestingly, the model allows for a quantitative study of the effect of stenotic and regurgitant valves. For instance, the model was already used to study the impact on brain haemodynamics of bilateral stenotic and regurgitant valves of the internal jugular veins \cite{Toro:2015aa,Cristini:2014a}.

Here we aim to analyse the lymphodynamical effect of stenotic and regurgitant lymphatic valves. We modelled a single collector with three lymphangions and two valves. Here we used $M=20$ computational cells to discretize the one-dimensional lymph vessel. We simulated a collector cannulated at each end, that is,  we imposed a fixed pressure at the leftmost and rightmost interfaces of the collector, as described in \ref{sec:imposedPressure}. The imposed pressures were $P_{in}=3$ cmH$_2$O and $P_{out}=4$ cmH$_2O$. A healthy valve is characterized by parameters $M_{st}=1$ and $M_{rg}=0$. A stenotic valve can be simulated by reducing the maximum effective area $A_{eff,max}=M_{st}A_0$, and this can be done by decreasing parameter $M_{st}$. Thus, a severe stenosis can be modelled by setting $M_{st}\approx 0$. Viceversa, a regurgitant valve can be simulated by increasing the minimal effective area $A_{eff,min}=M_{rg}A_0$, and this can be done by increasing parameter $M_{rg}$. A severe regurgitant valve can be modelled by setting $M_{rg}\approx 1$. For further details, see the Eqs. reported in \ref{sec:valves}. See also for \cite{Mynard:2012a,Mynard:2015a,Toro:2015aa}.

We consider four possible situations: a left stenotic valve, a right stenotic valve, a left regurgitant valve and a right regurgitant valve. 
The numerical results of the centred lymphangion are shown in Fig. \ref{fig:incompetentValves}. From top to bottom we show in order: the PA loops, the time-varying lymphatic pressure, diameter and flow at the centre of the lymphangion. The boundary pressures $P_{in}$ and $P_{out}$ are shown in the PA loops and in time-varying pressure plots. The first two columns show results for the left and right stenotic valves, while the remaining two columns show results for the left and right regurgitant valves. Parameters $M_{st}^{L/R}$ and $M_{rg}^{L/R}$ were varied from 0 to 1. 

A {\bf left stenotic valve} diminishes the inflow from the upstream valve. This results in the following: the greater the severity of the left stenosis, the greater the time required to fill the centre lymphangion after a contraction. For the tests considered here, contractions occur at a frequency of $\approx$ 7 min$^{-1}$, which means approximately every $\approx$ 8.6 s. For a severe left stenosis ($M_{st}^L$ = 0.025), the time required to fill the lymphangion is $\approx$ 16.5 s. A severe reduction of the EDD can happen when the lymphangion does not have enough time to fill itself, and this may happen when the contraction period is less than 8.6 s. To verify this hypothesis, we performed additional simulations with a left stenotic valve, varying the frequency and for different severities of the stenosis. We set $f_{min}=f_{Ca}$ from 4 to 24 min$^{-1}$ and calculated the resulting frequency, CPF and EDD. The numerical results are shown in Fig. \ref{fig:frequencyMstL}. For a mild stenosis ($M_{st}^L>0.1$) and low frequencies, the CPF does not suffer any changes, but as soon as the frequency increases (e.g. above approximately 10 min$^{-1}$ for $M_{st}^L=0.1$), the CPF decreases depending on the severity of the stenoses. At the frequency of $f=25$ min$^{-1}$ and $M_{st}^L=0.1$, the CPF reduces from $\approx$ 67.4 $\mu$L h$^{-1}$ to $\approx$ 33.0 $\mu$L h$^{-1}$, that is it reduces of the 51 $\%$. For even more severe left stenoses ($M_{st}^L<0.1$), the CPF drastically decreases and the lymphangion becomes unable to push the lymph forward. At the frequency of $f=25$ min$^{-1}$ and $M_{st}^L=0.025$, the CPF reduces of the 93 $\%$, namely it reduces to $\approx$ 4.1 $\mu$L h$^{-1}$. This comes from a decrease of the EDD for high frequencies. The PA loops for different frequencies and a severe left stenosis are also shown. The higher the frequencies, the greater the shrinkage of the PA loops. Overall, a left stenosis causes a decrease of the CPF for high frequencies of contractions. 

A {\bf right stenotic valve} drastically increases the ESP and the ESD. This comes from the difficulties for the lymph to be pushed downstream through a stenotic passage. As a matter of fact, the outflow greatly decreases. 

From simulating right and left stenotic valves we speculate the following: a stenotic valve causes an increase of the systolic peaks in the upstream lymphangions and maintains almost unchanged the downstream pressures. Moreover, it causes a reduction of the CPF for high frequencies of contractions in the downstream lymphangions. It would be interesting to  experimentally quantify and compare the effect of stenotic lymphatic valves in regions where the frequency of contraction is high and low. 

A {\bf left regurgitant valve} has a great impact on the effective pump flow, namely the real amount of flow ejected from the lymphangion. As the severity of the left regurgitant valve increases, backflows increase and this can be seen with the negative values of flows. This means that during contractions, the lymph is ejected backwards into the upstream lymphangion, and not forward into the downstream one. Moreover, the ESD diameter decreases and for a severe left regurgitant valve the downstream valve stays closed most of the time (result not shown) insofar as the ESP decreases. 

To conclude, a {\bf right regurgitant valve} increases the leakage from the downstream valve, even for small values of $M_{st}^R$. This results in increasing the EDP from 3 cmH$_2$O to 6 cmH$_2$O, which corresponds to the downstream boundary pressures $P_{out}$. For severe right regurgitant valves, the upstream valve does not open during the lymphatic cycle (results not shown).

The effects of regurgitant and stenotic valves are summarised in the lymphodynamical indexes shown in Table \ref{table:Valves}. The table shows from the left to the right column: the indexes, the stenotic valve case, the regurgitant valve case and the healthy control. For each valve dysfunction, results for a moderate and a severe case are shown.
For a left stenotic case, there are almost no variations in any of the indexes. As discussed before, problems arise for high frequencies of contractions. For a right stenotic valve, EF, SV, FPF and CPF halve, while the ESP increases. The cases of regurgitant valves show interesting properties on some of the lymphodynamical indexes. As a matter of fact, EF, SV, FPF and CPF indexes do not indicate any reduction in the pumping performance. Instead, based on the results in Table \ref{table:Valves}, it seems that the pumping action of the lymphangion has undergone improvements with the incompetence of the valves. For instance, with a left regurgitant valve case, EF, SV, FPF and CPF increase. The same happens for the right regurgitant valve case for indexes SV, FPF and CPF, though FPF and CPF might have increased because the frequency was increased. This is obviously misleading: since a significant amount of lymph is flowing retrograde due to the deficit, the effective time-averaged flow is approximately zero. Thus, we would expect CPF to be zero. As it was pointed out by Scallan et al. \cite{Scallan:2016a}, indexes EF, SV, FPF and CPF are usually assumed to represent forward lymph flow and therefore they can give inaccurate results for dysfunctional valves. Interestingly, CPFI gives unrealistic results and shows the inaccuracy of CPF. This index should be bounded between 0 and 1, but since the mean flow reduces to zero while CPF even increases depending on the severity and on the side of the regurgitant valve, the CPFI results to be much greater than 1 (CPFI $\approx$ 103 for a severe right regurgitant valve). Thus, we speculate that index CPFI can indicate a dysfunctional regurgitant valve.

\afterpage{
 \begin{landscape}
\begin{figure}[h]
\centering
\subfloat{\label{fig:a}\includegraphics[width=0.28\textwidth,height=0.25\textheight,keepaspectratio]{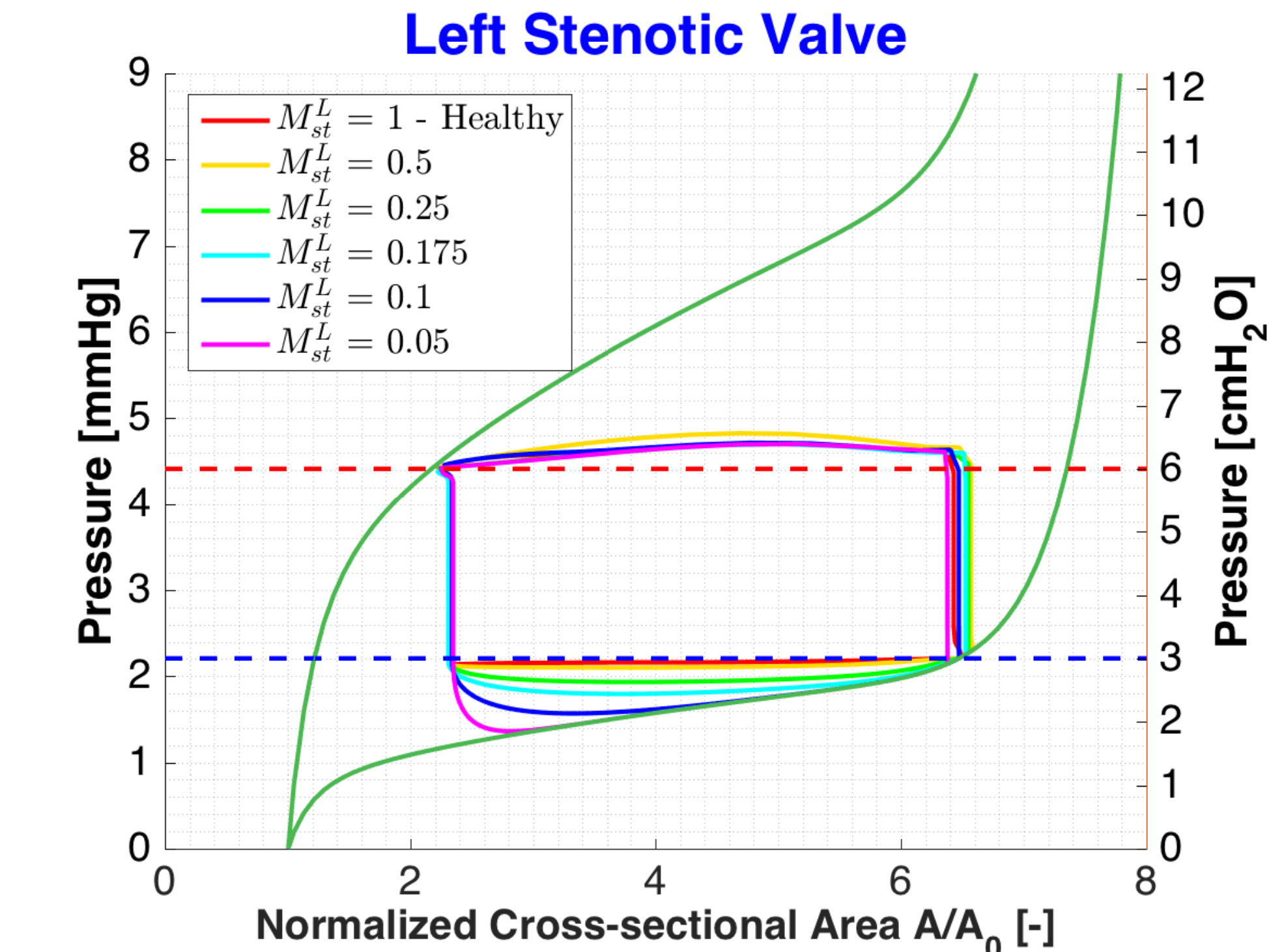}}
\subfloat{\label{fig:a}\includegraphics[width=0.28\textwidth,height=0.25\textheight,keepaspectratio]{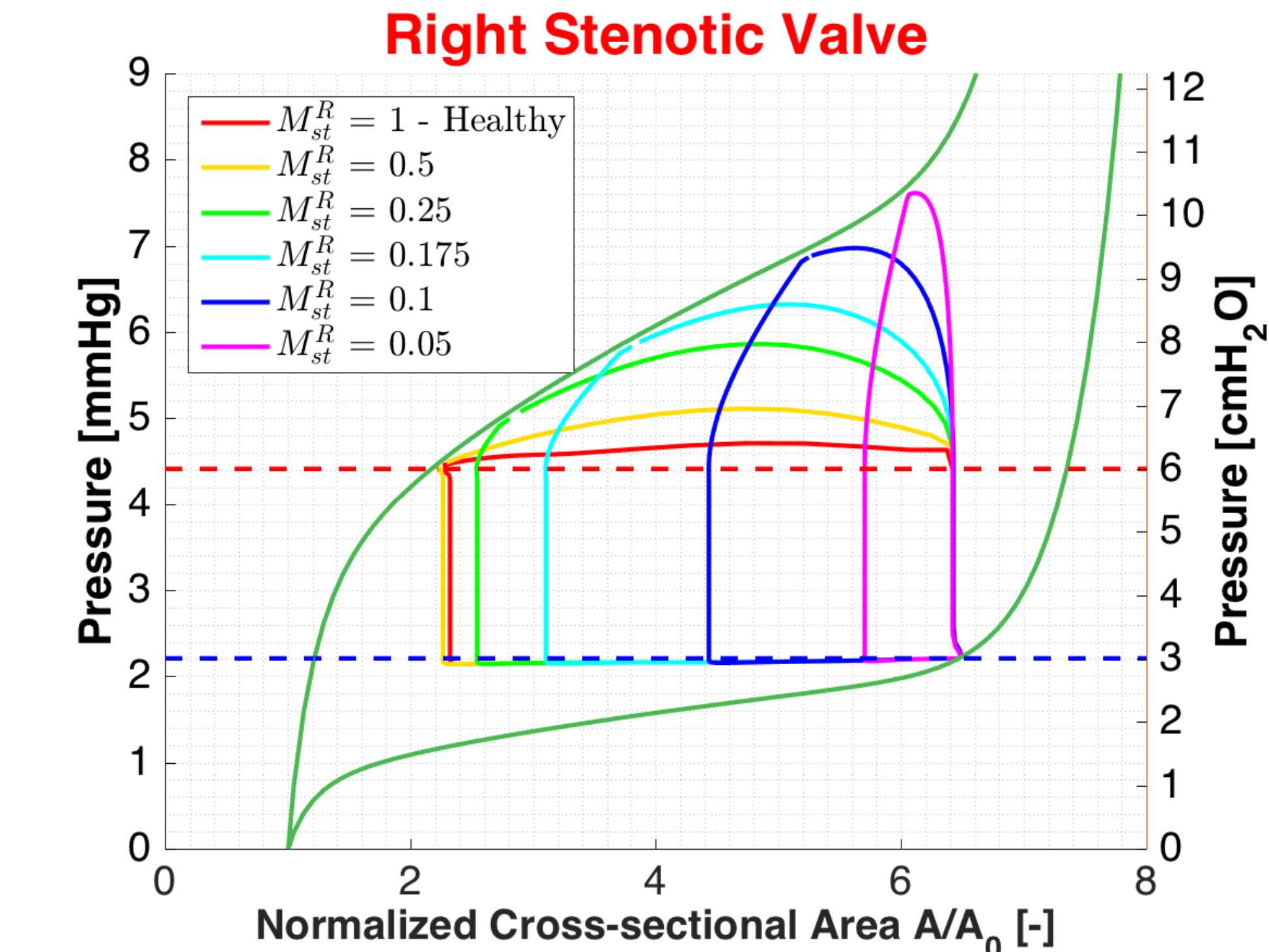}}
\subfloat{\label{fig:a}\includegraphics[width=0.28\textwidth,height=0.3\textheight,keepaspectratio]{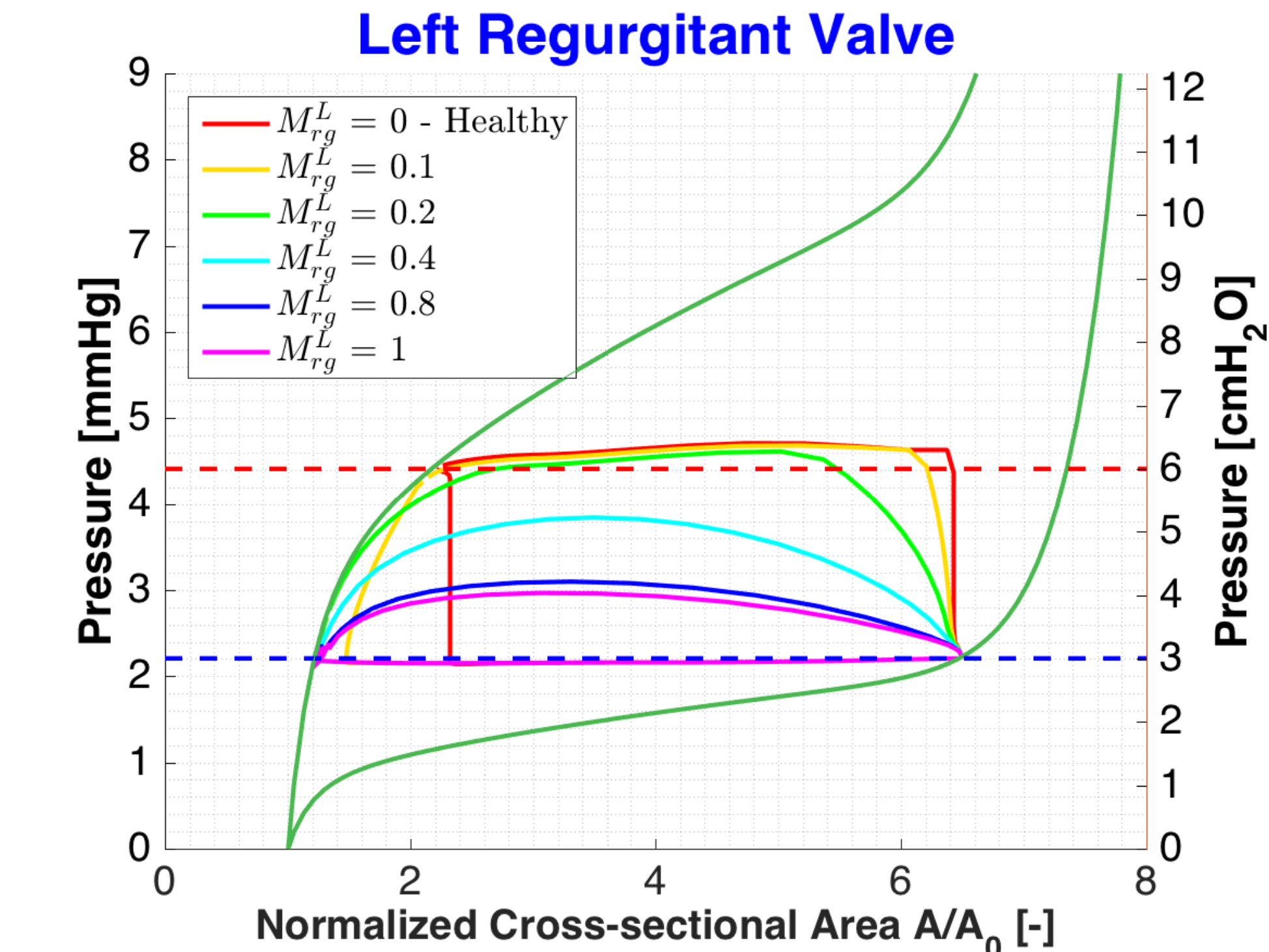}}
\subfloat{\label{fig:a}\includegraphics[width=0.28\textwidth,height=0.3\textheight,keepaspectratio]{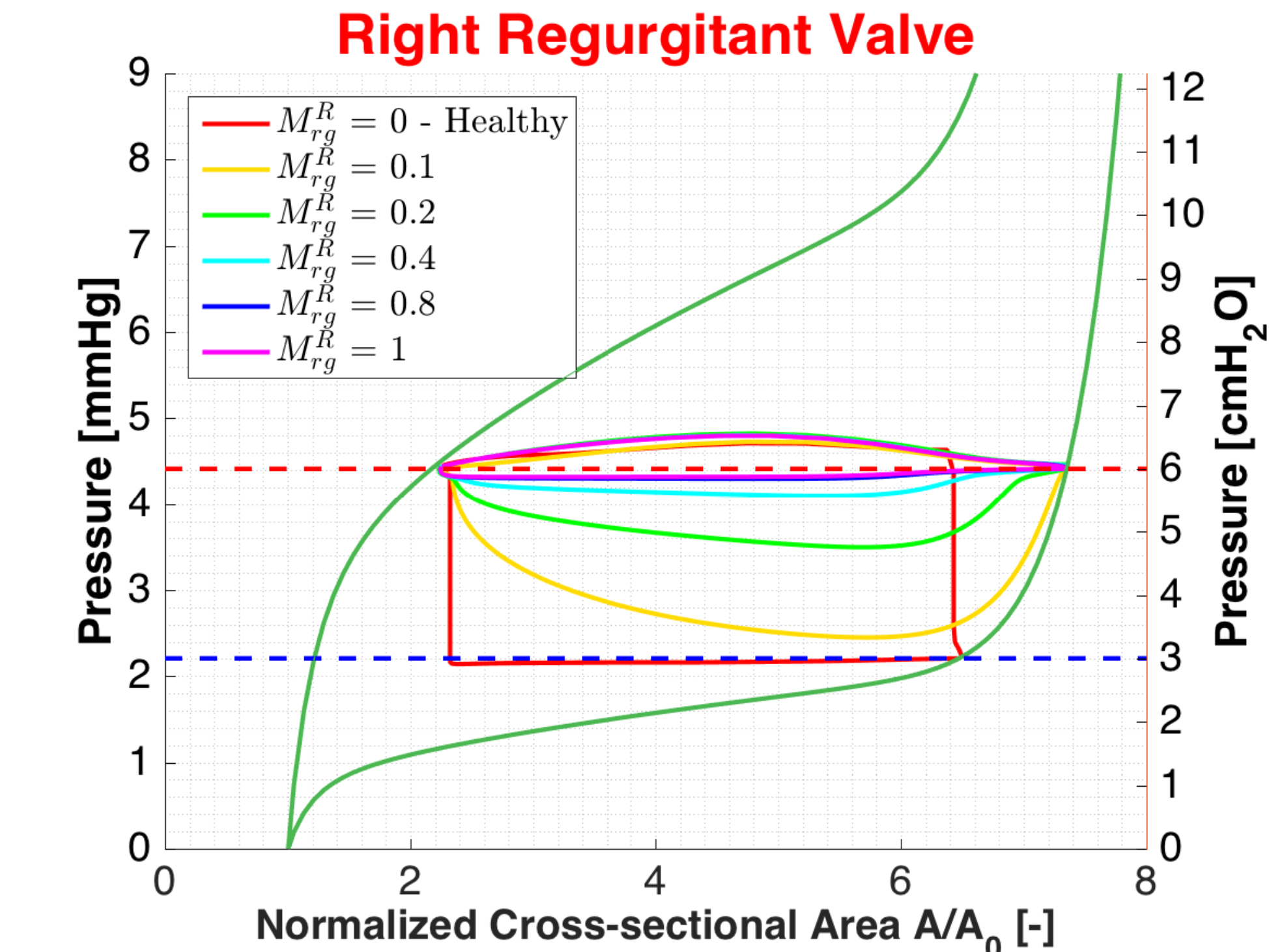}} \\
\subfloat{\label{fig:a}\includegraphics[width=0.28\textwidth,height=0.3\textheight,keepaspectratio]{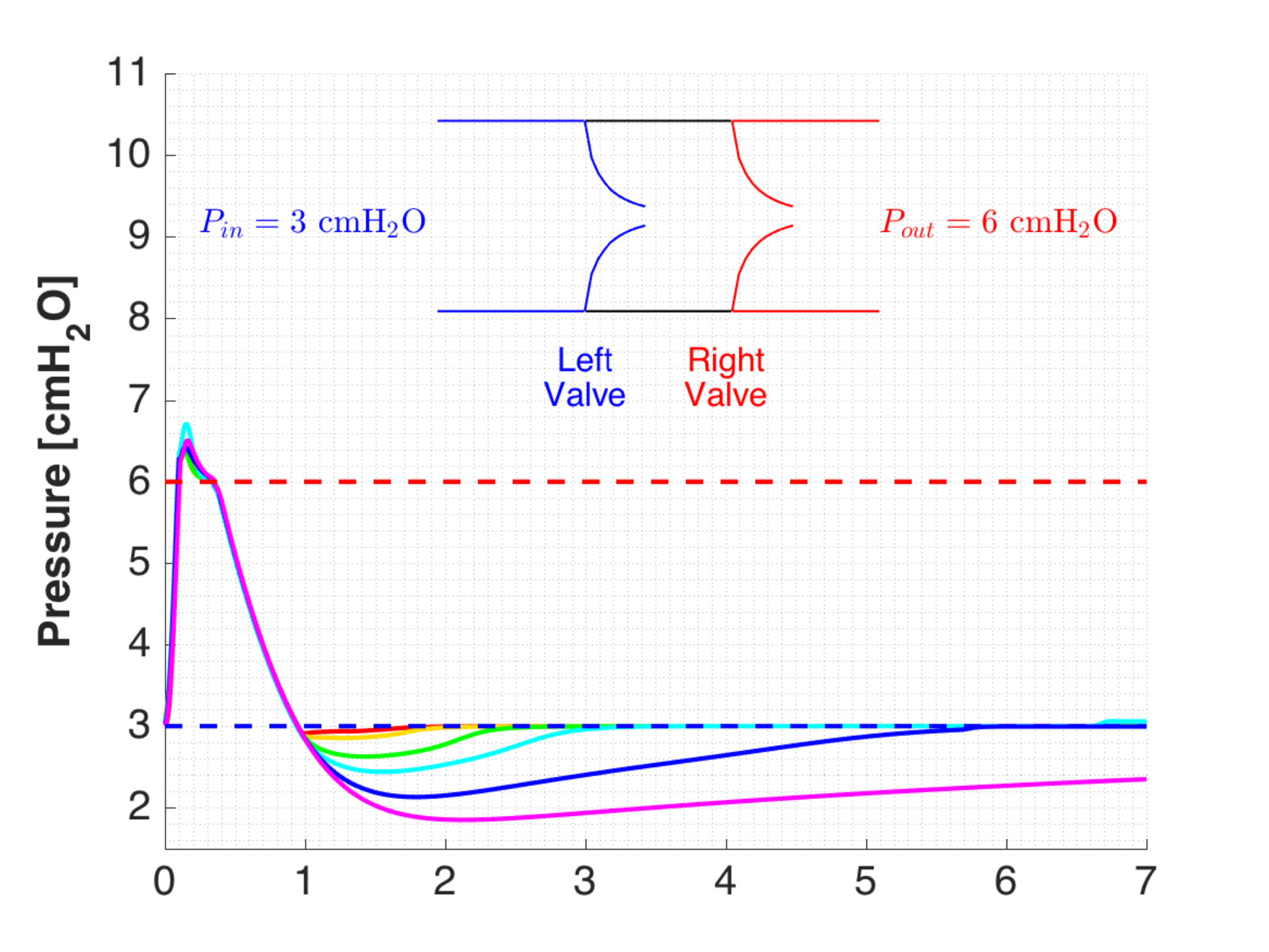}}
\subfloat{\label{fig:a}\includegraphics[width=0.28\textwidth,height=0.3\textheight,keepaspectratio]{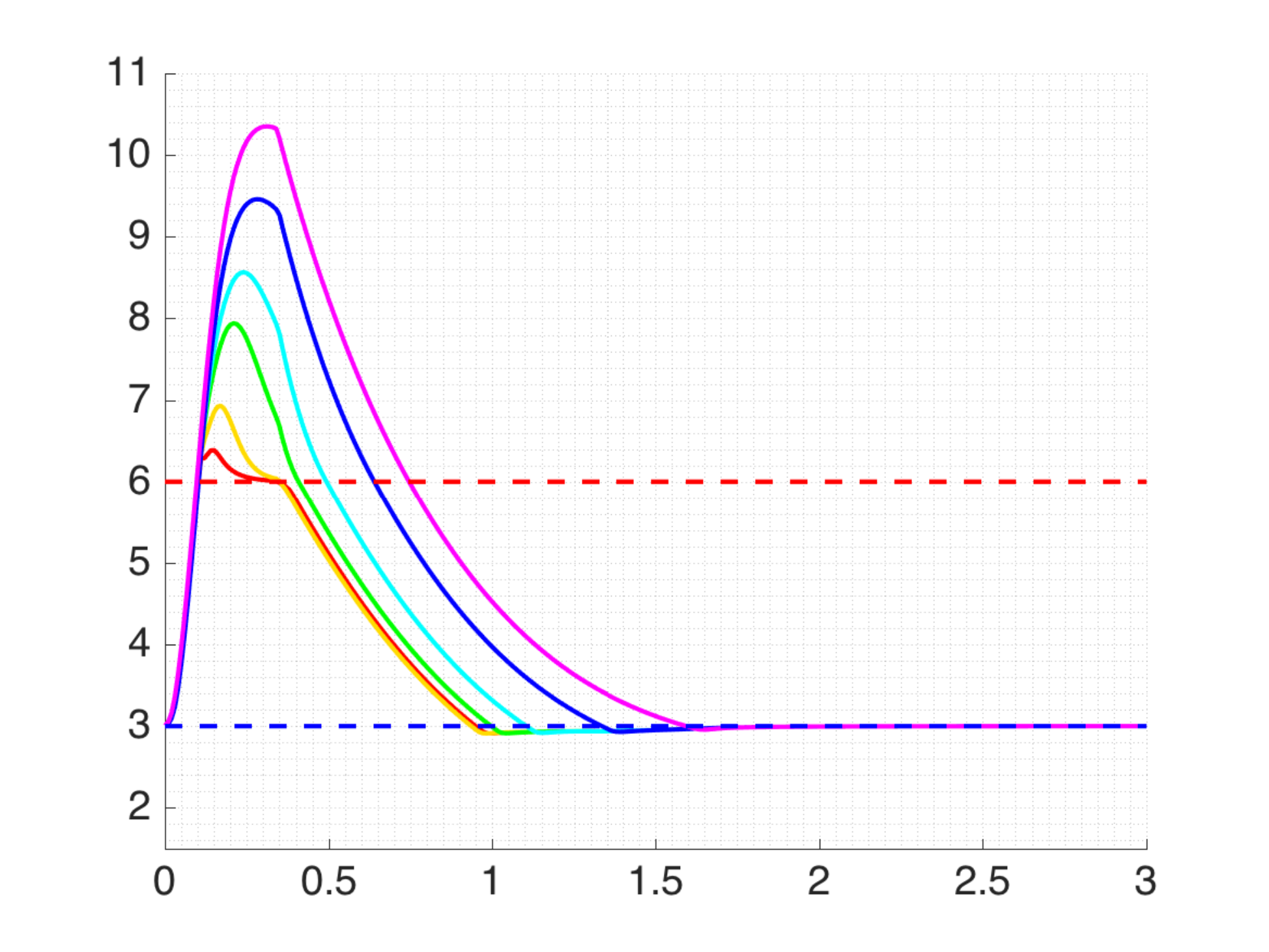}}
\subfloat{\label{fig:a}\includegraphics[width=0.28\textwidth,height=0.3\textheight,keepaspectratio]{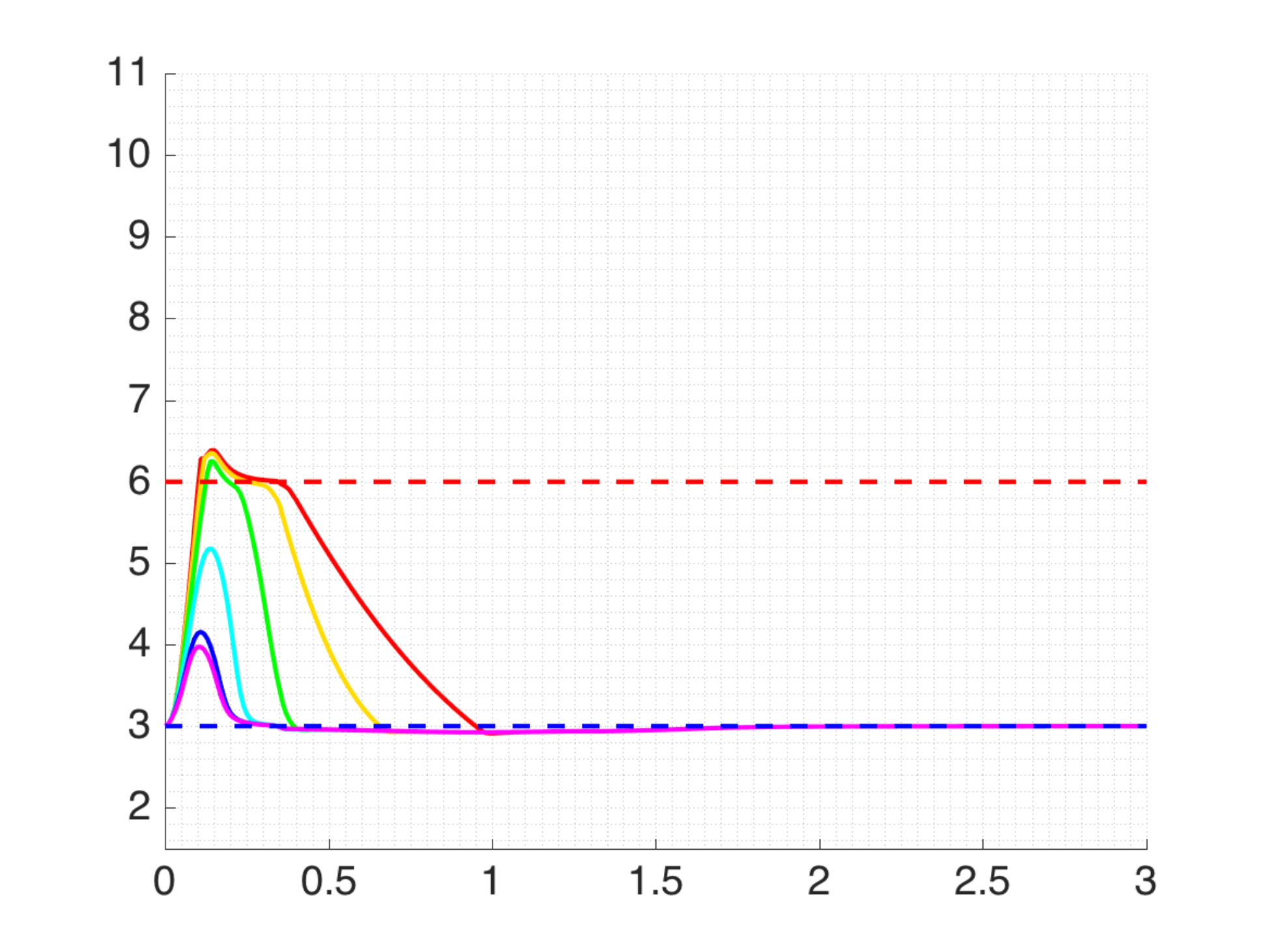}}
\subfloat{\label{fig:a}\includegraphics[width=0.28\textwidth,height=0.3\textheight,keepaspectratio]{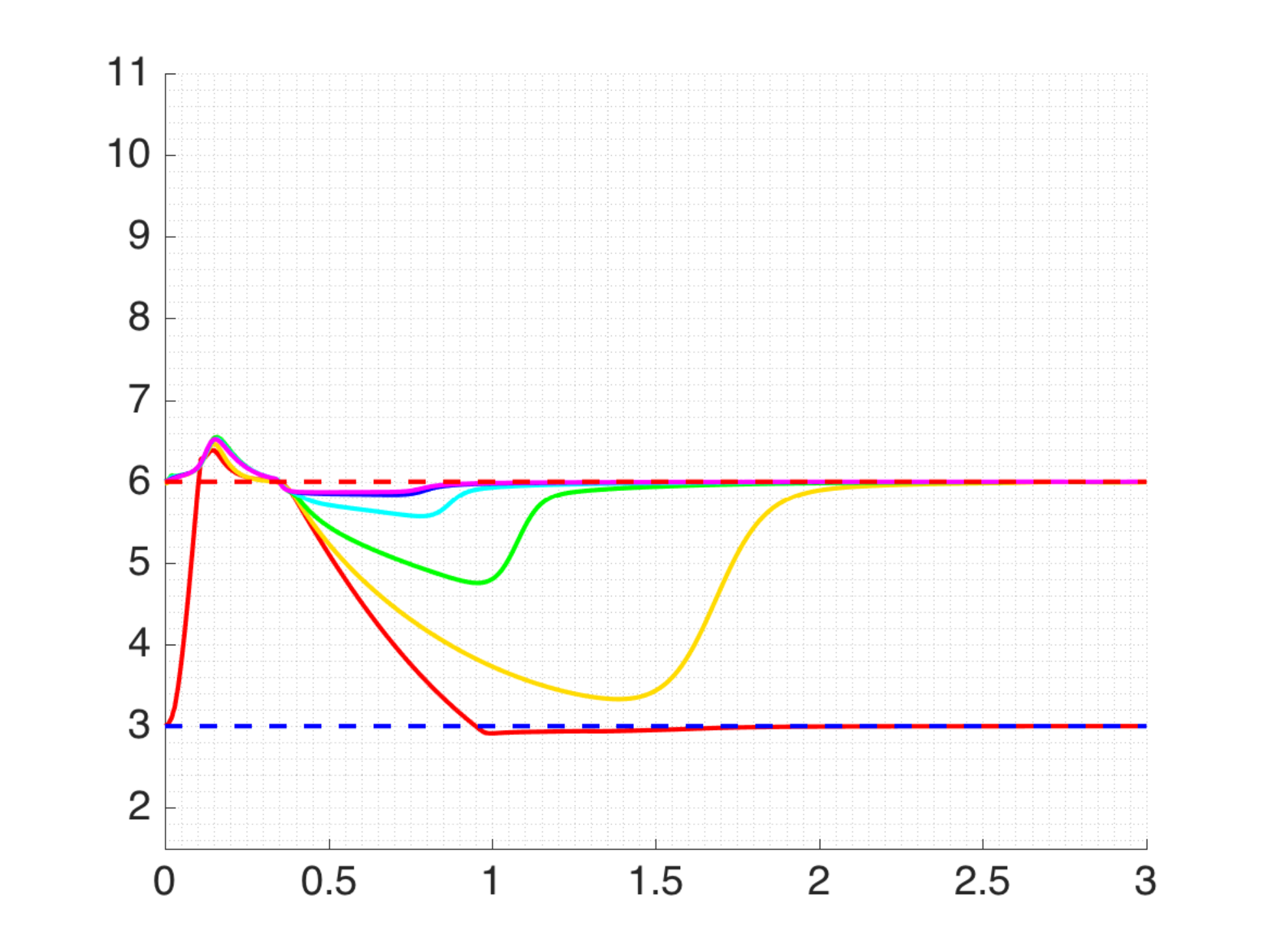}} \\
\subfloat{\label{fig:a}\includegraphics[width=0.28\textwidth,height=0.25\textheight,keepaspectratio]{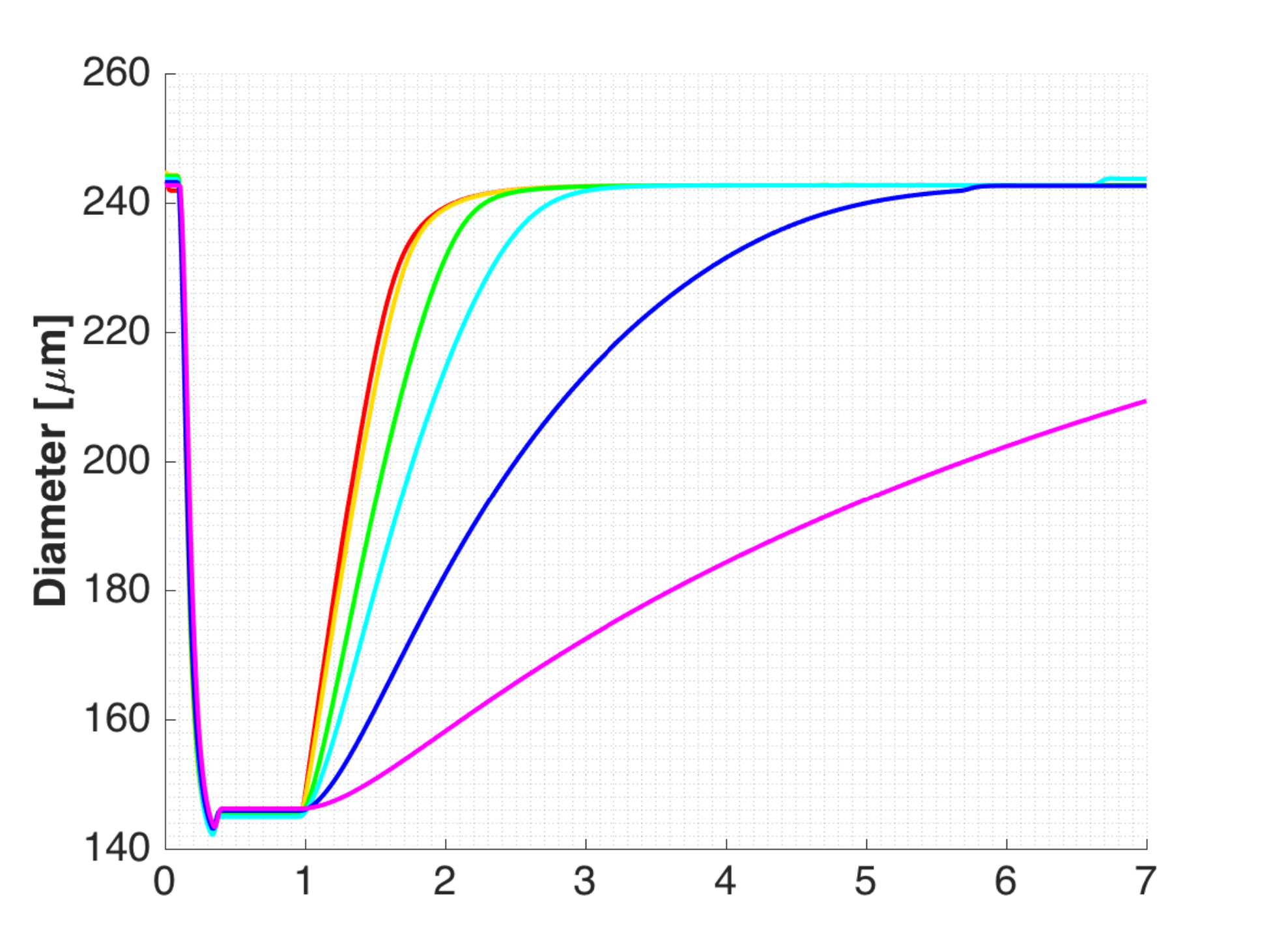}}
\subfloat{\label{fig:a}\includegraphics[width=0.28\textwidth,height=0.25\textheight,keepaspectratio]{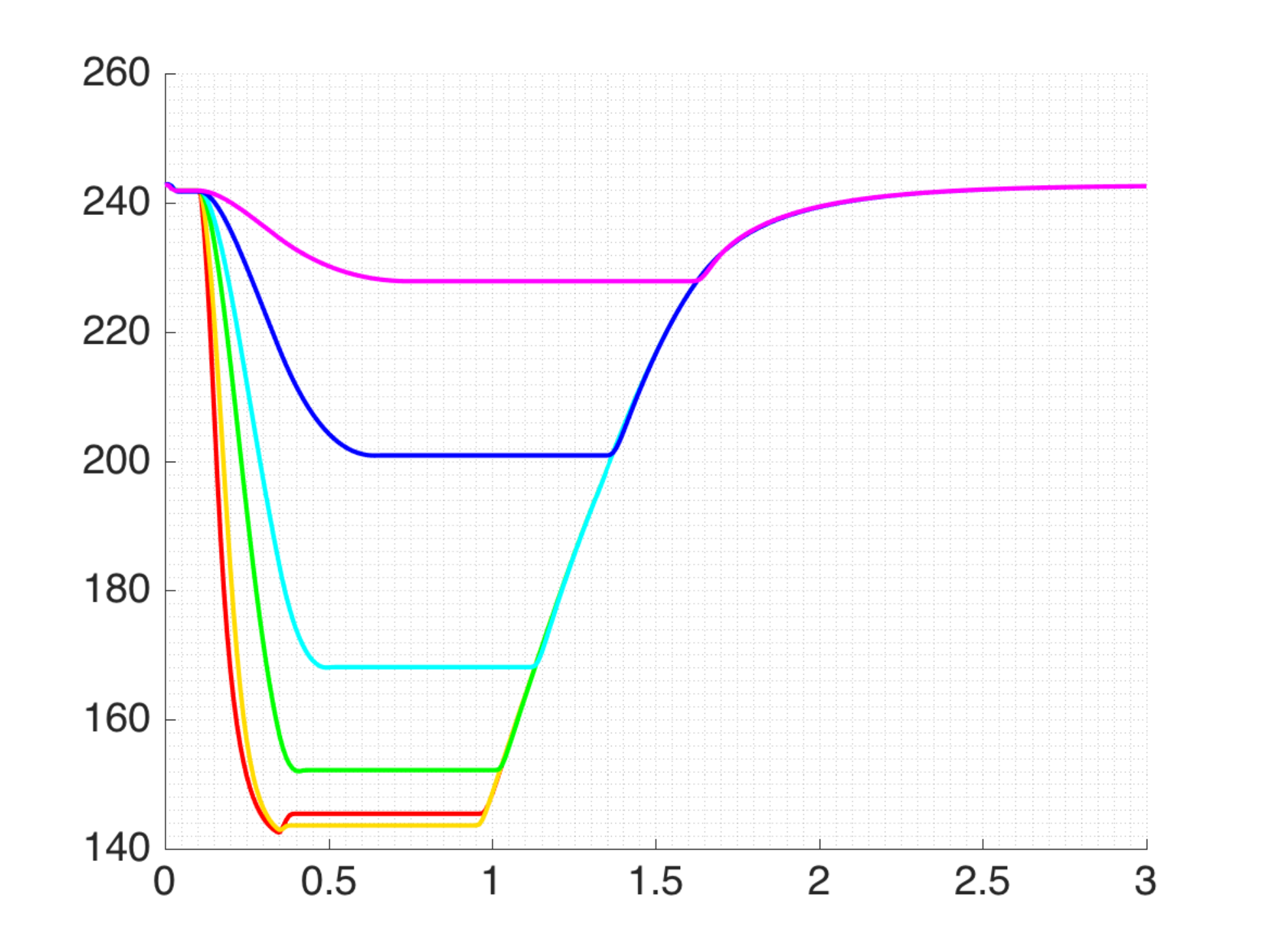}}
\subfloat{\label{fig:a}\includegraphics[width=0.28\textwidth,height=0.25\textheight,keepaspectratio]{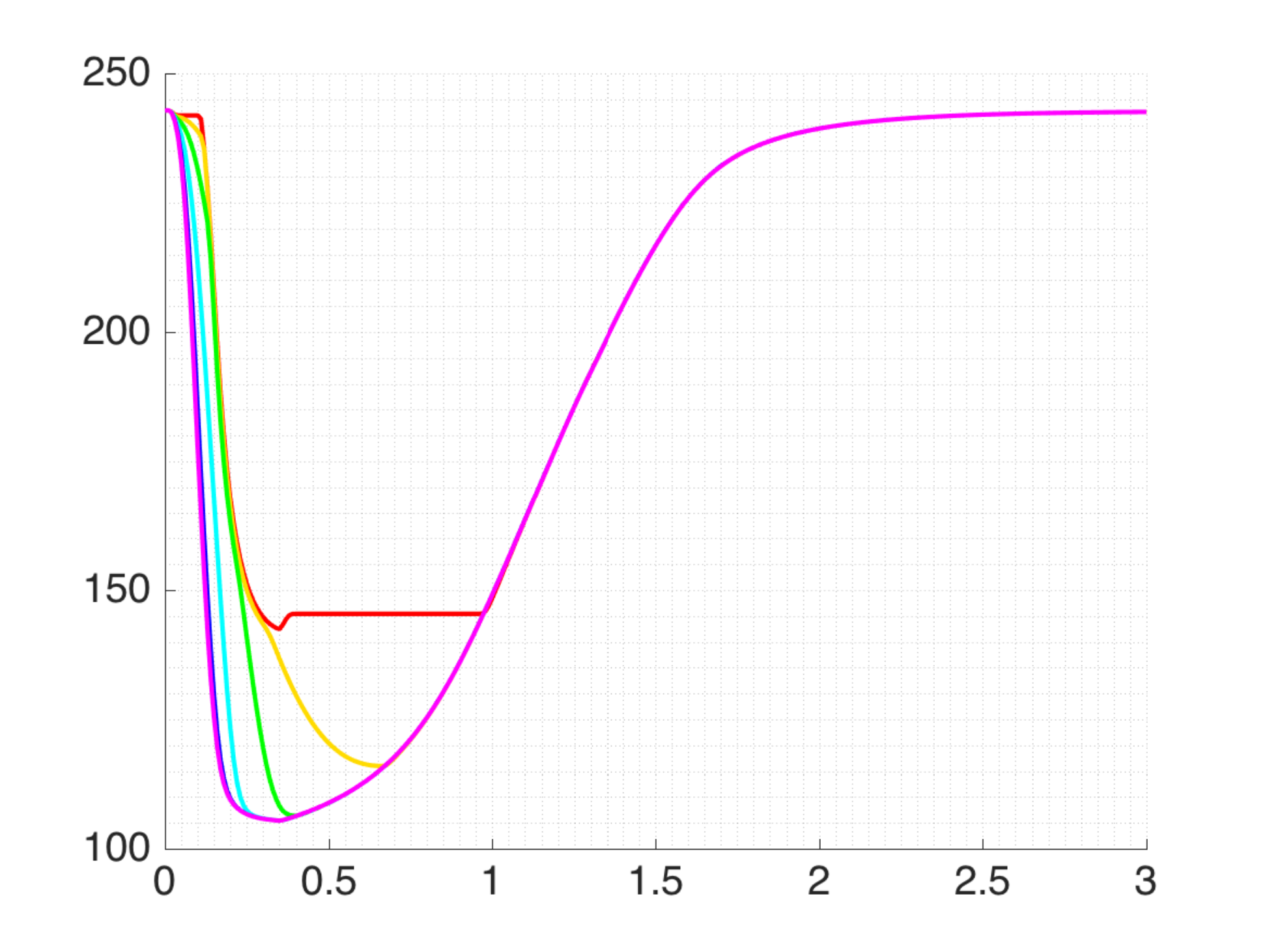}}
\subfloat{\label{fig:a}\includegraphics[width=0.28\textwidth,height=0.25\textheight,keepaspectratio]{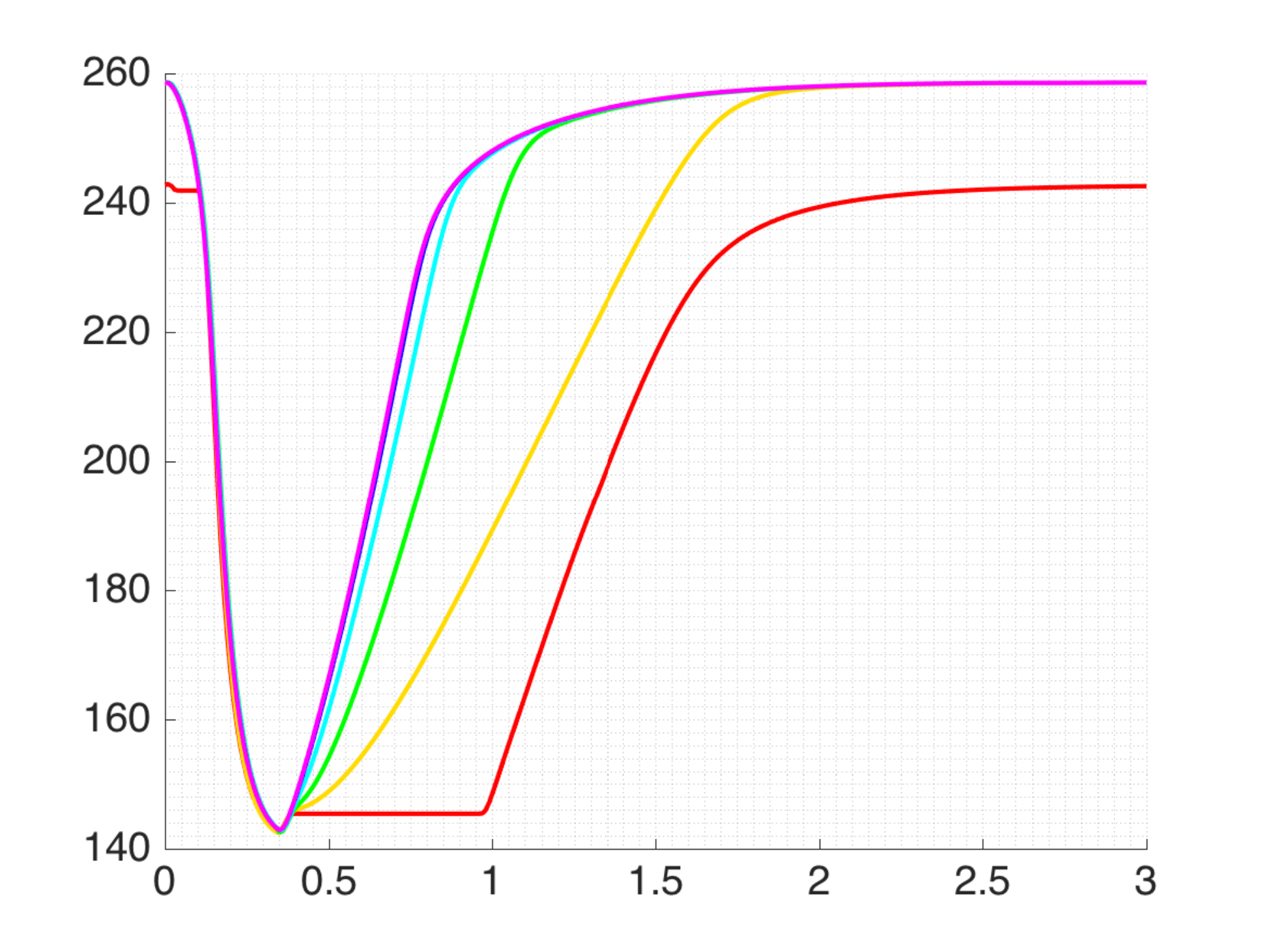}} \\
\subfloat{\label{fig:a}\includegraphics[width=0.28\textwidth,height=0.25\textheight,keepaspectratio]{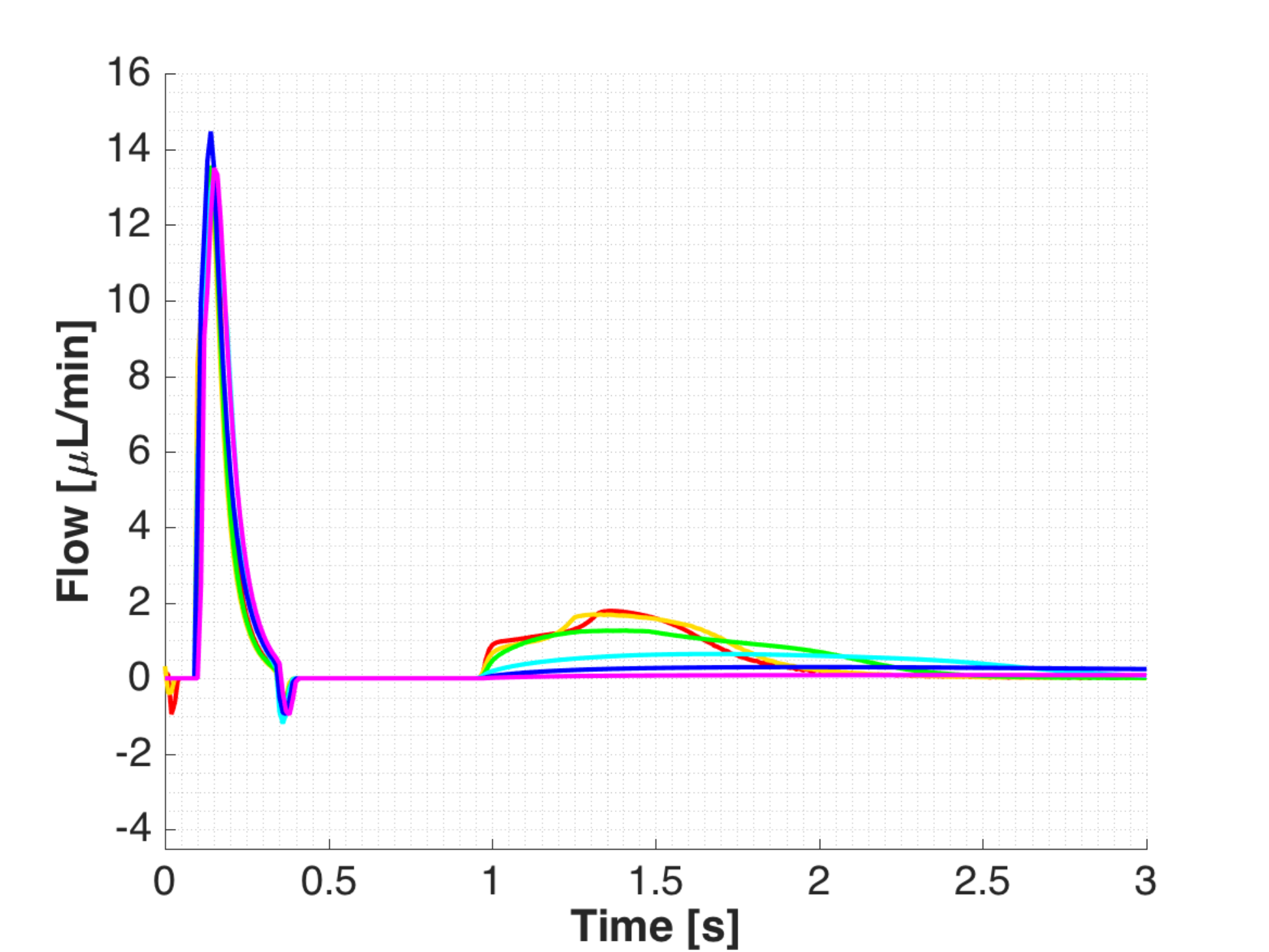}}
\subfloat{\label{fig:a}\includegraphics[width=0.28\textwidth,height=0.25\textheight,keepaspectratio]{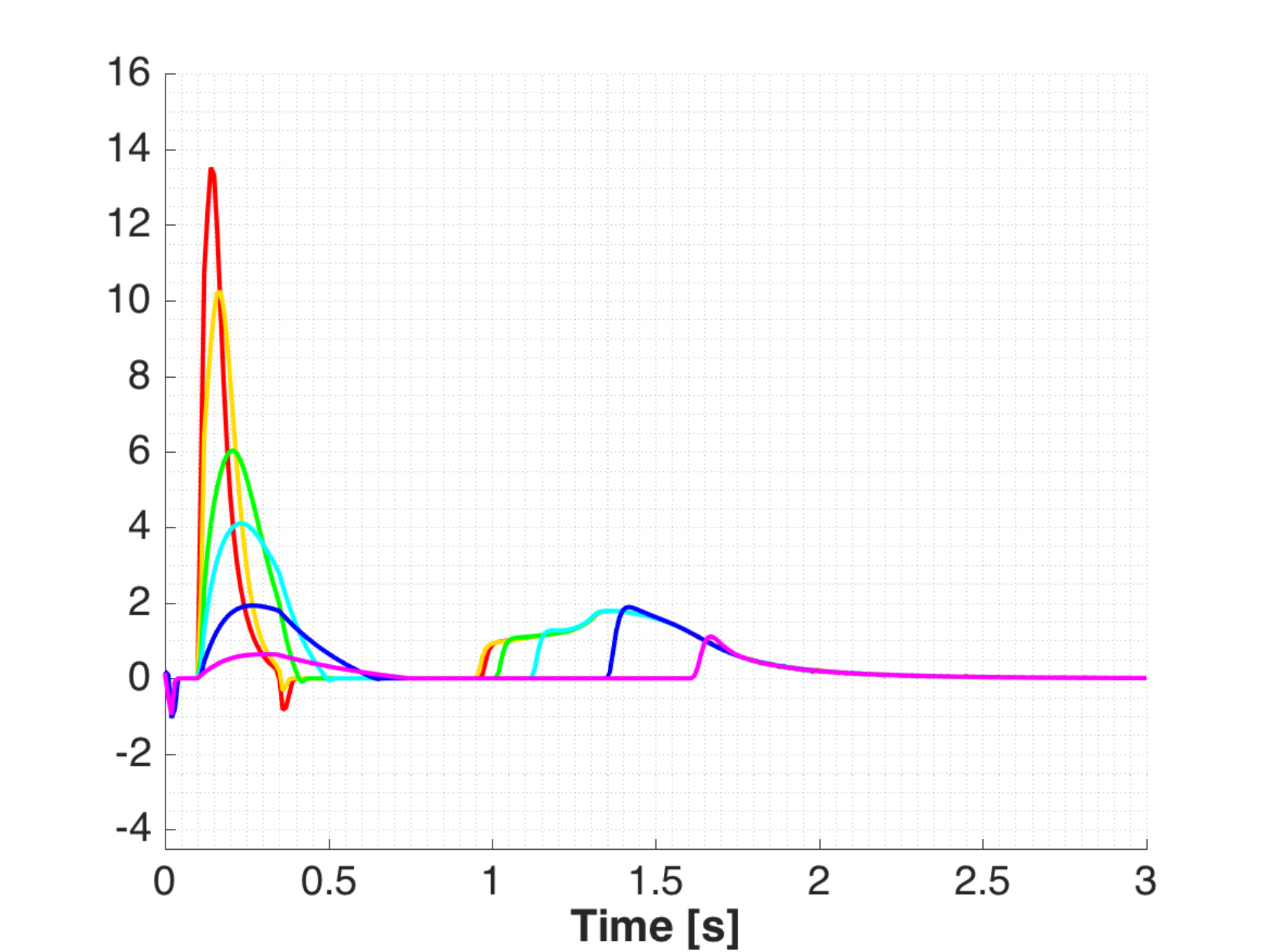}}
\subfloat{\label{fig:a}\includegraphics[width=0.28\textwidth,height=0.25\textheight,keepaspectratio]{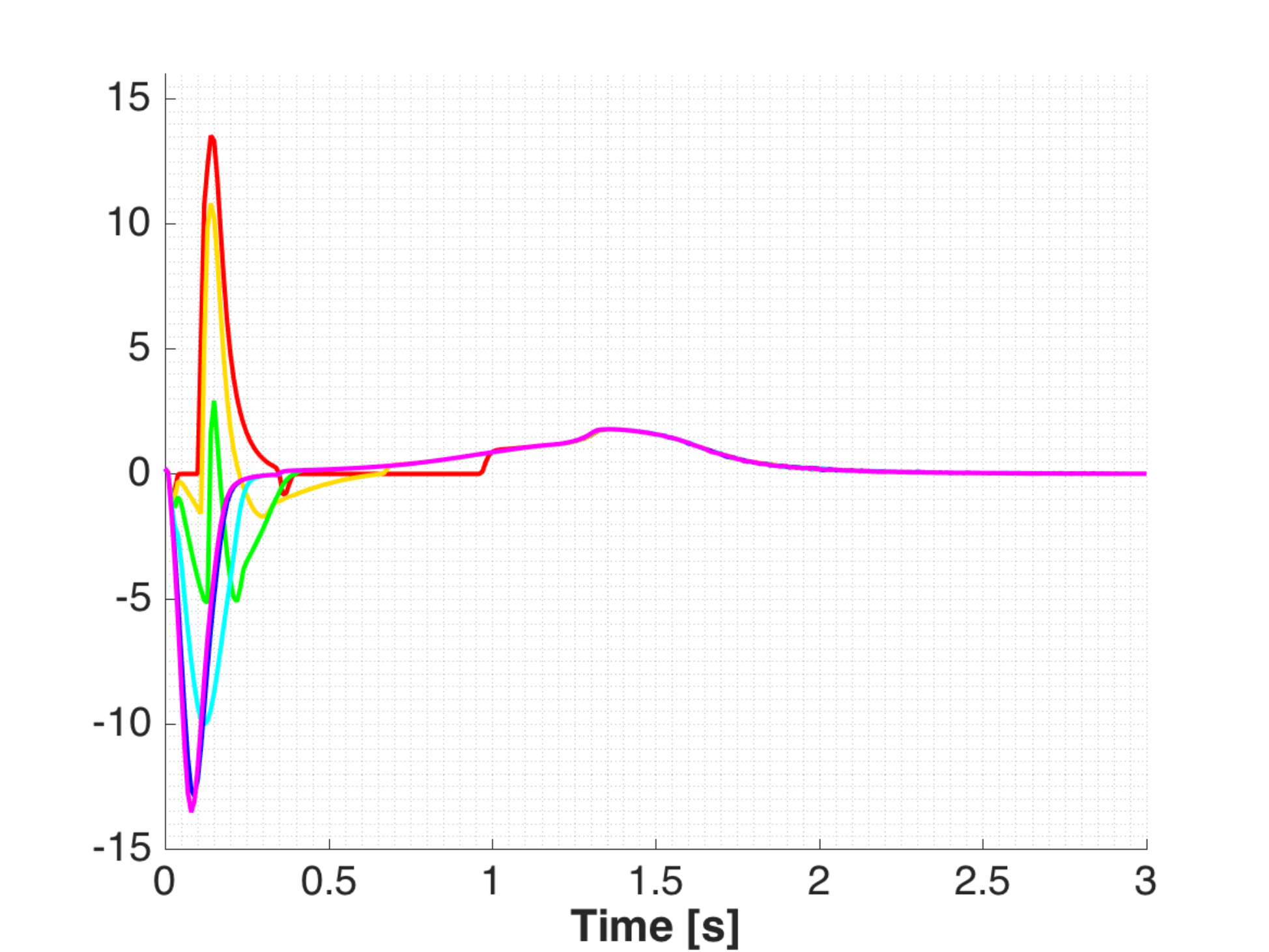}}
\subfloat{\label{fig:a}\includegraphics[width=0.28\textwidth,height=0.25\textheight,keepaspectratio]{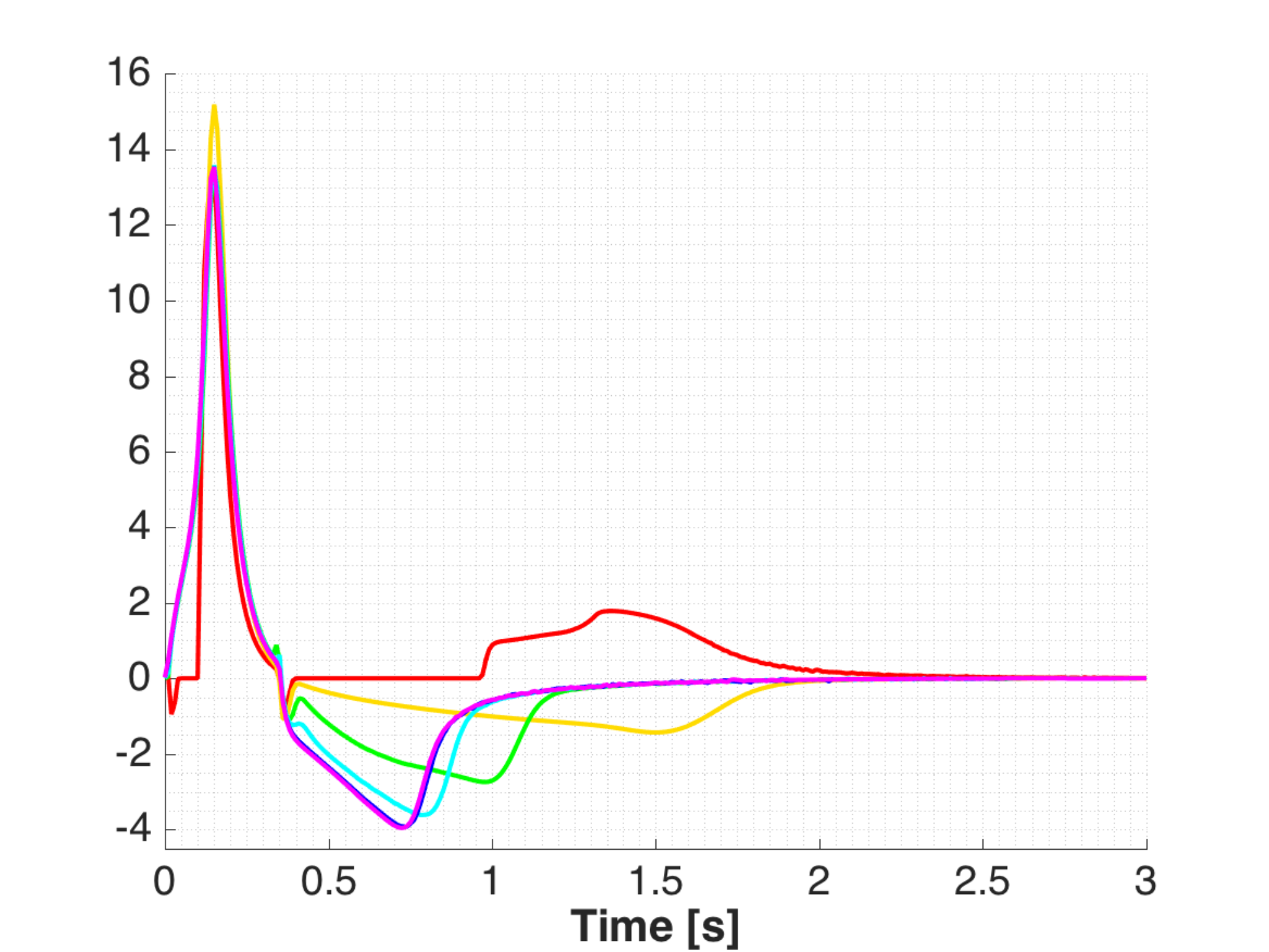}}
\caption{\scriptsize {\bf Effect of stenotic and regurgitant lymphatic valves}. From top to bottom we show in order: the PA loops, the lymphatic pressure, diameter and flow at the centre of the lymphangion. The boundary pressures $P_{in}$ and $P_{out}$ are shown in the PA loops and in pressure plots. The first two columns show results for the left and right stenotic valves, while the remaining two columns show results for the left and right regurgitant valves. Parameters $M_{st}^{L/R}$ and $M_{rg}^{L/R}$ were varied from 0 to 1.} \label{fig:incompetentValves}
\end{figure}
\end{landscape}
}

\section{Concluding Remarks} \label{Sec:conclusions}
In this paper, we have proposed a one-dimensional model for collecting lymphatics coupled to a novel Electro-Fluid-Mechanical Contraction (EFMC) model for dynamical contractions based on a modified FitzHugh-Nagumo model for action potentials and to a modern lumped model for valve dynamics. 
We performed an analysis of lymphodynamical indexes for a wide range of upstream and downstream pressure combinations. Numerical experiments showed that the contraction frequency resulting from the mathematical model strongly depends on the baseline frequency contraction, on the intramural pressure and on the wall-shear stress. 
We also performed numerical tests inspired by experiments where terminal lymphangions were cannulated, and the numerical results showed a good agreement with the experimental trend. The most influential model parameters were found by performing two sensitivity analyses for positive and negative pressure gradients.
We then quantified the effect of stenotic and regurgitant valves. A stenotic valve caused an increase of the systolic peaks in the upstream lymphangions and maintained almost unchanged the downstream pressures. Moreover, it caused a reduction of the CPF in the downstream lymphangions for high frequencies of contractions (up to 93$\%$ for a severe stenosis), while the CPF remained unaltered for low frequencies.
A regurgitant valve was unable to prevent backflows, and this resulted in zero net flows during the lymphatic cycle. The lymphodynamical indexes EF, SV, FPF and CPF were misleading for a regurgitant valve, insofar as from these indexes, it seemed that the pumping action of the lymphangion had undergone improvements with the incompetence of the valves. As a matter of fact, these indexes are usually assumed to represent forward lymph flow and therefore they can give inaccurate results for dysfunctional valves. To overcome this problem, we introduced the Calculated Pump Flow Index (CPFI), ideally bounded between 0 and 1, that indicates if the lymph is driven by contractions (CPFI = 1) or by a passive flow induced by a positive transaxial-pressure gradient (CPFI = 0). For a regurgitant valve, CPFI gave unrealistic results (CPFI $\approx$ 103 for a severe regurgitant valve). 

Several improvements can be made to the mathematical model. The EFMC model can be generalised to include a diffusion term. The resulting model would be able to simulate contraction waves travelling in the lymphatic wall and induced by pacemaker cells. In this way, gap-junctional communications would be included in the mathematical model. Then, the tube law can be improved in several ways. For instance, it can include a viscous term for the vessel wall. High-order methods can be used to improve the accuracy of the numerical solution, and these methods should properly couple the systems of ODEs with the one-dimensional model. 
Moreover, variable geometric parameters have been neglected, but they can be included in the simulation, as the mathematical formulation proposed in the present work allows for their presence. Finally, we also assumed an ex vivo setting, in which there is no interaction with the environment, such as skeleton muscle contraction or lymphatic contractions due to neuro-activities. 

The present EFMC model has been coupled to a one-dimensional model, but we speculate that it can be coupled to lumped parameter lymphatic models. Networks of collecting lymphatics can be simulated through the proposed model, but the lack of quantitative experimental measurements represents a great problem in validating the numerical results. We believe that the current mathematical model of collecting lymphatic can be coupled to multi-scale, closed-loop mathematical model of the cardiovascular system and can give quantitative lymphodynamical information in healthy and pathological cases.

\section*{Acknowledgment}
The authors gratefully acknowledge the suggestions given by Prof. Christian Vergara from the Department of Mathematics in Milan.
	
\section*{References}
\bibliographystyle{elsarticle-num.bst}      
\bibliography{MyBibFile.bib}

\end{document}